\documentclass[a4paper,10pt]{amsart}
\usepackage{amsmath,amsfonts,amssymb,amsthm, xypic,mathrsfs,stmaryrd,mathtools}
\usepackage{epsfig,psfrag}
\usepackage[all]{xy}
\usepackage{tikz-cd}
\usepackage{hyperref}

\def\1{\mathbf{1}}

\def\Id{\mathrm{Id}}  
\def\C{\mathbf{C}}    

\def\X{\mathbb{X}}    

\newtheorem{theo}{\bf{Theorem}}[section]
\newtheorem{lem}[theo]{Lemma}

\newtheorem{cor}[theo]{Corollary}
\newtheorem{prop}[theo]{Proposition}

\DeclareMathOperator{\GL}{GL}

\DeclareMathOperator{\Hom}{Hom}
\DeclareMathOperator{\PHom}{ParHom}
\DeclareMathOperator{\PCoh}{ParCoh}

\DeclareMathOperator{\PExt}{ParExt}
\DeclareMathOperator{\PPic}{ParPic}
\DeclareMathOperator{\PTor}{ParTor}
\DeclareMathOperator{\PBun}{ParBun}
\DeclareMathOperator{\PHiggs}{ParHiggs}

\DeclareMathOperator{\Par}{par \; deg}

\DeclareMathOperator{\Aut}{Aut}

\DeclareMathOperator{\Ext}{Ext}
\DeclareMathOperator{\Pic}{Pic}
\DeclareMathOperator{\Tor}{Tor}
\DeclareMathOperator{\Tr}{Tr}
\DeclareMathOperator{\Spec}{Spec}
\DeclareMathOperator{\Rep}{Rep}
\DeclareMathOperator{\Bun}{Bun}

\DeclareMathOperator{\Coh}{Coh}

\DeclareMathOperator{\Gal}{Gal}

\DeclareMathOperator{\rank}{rank}

\numberwithin{equation}{section}

\begin{document}

\title{Spherical Hall Algebras of weighted projective curves}

\author{Jyun-Ao Lin}

\begin{abstract}
In this article, we deal with the structure of the spherical Hall algebra $\mathbf{U}$ of coherent sheaves with parabolic structures on a smooth projective curve $X$ of abitraty genus $g$. We provide a shuffle-like presentation of the bundle part $\mathbf{U}^>$ and show the existence of generic spherical Hall algebra of genus $g$. We also prove the algebra $\mathbf{U}$ contains the characteristic functions on all the Harder-Narasimhan strata. These results together imply Schiffmann's theorem on the existence of Kac polynomials for parabolic vector bundles of fixed rank and multi-degree over $X$. On the other hand, the shuffle structure we obtain is new and we make links to the representations of quantum affine algebras of type $A$. 
\end{abstract}


\maketitle

\section{Introduction}

Hitchin introduced the Higgs bundles on Riemann surfaces of arbitrary genus $g$ as solutions of a two-dimensional reduction of the $4D$ Yang-Mills equations and computed the Betti numbers of the moduli spaces of stable Higgs bundles of rank $2$ and fixed degree via Morse theory \cite{Hit1,Hit2}. The Poincar\'e polynomials of these moduli spaces were then been computed in rank $3$ by Gothen \cite{Got}, in rank $4$ by Garcia-Prada, Heinloth and Schmitt \cite{GPHS} and finally by Schiffmann \cite{S5} in general cases via the numbers of (isomorphism classes of) geometrically indecomposable vector bundles on the smooth projective curves of genus $g$. The later's method was based on the combinatorial (i.e. shuffle) structure of the spherical Hall algebras of the category of coherent sheaves on curves. In order to extend Schiffmann's method to parabolic cases,  in this paper, we investigate the structure of the (spherical) Hall algebras of parabolic coherent sheaves. 

On the other hand, the Hall algebra of a curve has its own role playing in the theory of quantum groups and automorphic forms: initiated studied by Kapranov in \cite{Kap}, the spaces of unramified automorphic forms of $GL_r$ can be regarded as the subalgebra of the Hall algebra associated to vector bundles. In the case of $\mathbb{P}^1$, Kapranov translated the functional equations satisfied by Eisentein series associated to pairs of cuspidal forms into commutation relations between the corresponding generators in the Hall algebra (see also \cite{BK}). And such commutation relations bear a resemblance with those appearing in Drinfeld's loop-liked realization of quantum affine algebras \cite{Dr87}. Moreover, the classical Hecke operator associated to a point on automorphic forms is translated by the Hall multiplication with the characteristic function of the corresponding torsion sheaf. Such an apporach was later generalized to curves of arbitrary genus by Burban, Fratila, Kapranov, Schiffmann and Vasserot \cite{BS2, Fra, SV1, SV2}...etc and to weighted projective lines (i.e. parabolic coherent sheaves on $\mathbb{P}^1$) by Burban and Schiffmann \cite{BS, S2, S3}. Our study on the parabolic cases can be thought as a completion of the investigation of Hall algebras of hereditary categories. 

The strategy of studying the parabolic cases is somehow straightforward but many difficulties arise and also some new phenomenas appear. In order to apply the general construction of Hall algebra $\mathbf{H}$ to (quasi-)parabolic vector bundles over a smooth projective curve $X$ of arbitrary genus $g$, we have to borrow the notion of (quasi-)parabolic coherent sheaves, introduced by Heinloth in \cite{Hei}, that form a hereditary abelian category. As in the case of Hall algebras of quivers and curves, we are mainly interested in the natural subalgebra $\mathbf{U}$ (called spherical Hall algebra) generated by the characteristic functions associated to parabolic torsion sheaves and parabolic line bundles of fixed classes in the numerical Grothendieck group. We first prove the following important but rather expected result
\begin{theo}(Theorem \autoref{T:HN})
$\mathbf{U}$ contains all the characteristic functions of Harder-Narasimhan strata. 
\end{theo}

The main difficulty of studying the structure of $\mathbf{U}$ is that the subalgebra corresponding to parabolic torsion sheaves is not commutative anymore but several copies of the (positive halves of) quantum affine algebras of type $A$. An explicit Hecke action of the quantum affine algebras on the space of functions on parabolic vector bundles is necessary in order to compute the images of constant terms maps (i.e. the composition of the natural projection on the parabolic vector bundles part and the iterated comultiplication). Let $\mathbf{U}^>[1]$ denote the space of functions associated to rank $1$ parabolic bundles. We thus obtain
\begin{prop}(Corollary \ref{C5:indecomp})
$\mathbf{U}^>[1]$ carries a natural structure of $\mathbf{U}_v(\widehat{\mathfrak{sl}}_{w_1}) \otimes \cdots \otimes \mathbf{U}_v(\widehat{\mathfrak{sl}}_{w_N})$-modules. Moreover, $\mathbf{U}^>[1]$ is isomorphic to the affine extremal weight module $V(\varpi_{1,1}) \otimes \cdots \otimes V(\varpi_{N,1})$ of extremal weight $\varpi_{1,1}, \dots, \varpi_{N,1}$. 
\end{prop}
Here the $w_i$'s are certain positive integers defining the parabolic structure on vector bundles over $X$ and $\varpi_{i,1}$ are respectively the first fundamental weight of the quantum affine algebras $\mathbf{U}_v(\widehat{\mathfrak{sl}}_{w_i})$. It is known that such an affine extremal weight module is irreducible and admits crystal bases. We refer the reader to \cite{Kash} or \cite{BN} and references therein for details. It is also worth to mention that, contrary to the case of ordinary vector bunles, the ``Hecke algebra'' arising from the torsion subalgebra is not identified with the Iwahori-Hecke algebra appeared in the context of automorphic forms with tame ramifications on certain points (c.f. \cite{Hei}). We expect that there will be some analog Schur-Weyl duality to describe the relation between these two different ``Hecke actions'' on the (tensor product of) $\mathbf{U}^>[1]$. 

\vspace{.1in}

To state our main theorem on the structure of the spherical Hall algebra $\mathbf{U}^>$ of parabolic vector bundles, let us introduce some notations: let $S=\{ p_1,\dots,p_N \}$, $V_{p_i}$ be a $w_i$-dimensional $\mathbb{C}$-vector space with basis $\{v_{p_j}^{(i)} \}_{i=0,1,\dots, w_j-1}$, $V := \mathbb{C}[z^{\pm 1}] \otimes_\mathbb{C} \bigotimes_{p \in S} V_p$ and
\[
V^{\widehat{\otimes}r} := \mathbb{C}[z_1^{\pm 1}, \dots, z_r^{\pm 1} ] [[ z_1/z_2, \dots, z_{r-1}/z_r ]] \otimes \bigotimes_{p \in S} V_p^{\otimes r}.
\]
For any fixed rational function $h(z) \in \mathbb{C}(z)$, we define a linear oeprator $\Gamma_{h}(z/z') : V \otimes V \to V \widehat{\otimes} V$ by
\begin{equation*}
\begin{split}
&\Gamma_h(z/z')\big( (p(z)\otimes v_{p_1}^{(i_1)} \otimes \cdots \otimes v_{p_N}^{(i_N)} )  \otimes (q(z') \otimes v_{p_1}^{(j_1)} \otimes \cdots \otimes v_{p_N}^{(j_N)}\big) \\ 
= &\sum_{S' \subseteq S} h(z/z') \big(\frac{1-z/z'}{1-v^2z/z'} v\big)^{\# S'} ((p(z)\otimes v_{p_1}^{(i'_1)} \otimes \cdots \otimes v_{p_N}^{(i'_N)} )  \otimes (q(z') \otimes v_{p_1}^{(j'_1)} \otimes \cdots \otimes v_{p_N}^{(j'_N)}) 
\end{split}
\end{equation*}
where, in each summand,
\begin{equation*}
\begin{matrix}
i'_p = \begin{cases}
i_p + \overline{j_p - i_p} & \textit{ if } p \in S \setminus S' \\
i_p & \textit{ if } p \notin S \setminus S'
\end{cases}, &
j'_p = \begin{cases}
j_p - \overline{j_p - i_p} & \textit{ if } p \in S \setminus S' \\
j_p & \textit{ if } p \notin S \setminus S'
\end{cases}
\end{matrix}
\end{equation*}
and the $\overline{j_p - i_p}$ is known as the residue of $j_p - i_p$ modulo $w_p$. We define an associative algebra $\mathbf{Sh}_{h(z)}(V)$ as follows: $\mathbf{Sh}_{h(z)}(V)$ is isomorphic, as a vector space, to $\bigoplus_{r \geq 0} V^{\widehat{\otimes} r}$.
The product is defined on $\mathbf{Sh}_{h(z)}(V)$ by:
\begin{equation*}
\begin{split}
P(z_1,\dots,z_r) \otimes v \star Q(z_1,\dots ,z_s) \otimes w =
 \sum_{\sigma \in Sh_{r,s}} \prod^{\longrightarrow}_{(i,j) \in I_\sigma} \Gamma^{ij}_h (z_i/z_j)\sigma\big( P(z_1,\dots,z_r) Q(z_{r+1},\dots,z_{r+s})  \otimes v \otimes w \big) 
\end{split}
\end{equation*}
where, $\Gamma_h^{i,j}(z_i/z_j)$ is known as acting on the pair of $i$th and $j$th compononts, the oriented product is with respect to the lexicographical order on the pairs $(i,j)$ and the rational functions $h(z_i/z_j)\big(\frac{1-z_i/z_j}{1-v^2z_i/z_j} v\big)^{\# S'} $ are developed as Laurent series in $z_1/z_2, \dots , z_{r-1}/z_r$. Let $\mathbf{S}_{h(z)} $ be the subalgebra of $\mathbf{Sh}_{h(z)}(V)$ generated by $V$. 
\begin{theo}(Theorem \ref{T:maintheo})
Set $h_X(z) = v^{2g-2} \frac{\zeta_X(z)}{\zeta_x(v^2z)}$. Then there is an algebra homomorphism $J: \mathbf{U}^> \to \mathbf{Sh}_{h_X(z)}(V)$ such that $J(\mathbf{U}^>) \simeq \mathbf{S}_{h_X(z)}$. 
\end{theo}
The above theorem immediately implies the existence of a ``generic form'' (c.f. Corollary \ref{S5:maincor}) and provides the existence of the Kac polynomials $\mathcal{A}_{g,\mathbf{w},r,\mathbf{d}}(\sigma_X)$ for counting the number of (isoclasses of) geometrically indecomposable quasi-parabolic vector bundles of fixed rank and ``multi-degree'' on the curve $X$ over a finite field (\cite[Theorem 7.1]{S5}). The link between the Kac polynomials and the betti number of the moduli spaces of stable parabolic Higgs bundles of the same fixed rank and multi-degree is the following 
\begin{theo}(Theorem \ref{T:Higgs})
Let $\PHiggs^{ss,vec}_{r,\mathbf{d}}(X)$ be the moduli space of stable parabolic Higgs bundles of class $(r,\mathbf{d})$ over $X$ with $(r,\mathbf{d})$ generic. Then 
\begin{equation*}
\mathcal{A}_{g,\mathbf{w},r,\mathbf{d}}(\sigma_X)(q) = v^{\dim \PHiggs^{ss,vec}_{r,\mathbf{d}}(X)} \# \PHiggs^{ss,vec}_{r,\mathbf{d}}(X)(\mathbb{F}_q).
\end{equation*}
\end{theo}
In the recent preprint \cite{DGT} Dobrovolska, Ginzburg and Trakvin provide a much stronger statement (for even non-generic cases) of the above theorem using a different (and more general) approach. However, to derive an explicit combinatorial formula for the Kac polynomial we still need to generalize Schiffmann's computations and this subject will be treated in a companion paper. 

Finally let us come back to the spherical Hall algebra $\mathbf{U}^>$. The shuffle like structure of $\mathbf{U}^>$ seems \emph{new} (we are not able to identify it with other known algebra structures in the literatures) but there are some evidences on the importance of studying $\mathbf{U}^>$ and its representations: In \cite{SV2} Schiffmann and Vasserot identify the spherical Hall algebra of ordinary vector bundles with the $K$-theoretic spherical Hall algebra of genus $g$ commuting varieties and explain the connection via geometric Langlands correspondence of $\GL_r$. There is also a tamely ramified version of Langlands correspondence (c.f. \cite{Dr83, Dr84}) and $\mathbf{U}^>$ fits into the automorphic side of the picture. On the other hand, Hausel and Rodriguez-Villegas derived a conjectural formula for the Poincar\'e polynomial of the moduli spaces of stable (parabolic) Higgs bundles via the arithmetic harmonic analysis of the genus $g$ character varieties \cite{HV} (see also \cite{HLV,HLV2}) and the laters are the categorical quotients of the genus $g$ commuting varieties in the case of ordinary vector bundles. The natural candidate to complete the picture is to consider the $K$-theoretic spherical Hall algebra of these character varieties (or the related quiver varieties). Note that in the case of $\mathbb{P}^1$ and only one marked point, $\mathbf{U}^>$ is identified with the tensor space $T(V(\varpi_1))$ of the affine extremal weight $\mathbf{U}_v(\widehat{\mathfrak{sl}}_w)$-module and the later can be realized as the equivariant $K$-theory of the Nakajima varietires of type $A$ quiver which is indeed the genus $0$ (twised) character varieties. We plan to continue on this direction in the future.

\vspace{.1in}

The content of the paper is as follows: In Sect. 1 we clarify the notion of parabolic coherent sheaves and prove some properties needed but cannot address the references. In Sect. 2 there are some generalities of the Hall algebras of curves. We state and prove our main theorems in Sect. 3. and leave the complicated computation in the Appendix A..

\subsection*{Acknowlegments} The article is part of the author's Ph.D. thesis. The author is grateful to his supervisor Olivier Schiffmann for suggestiing this problem and for his patience and guidance for a long time. Furthermore, the author would like to thank Pierre Baumann, Francesco Sala and Sergey Mozgovoy for their patience to verify the computations and for their useful comments through all the content, Hau-Feng Zhang and Jae-Hoon Kwon for their knowledges on the extremal weight modules of the quantum affine algebras.

\section{Parabolic coherent sheaves on a smooth projective curve}\label{C:curve}

\subsection{Reminder of coherent sheaves on a curve}\label{S2:curve} We fix a smooth projective curve $X$ of genus $g$ defined over a finite field $\mathbf{k}=\mathbb{F}_q$ of $q$ elements and denote by $\overline{X} = X \times_{\Spec \mathbf{k}} \Spec \overline{\mathbf{k}}$ the extension of scalars of $X$ to the algebraic closure $\overline{\mathbf{k}}$ of $\mathbf{k}$. The Galois group $\Gal(\overline{\mathbf{k}} / \mathbf{k})$ acts on the set of all points $\overline{X}(\overline{\mathbf{k}})$. By a \textit{closed point} of $X$ we will mean a $\Gal(\overline{\mathbf{k}} / \mathbf{k})$-orbit of points in $\overline{X}(\overline{\mathbf{k}})$ and we denote by $|X|$ the set of closed points of $X$. The \textit{degree} $\deg (x)$ of a closed point $x \in |X|$ is the number of elements in the associated orbit. Equivalently, if $\mathscr{O}_{X,x}$ stands for the local ring at $x$ of the structure sheaf $\mathscr{O}_X$ of $X$ with the maximal ideal $\mathbf{m}_x$ and if $\mathbf{k}_x = \mathscr{O}_{X,x} / \mathbf{m}_x$ is the residue field, then $\deg(x) = [ \mathbf{k}_x : \mathbf{k} ]$. 

For any $d \geq 1$, let $X(\mathbb{F}_{q^d})$ be the set of $\mathbb{F}_{q^d}$-rational points of $X$. The Galois group $\Gal( \mathbb{F}_{q^d} / \mathbb{F}_q)$ acts on $X(\mathbb{F}_{q^d})$ and the orbits correspond to the closed points whose degree is a divisor of $d$. We have
\begin{equation*}
\# X(\mathbb{F}_{q^d}) = \sum_{n | d} \sum_{\begin{matrix} x \in |X| ,\\ \deg(x) = n \end{matrix}} n.
\end{equation*}

Let us fix a prime number $l$ corpime to $q$ and consider the $l$-adic cohomology group $H^i(\overline{X}, \overline{\mathbb{Q}}_l)$. Then
\begin{equation*}
 \begin{split}
 \dim_{\overline{\mathbb{Q}}_l} H^0( \overline{X}, \overline{\mathbb{Q}}_l) &=1,\\
 \dim_{\overline{\mathbb{Q}}_l} H^1( \overline{X}, \overline{\mathbb{Q}}_l) &= 2g, \\
 \dim_{\overline{\mathbb{Q}}_l} H^2( \overline{X}, \overline{\mathbb{Q}}_l) &= 1.
 \end{split}
\end{equation*}

Let $Fr$ be the geometric Frobenius. $Fr$ acts on $H^i(\overline{X}, \overline{\mathbb{Q}}_l)$ and we denote by $\alpha_1, \dots, \alpha_{2g}$ its eigenvalues in $H^1(\overline{X}, \overline{\mathbb{Q}}_l)$. We will fix an embedding $\overline{\mathbb{Q}}_l \subset \mathbb{C}$ and consider $\alpha_1,\dots, \alpha_{2g}$ as complex numbers. It is known that all the $\alpha_i$ are algebraic numbers satisfying $|\alpha_i| = q^{1/2}$ and these may be reordered in such way that $\alpha_i \alpha_{2g+1-i} = q$ for all $i$. The number of points in $X(\mathbb{F}_{q^d})$ may be expressed in terms of these Frobenius eigenvalues as 
\begin{equation*}
\# X(\mathbb{F}_{q^d}) = \sum_{i=0}^{2} (-1)^i \Tr( Fr , H^i(\overline{X}, \overline{\mathbb{Q}}_l) ) = 1 - \sum_{i} \alpha_i^d + q^d.
\end{equation*}

A convenient way to formulate the above is to introduce the zeta function
\begin{equation*}
\zeta_X (t) = \exp \big( \sum_{d \geq 1} \# X(\mathbb{F}_{q^d}) \frac{t^d}{d} \big).
\end{equation*}
Then Weil's conjecture says
\begin{equation}\label{C:zeta}
\zeta_X (t) = \frac{\prod_{i=1}^{2g} (1 - \alpha_i t)}{(1-t)(1-qt)}.
\end{equation}

Let $\Coh(X)$ be the category of coherent sheaves on $X$. It is of global dimension one, and by a classical theorem of Serre we have $\dim \Hom(\mathscr{F},\mathscr{G}) < \infty$ and $\dim \Ext^1(\mathscr{F},\mathscr{G}) < \infty$ for any two coherent sheaves $\mathscr{F}$ and $\mathscr{G}$. Thus $\Coh(X)$ is finitary. Let $\omega_X = T^*X$ be the canonical bundle of $X$. Then we have $g = \dim H^0(\omega_X)$. The Serre duality is a functorial isomorphism
\begin{equation}\label{C:Serre}
\Ext^*(\mathscr{F}, \mathscr{G}) \cong \Hom( \mathscr{G} \otimes \omega_X^{-1} , \mathscr{F})
\end{equation}
for any two coherent sheaves $\mathscr{F}$ and $\mathscr{G}$. Moreover, any coherent sheaf $\mathscr{F}$ has a canonical torsion subsheaf $\mathscr{F}_t \subset \mathscr{F}$ and canonical quotient vector bundle $\mathscr{F}_v = \mathscr{F} / \mathscr{F}_t$. The exact sequence
\begin{equation*}
0 \to \mathscr{F}_t \to \mathscr{F} \to \mathscr{F}_v \to 0
\end{equation*}
splits, i.e., any coherent sheaf $\mathscr{F}$ can be decomposed as a direct sum $\mathscr{F} = \mathscr{F}_t \oplus \mathscr{F}_v$. 

Let $\Tor(X)$ stand for the full subcategory of $\Coh(X)$ consisting of torsion sheaves. It decomposes into a direct product of blocks 
\begin{equation*}
\Tor(X) = \prod_{ x \in |X| } \Tor_x
\end{equation*}
where $\Tor_x$ is the category of torsion sheaves supported at $x$. This category is Morita equivalent to the category of finite-dimensional modules over the local ring at $x$ which, by the smoothness of $X$, is a discrete valuation ring. In particular, the skyscraper sheaf at $x$, denoted also by $\mathscr{O}_{X,x}$, is the unique simple sheaf supported at $x$. More generally, for each positive integer $n$, there is a unique indecomposable torsion sheaf $\mathscr{O}^{(n)}_{X,x}$. By Krull-Schmidt Theorem, the set of isomorphism classes of torsion sheaves supported at $x$ is in bijection with the set $\Pi$ of all partitions via
\begin{equation*}
(\lambda_1 \geq \cdots \geq \lambda_t) = \lambda \mapsto \mathscr{O}_{X,x}^{(\lambda)} = \mathscr{O}_{X,x}^{(\lambda_1)} \oplus \cdots \oplus \mathscr{O}_{X,x}^{(\lambda_t)}.
\end{equation*}

Let $\Bun(X)$ be the exact subcategory of $\Coh(X)$ consisting of vector bundles on $X$. The \textit{rank} of a coherent sheaf $\mathscr{F}$ is the rank of its canonical quotient vector bundle $\mathscr{F}_v$. The \textit{degree} of a sheaf $\mathscr{F}$ is the only invariant satisfying $\deg \mathscr{O}_X = 0, \deg \mathscr{O}_{X,x} = \deg(x) = [\mathbf{k}_x : \mathbf{k}]$ and which is additive on short exact sequences. Let $\Pic(X)$ be the Picard group of line bundles on $X$. The kernel of the degree map $\deg: \Pic(X) \to \mathbb{Z}$ is denoted by $\Pic^0(X)$ which is a finite abelian group.

A consequence of the Serre duality is that for any $\mathscr{T} \in \Tor(X)$ and $\mathscr{E} \in \Bun(X)$, we have
\begin{equation*}
\Hom(\mathscr{T},\mathscr{E}) = \Ext^1(\mathscr{E},\mathscr{T}) = \{ 0 \}.
\end{equation*}
Moreover, $\Ext^1(\mathscr{O}_{X,x}, \mathscr{O}_X) \simeq \mathbf{k}^{\deg(x)} $ for any closed point $x \in |X|$ and there is a unique line bundle extension of $\mathscr{O}_{X,x}$ by $\mathscr{O}_X$, for which we denote by $\mathscr{O}_X(x)$.

Finally, consider the Grothendieck group $K(X) = K(\Coh(X))$ of $\Coh(X)$. The degree map descends to a group homomorphism $\deg: K(X) \to \mathbb{Z}$ which satisfies $\deg \mathscr{O}_X(x) = \deg(x)$ for all $x \in |X|$. Together with the rank, these maps define a natural morphism
\begin{equation*}
\begin{split}
\varphi: K(X) &\to \mathbb{Z}^2 \\
[\mathscr{F}] &\mapsto ( \rank \mathscr{F}, \deg \mathscr{F} ).
\end{split}
\end{equation*}
By the Riemann-Roch formula
\begin{equation}\label{C:RR}
\langle \mathscr{F} , \mathscr{G} \rangle = (1-g)\rank \mathscr{F} \rank \mathscr{G} + \rank \mathscr{F} \deg \mathscr{G} - \deg \mathscr{F} \rank \mathscr{G},
\end{equation}
the Euler form factors through the map $\varphi$ and the kernel of $\varphi$ is precisely equal to the kernel of the Euler form. Thus we may consider $\mathbb{Z}^2$ as the \textit{numerical} Grothendieck group $K'(\Coh(X))$ of $\Coh(X)$.

\subsection{Parabolic coherent sheaves}\label{S3:defpara} Let us fix a finite set $S= \{ p_1, \dots, p_N \} \subset X(\mathbf{k})$ of $N$ rational points and fix a collection $\mathbf{w} = ( w_p)_{p \in S}$ of positive integers. A coherent sheaves $\mathscr{F}_\bullet$ on $X$ with \textit{$\mathbf{w}$-step (quasi-)parabolic structures on $S$} (also called a quasi-parabolic coherent sheaf for short) is a $N$-dimensional commutative diagram of objects $\mathscr{F}_{\mathbf{a}} \in \Coh(X)$:
\begin{equation}\label{D:par}
\xymatrix{
\mathscr{F}_{\mathbf{a}} \ar[r]^-{\varphi_{\mathbf{a}, \epsilon_p}} \ar[d]^-{\mathbf{a}, \epsilon_q} & \mathscr{F}_{\mathbf{a} + \epsilon_p} \ar[r]^-{\varphi_{\mathbf{a}+\epsilon_p,\epsilon_p}} \ar[d] & \cdots \ar[r] \ar[d] & \mathscr{F}_{\mathbf{a} + (w_p-1)\epsilon_p} \ar[r]^-{\varphi_{\mathbf{a} + (w_p-1)\epsilon_p , \epsilon_p}} \ar[d]^-{\varphi_{\mathbf{a} + (w_p -1) \epsilon_p, \epsilon_q}} & \mathscr{F}_{\mathbf{a}+ w_p \epsilon_p} \ar[d] \\
\mathscr{F}_{\mathbf{a} + \epsilon_q} \ar[r]^-{\varphi_{\mathbf{a}+\epsilon_q, \epsilon_p}} \ar[d] & \mathscr{F}_{\mathbf{a}+\epsilon_p+\epsilon_q} \ar[r] \ar[d] & \cdots \ar[r] \ar[d] & \mathscr{F}_{\mathbf{a} + (w_p-1)\epsilon_p + \epsilon_q} \ar[r] \ar[d] & \mathscr{F}_{\mathbf{a} + w_p \epsilon_p + \epsilon_q} \ar[d] \\
\vdots \ar[d] & \vdots \ar[d] & \vdots \ar[d] & \vdots \ar[d] & \vdots \ar[d] \\
\mathscr{F}_{\mathbf{a} + w_q \epsilon_q} \ar[r] & \mathscr{F}_{\mathbf{a} + w_q \epsilon_q + \epsilon_p} \ar[r] & \cdots \ar[r] & \mathscr{F}_{\mathbf{a} + w_q \epsilon_q + (w_p-1)\epsilon_p} \ar[r] & \mathscr{F}_{\mathbf{a} + w_p \epsilon_p + w_q \epsilon_q}
}
\end{equation}
satisfying the periodic properties $\mathscr{F}_{\mathbf{a} + w_p \epsilon_p } = \mathscr{F}_{\mathbf{a}}(p)$, $\varphi_{\mathbf{a}+ w_p\epsilon_p, \epsilon_q} = \varphi_{\mathbf{a},\epsilon_q}(p)$ and such that the composition
\begin{equation*}
\varphi_{\mathbf{a} + (w_p-1)\epsilon_p, \epsilon_p}  \circ \cdots \circ \varphi_{\mathbf{a}, \epsilon_p} : \mathscr{F}_{\mathbf{a}} \to \mathscr{F}_{\mathbf{a}}(p)
\end{equation*}
is the natural morphism for all $\mathbf{a} \in \mathbb{Z}^S$ and all $p \in S$, where $\epsilon_p$ is the canonical basis of $\mathbb{Z}^S$. A morphism $f_\bullet: \mathscr{F}_\bullet \to \mathscr{G}_\bullet$ between two quasi-parabolic coherent sheaves $\mathscr{F}_\bullet $ and $ \mathscr{G}_\bullet $ is defined in an obvious way: namely, a collection $f_\bullet$ of morphisms $f_{\mathbf{a}} : \mathscr{F}_{\mathbf{a}} \to \mathscr{G}_{\mathbf{a}}$ such that each diagram
\begin{equation*}
\xymatrix{
\mathscr{F}_{\mathbf{a}} \ar[r]^-{\varphi_{\mathbf{a},\epsilon_p}} \ar[d]^-{f_{\mathbf{a}}} & \mathscr{F}_{\mathbf{a}+\epsilon_p} \ar[d]^-{f_{\mathbf{a}+\epsilon_p}} \\
\mathscr{G}_{\mathbf{a}} \ar[r]^-{\varphi'_{\mathbf{a},\epsilon_p}} & \mathscr{G}_{\mathbf{a}+\epsilon_p} 
}
\end{equation*}
commutes for any $\mathbf{a} \in \mathbb{Z}^S$ and any $p \in S$. 

If $\mathscr{F}_{\mathbf{a}+i\epsilon_p}$ is not torsion free at some point $p$ for some $i$, then the natural map $\mathscr{F}_{\mathbf{a}} \to \mathscr{F}_{\mathbf{a}}(p)$ is not injective, hence at least one of the $\varphi_\bullet$ is not injective. By the same argument we can easily see that all the $\varphi_{\mathbf{a},\epsilon_p}$'s are isomorphisms on $X \setminus \{ p \}$ and all the $\mathscr{F}_{\mathbf{a} + i\epsilon_p}$'s have the same generic rank. Thus we define the \textit{rank} $\rank \mathscr{F}_\bullet$ of a quasi-parabolic coherent sheaf $\mathscr{F}_\bullet$ as the rank of $\mathscr{F}_{\mathbf{0}}$. A quasi-parabolic coherent sheaf $\mathscr{F}_\bullet$ is called \textit{locally free} or a \textit{quasi-parabolic vector bundle} if all the components $\mathscr{F}_{\mathbf{a}}$ are locally free and $\mathscr{F}_\bullet$ is called \textit{torsion} if $\rank \mathscr{F}_\bullet = 0 $. Any quasi-parabolic coherent sheaf $\mathscr{F}_\bullet$ can be decomposed as a direct sum $\mathscr{F}^{tor}_{\bullet} \oplus \mathscr{F}^{vec}_{\bullet}$ of a canonical quasi-parabolic torsion subsheaf $\mathscr{F}^{tor}_{\bullet} \subset \mathscr{F}_\bullet$ and a canonical quasi-parabolic quotient bundle $\mathscr{F}^{vec}_{\bullet} = \mathscr{F}_\bullet / \mathscr{F}^{tor}_{\bullet}$ (c.f. \cite[Remark 3.3]{Hei}).

For a quasi-parabolic coherent sheaf $\mathscr{F}_\bullet$, the degrees of the middle terms in the diagram \eqref{D:par} are determined by the degrees of its frames: the sheaf $\mathscr{F}_{\mathbf{a} + i\epsilon_p + j\epsilon_q}$ is isomorphic to $\mathscr{F}_{\mathbf{a} + i\epsilon_p}$ on $X \setminus \{ q \}$ and is isomorphic to $\mathscr{F}_{\mathbf{a} + j \epsilon_q}$ on $X \setminus \{ p \}$. These sheaves glue together since they both isomorphic to $\mathscr{F}_{\mathbf{a}}$ on $X \setminus \{p,q \}$. Therefore we define the \textit{multi-degree} $\deg \mathscr{F}_\bullet$ of $\mathscr{F}_\bullet$ as the collection
\begin{equation*}
\deg \mathscr{F}_\bullet = ( \deg \mathscr{F}_{\mathbf{0}}, \deg \mathscr{F}_{n_p \epsilon_p})_{p \in S, 1 \leq n_p \leq w_p-1}.  
\end{equation*}

We fix a collection $\underline{a} = (a_{p,i})_{p \in S, 1 \leq i \leq w_p-1}$ of real numbers such that 
\begin{equation}\label{E:stab}
0 \leq a_{p,w_p-1} < a_{p,w_p-2} < \cdots < a_{p,1} <1
\end{equation}
for all $p \in S$. For any quasi-parabolic coherent sheaf $\mathscr{F}_\bullet$, we define the \textit{parabolic degree} $\Par \mathscr{F}_\bullet$ of $\mathscr{F}_\bullet$ as
\begin{equation*}
\Par \mathscr{F}_\bullet = \deg \mathscr{F}_{\mathbf{0}} + \sum_{p \in S} \sum_{1 \leq i \leq w_p-1} a_{p,i}( \deg \mathscr{F}_{i\epsilon_p} - \deg \mathscr{F}_{(i-1)\epsilon_p} ).
\end{equation*}
The \textit{parabolic slope} of $\mathscr{F}_\bullet$ is defined as 
\begin{equation*}
\mu(\mathscr{F}_\bullet) = \frac{\Par \mathscr{F}_\bullet}{\rank \mathscr{F}_\bullet} \in \mathbb{Q} \cup \{ \infty \}.
\end{equation*}
As usual, an $\mathscr{F}_\bullet$ is said to be \textit{semi-stable} (resp. \textit{stable}) of slope $\nu$ if $\mu(\mathscr{F}_\bullet) = \nu$ and if $\mu(\mathscr{G}_\bullet) \leq \nu$ (resp. $\mu(\mathscr{G}_\bullet) < \nu$) for any proper subsheaf $\mathscr{G}_\bullet \subset \mathscr{F}_\bullet$. 

\subsection{Homological properties}\label{S3:PCoh} Let $\PCoh(X)$ be the category of coherent sheaves on $X$ with $\mathbf{w}$-step quasi-parabolic structure on $S$ within the chosen stability condition $\underline{a}$. We will call an object $\mathscr{F}_\bullet$ of $\PCoh(X)$ a \textit{parabolic coherent sheaf} for convenience and for any pair of parabolic coherent sheaves $(\mathscr{F}_\bullet,\mathscr{G}_\bullet)$ we denote by $\PHom(\mathscr{F}_\bullet, \mathscr{G}_\bullet)$ the set of morphisms from $\mathscr{F}_\bullet$ to $\mathscr{G}_\bullet$ (and similarly for $\PExt^i(\mathscr{F}_\bullet, \mathscr{G}_\bullet)$ for all $i$). It is obvious the $\PCoh(X)$ is an abelian category where the kernel and the cokernel of a morphism can be defined componentwisely. The compatibilities thus follow from the corresponding ones of $\Coh(X)$. Moreover we have 
\begin{lem}\cite[Lemma 2.1, 3.4]{Hei}\cite[Lemma 4.2, Theorem 4.3]{Len98}
$\PCoh(X)$ is connected, finitary, abelian, noetherian and of global dimension one. Moreover, $\PCoh(X)$ has enough injectives. 
\end{lem}

\noindent
\textbf{Remarks.} 
\begin{enumerate}
\item[(1)] If $\mathscr{E}_\bullet \in \PCoh(X)$ is a parabolic vector bundle, all the arrows in the diagram \eqref{D:par} are inclusions. We may hence interpret $\mathscr{E}_\bullet$ as a filtration
\[
\mathscr{E}_{\mathbf{0}} \subseteq \mathscr{E}_{\epsilon_p} \subseteq \mathscr{E}_{2\epsilon_p} \subseteq \cdots \subset \mathscr{E}_{(w_p-1)\epsilon_p} \subseteq \mathscr{E}_{w_p \epsilon_p} = \mathscr{E}_{\mathbf{0}}(p)
\]
or equivalently as a filtration
\[
0 = \mathscr{E}_{\mathbf{0}} / \mathscr{E}_{\mathbf{0}} \subseteq \mathscr{E}_{\epsilon_p} / \mathscr{E}_{\mathbf{0}} \subseteq \cdots \subseteq \mathscr{E}_{(w_p-1)\epsilon_p} / \mathscr{E}_{\mathbf{0}} \subseteq \mathscr{E}_{\mathbf{0}}(p) / \mathscr{E}_{\mathbf{0}}
\]
of the fiber $(\mathscr{E}_{\mathbf{0}})_p = \mathscr{E}_{\mathbf{0}}(p) / \mathscr{E}_{\mathbf{0}}$ of $\mathscr{E}_{\mathbf{0}}$ at $p$. In other word, our definition of a parabolic vector bundle coincides with the one in the usual sense of Mehta and Seshadri(c.f. \cite{MS}). The presented generalization to coherent sheaves is due to Heinloth (see \cite{Hei}). 
\item[(2)] It is easy to see that the morphisms in $\PCoh(X)$ between parabolic bundles are necessary the \textit{strongly} ones since we equip all the quasi-parabolic vector bundles with the same \textit{weights} (cf. \cite{GPGM}).

\end{enumerate}

There is a natural $\mathbb{Z}^S$-action on $\PCoh(X)$ as the shift automorphism $\mathbf{a}: \PCoh(X) \to \PCoh(X)$ defined by 
\[
\mathbf{a} \cdot \mathscr{F}_\bullet = \mathbf{a} \cdot (\mathscr{F}_{\mathbf{b}})_{\mathbf{b} \in \mathbb{Z}^S} = (\mathscr{F}_{\mathbf{b} + \mathbf{a}})_{\mathbf{b} \in \mathbb{Z}^S} =: \mathscr{F}_\bullet(\mathbf{a})
\]
for all $\mathbf{a} \in \mathbb{Z}^S$. The Picard group $\Pic(X)$ of $\Coh(X)$ also acts on $\PCoh(X)$ as automorphisms:
\begin{equation*}
\begin{split}
\Pic(X) \times \PCoh(X) &\to \PCoh(X) \\
(\mathscr{L}, \mathscr{F}_\bullet = (\mathscr{F}_{\mathbf{a}})_{\mathbf{a} \in \mathbb{Z}^S} ) &\mapsto \mathscr{L} \cdot \mathscr{F}_\bullet := (\mathscr{F}_{\mathbf{a}} \otimes \mathscr{L})_{\mathbf{a} \in \mathbb{Z}^S}.
\end{split}
\end{equation*}
Clearly these two actions commute and satisfy the relations 
\begin{equation*}
\mathscr{O}_X(p) \cdot \mathscr{F}_\bullet = w_p \epsilon_p \cdot \mathscr{F}_\bullet = \mathscr{F}_\bullet(w_p \epsilon_p)
\end{equation*}
for all $p \in S$.

We will say a parabolic coherent sheaf $\mathscr{F}_\bullet \in \PCoh(X)$ is \textit{constant} if $\mathscr{F}_\bullet$ is of the form
\begin{equation*}
\xymatrix{
\mathscr{F}_{\mathbf{0}} = \mathscr{F} \ar[r] \ar[d] & \mathscr{F} \ar[r] \ar[d] & \cdots \ar[r] \ar[d] & \mathscr{F} \ar[r] \ar[d] & \mathscr{F}(p) \ar[d] \\
\mathscr{F} \ar[r] \ar[d] & \mathscr{F} \ar[r] \ar[d] & \cdots \ar[r] \ar[d] & \mathscr{F} \ar[r] \ar[d] & \mathscr{F}(p) \ar[d] \\
\vdots \ar[r] \ar[d] & \vdots \ar[r] \ar[d] & \vdots \ar[r] \ar[d] & \vdots \ar[r] \ar[d] & \vdots \ar[d] \\
\mathscr{F} \ar[r] \ar[d] & \mathscr{F} \ar[r] \ar[d] & \cdots \ar[r] \ar[d] & \mathscr{F} \ar[r] \ar[d] & \mathscr{F}(p) \ar[d] \\
\mathscr{F}(q) \ar[r] & \mathscr{F}(q) \ar[r] & \cdots \ar[r] & \mathscr{F}(q) \ar[r] & \mathscr{F}(p+q)
}
\end{equation*}
for some coherent sheaf $\mathscr{F} \in \Coh(X)$ and we denote by $(\mathscr{F})_\bullet$ such a parabolic coherent sheaf. The category $\Coh(X)$ can be embedded into $\PCoh(X)$ by assigning a coherent sheaf $\mathscr{F} \in \Coh(X)$ to the constant parabolic sheaf $(\mathscr{F})_\bullet$. Conversely there is a natural functor $( \; )_{\mathbf{a}} : \PCoh(X) \to \Coh(X)$ for any $\mathbf{a} \in \mathbb{Z}^S$ by sending $\mathscr{F}_\bullet = (\mathscr{F}_{\mathbf{a}})_{\mathbf{a} \in \mathbb{Z}^S}$ to $\mathscr{F}_{\mathbf{a}}$. These two functors $( \;  )_\bullet$ and $( \; )_{\mathbf{a}}$ satisfy the following adjunction properties:
\begin{equation}\label{E:adjhom}
\begin{split}
\PHom((\mathscr{F})_\bullet, \mathscr{G}_\bullet) &= \Hom(\mathscr{F}, \mathscr{G}_{\mathbf{0}} ) \\
\PHom(\mathscr{G}_\bullet, (\mathscr{F})_\bullet) &= \Hom(\mathscr{G}_{\mathbf{w}}, \mathscr{F})
\end{split}
\end{equation}
where we set $\mathbf{w} := \sum_{p \in S} (w_p - 1) \epsilon_p$ and hence $\mathscr{G}_{\mathbf{w}} = \mathscr{G}_{\sum_{p \in S} (w_p -1)\epsilon_p}$. Similarly for the extension spaces 
\begin{equation}\label{E:adjext}
\begin{split}
\PExt^1((\mathscr{F})_\bullet, \mathscr{G}_\bullet) &= \Ext^1(\mathscr{F}, \mathscr{G}_{\mathbf{0}} ) \\
\PExt^1(\mathscr{G}_\bullet, (\mathscr{F})_\bullet) &= \Ext^1(\mathscr{G}_{\mathbf{w}}, \mathscr{F})
\end{split}
\end{equation}

\subsection{Tensor product by parabolic line bundles}\label{S3:tensor} We call an $\mathscr{L}_\bullet \in \PCoh(X)$ a \textit{parabolic line bundle} if $\mathscr{L}_\bullet$ is locally free and of rank one. Observe that any parabolic line bundle $\mathscr{L}_\bullet$ is isomorphic to a suitable shift of a constant parabolic line bundle, i.e., $\mathscr{L}_\bullet \simeq (\mathscr{L'})_\bullet(\mathbf{a})$ for some line bundle $\mathscr{L}' \in \Pic(X)$ and $\mathbf{a} \in \mathbb{Z}^S$. 

There is a well-behaved notion of tensor product of two parabolic vector bundles, which is best understood in terms of $\mathbb{R}$-filtered sheaves (c.f. \cite{Yo}). However, we shall only use the case of tensoring a parabolic coherent sheaf $\mathscr{F}_\bullet$ by a parabolic line bundle $\mathscr{L}_\bullet$ and it is translated in our language as follows: let $\mathscr{L}_\bullet = (\mathscr{L})_\bullet(\mathbf{a})$ for some $\mathscr{L} \in \Pic(X)$ and some $\mathbf{a} \in \mathbb{Z}^S$. Then 
\begin{equation*}
\mathscr{F}_\bullet \otimes (\mathscr{L})_\bullet(\mathbf{a}) := (\mathscr{F}_{\mathbf{b}} \otimes \mathscr{L})_{\mathbf{b} \in \mathbb{Z}^S} (\mathbf{a}).
\end{equation*}

A direct consequence of the Serre duality (c.f. \eqref{C:Serre}) and the adjunction properties \eqref{E:adjhom} and \eqref{E:adjext} is that for any $\mathscr{F} \in \Coh(X)$ and any $\mathscr{G}_\bullet \in \PCoh(X)$, we have 
\begin{equation}\label{C:weakSerre}
\begin{split}
\PExt^1( (\mathscr{F})_\bullet, \mathscr{G}_\bullet)^* &\simeq \PHom( \mathscr{G}_\bullet \otimes \omega_{\bullet}^{-1}, (\mathscr{F})_\bullet )\\
\PExt^1( \mathscr{G}_\bullet, (\mathscr{F})_\bullet)^* &\simeq \PHom( (\mathscr{F})_\bullet \otimes \omega_{\bullet}^{-1}, \mathscr{G}_\bullet )
\end{split}
\end{equation}
where we set $\omega_{\bullet}^k = (\omega_X^{\otimes k})_\bullet(k\mathbf{w})$ for each $k \in \mathbb{Z}$. It will be also convenient to set $\mathscr{O}_{\bullet} := (\mathscr{O}_X)_\bullet$.

Let $\PBun(X)$ be the full subcategory of $\PCoh(X)$ consisting of parabolic vector bundles. For any parabolic vector bundle $\mathscr{E}_\bullet \in \PBun(X)$, there is also a well-defined notion of the \textit{dual} parabolic bundle $\mathscr{E}_\bullet^*$ which can be also translated in our language as follows: we first consider the ``naive'' dual $\mathscr{E}'^*_\bullet = (\mathscr{E}'^*_{\mathbf{a}})_{\mathbf{a} \in \mathbb{Z}^S}$ with $\mathscr{E}'^*_{\mathbf{a}} = \mathscr{E}^*_{-\mathbf{a}}$ for all $\mathbf{a} \in \mathbb{Z}^S$. The dual of a parabolic vector bundle $\mathscr{E}_\bullet$ is defined to be $\mathscr{E}^*_\bullet : = \mathscr{E}'^*_\bullet(\mathbf{-w})$. It is easy to check that $\mathscr{E}^{**}_\bullet =\mathscr{E}_\bullet$ and $\Par \mathscr{E}^*_\bullet = - \Par \mathscr{E}_\bullet$. Moreover, we have the following parabolic version of Serre duality:

\begin{prop}[Serre duality]\cite[Proposition 3.7]{Yo}\label{P:Serre}
For any pair $\mathscr{F}_\bullet, \mathscr{G}_\bullet \in \PCoh(X)$, we have
\begin{equation*}
\PExt(\mathscr{F}_\bullet, \mathscr{G}_\bullet)^* \simeq \PHom(\mathscr{G}_\bullet \otimes \omega^{-1}_\bullet, \mathscr{F}_\bullet)
\end{equation*}
\end{prop}

\begin{prop}\label{P:vecfilt}
Any parabolic vector bundle $\mathscr{E}_\bullet$ has a filtration
\begin{equation}\label{E:vecfilt}
0 = \mathscr{E}^0_{\bullet} \subset \mathscr{E}^1_\bullet \subset \cdots \subset \mathscr{E}^n_\bullet = \mathscr{E}_\bullet
\end{equation}
whose factors $\mathscr{E}^i_\bullet / \mathscr{E}^{i-1}_\bullet$ are parabolic line bunles.
\end{prop}
\noindent
\textit{Proof.} We proceed by induction on the $\rank \mathscr{E}_\bullet$. By the adjunction properties \eqref{E:adjhom}, we may assume by a suitable shift that $\PHom(\mathscr{L}_\bullet, \mathscr{E}_\bullet) \neq 0$, where $\mathscr{L}_\bullet$ is a parabolic line bundle. Hence $\mathscr{L}_\bullet$ may be viewed as a subsheaf of $\mathscr{E}_\bullet$. Note that every parabolic coherent sheaf $\mathscr{F}_\bullet$ can be decomposed into a direct sum of a parabolic vector bundle $\mathscr{F}^{vec}_{\bullet}$ and a parabolic torsion sheaf $\mathscr{F}^{tor}_{\bullet}$. If $\mathscr{E}^1_\bullet \subset \mathscr{E}_\bullet$ is chosen such that $\mathscr{E}^1_\bullet / \mathscr{L}_\bullet \simeq (\mathscr{E}_\bullet / \mathscr{L}_\bullet)^{tor}$, then $\mathscr{E}^1_\bullet$ is clearly a parabolic line bundle and $\mathscr{E}_\bullet / \mathscr{E}^1_\bullet$ is again a parabolic vector bundle of rank less than $\rank \mathscr{E}_\bullet$. The assertion now follows from the induction hypothesis.

\qed

\subsection{Parabolic torsion sheaves}\label{S3:tor} We denote by $\PTor(X)$ the subcategory of $\PCoh(X)$ consisting of parabolic torsion sheaves. As the usual case without parabolic structure, it splits into a direct product of blocks
\begin{equation}\label{E:torsiondec}
\PTor(X) = \prod_{x \in |X|} \PTor_x(X)
\end{equation}
where $\PTor_x(X)$ is the subcategory of parabolic torsion sheaves supported at a single point $x \in |X|$. If $x \notin S$ and $\mathscr{T}_\bullet \in \PTor_x(X)$, then all the $\mathscr{T}_{\mathbf{a}}$ are isomorphic and thus $\mathscr{T}_\bullet = (\mathscr{T}_{\mathbf{0}})_\bullet$. Thus $\PTor_x(X)$ is equivalent to $\Tor_x(X)$, the subcategory of torsion sheaves (without extra structure) on $X$ supported at a single point $x \in |X|$. Hence $\PTor_x(X)$ is equivalent to the category of $\mathscr{O}_{X,x}$-mod of finite-length modules over the local ring $\mathscr{O}_{X,x}$ of $x$, which is a discrete valuation ring and we denote by $\mathscr{O}_{x,\bullet}^{(\lambda)} := (\mathscr{O}_{X,x}^{(\lambda)})_\bullet$. If $p \in S$, any indecomposable parabolic torsion sheaf supported at $p$ will be of the form $\mathscr{O}_{p,\bullet}^{(k)}(i\epsilon_{p}) := \mathscr{O}_\bullet(i\epsilon_{p}) / \mathscr{O}_\bullet((i-k)\epsilon_{p})$ for some $0 \leq i < w_p$ and $k > 0$ (c.f. \cite[Section 3]{Hei}). 

\vspace{.1in}
It will be convenient for us to interpret the categories $\PTor_x(X)$ as the categories of representations of quivers: let $C_n$ denote the quiver of type $A_{n-1}^{(1)}$ with cyclic orientation. We also let $C_1$ be the Jordan quiver, i.e., the quiver with one vertex and one loop. A representation of $C_n$ over a field $\mathbf{k}$ is a collection of $\mathbf{k}$-linear vector spaces $V_i$ with $i \in \mathbb{Z} / n \mathbb{Z}$ together with a collection of linear maps $x_i: V_i \to V_{i-1}$. Morphisms between two representations $(\underline{V},\underline{x})$ and $(\underline{W},\underline{y})$ are $\mathbf{k}$-linear maps $\phi_i: V_i \to W_i $ such that $\phi_i x_i = y_i \phi_i$. Finally, a representation $(\underline{V},\underline{x})$ is called \textit{nilpotent} if there exists $M \gg 0$ such that $x_{i+M} \cdots x_i=0$ for all $i$. 

The set of isomorphism classes of indecomposable nilpotent representations of $C_1$ over $\mathbf{k}$ is in bijection with $\mathbb{N}^*$ and we denote by $|m)$ the class of the indecomposable nilpotent representations of dimension $m$. The set of isomorphism classes of nilpotent representations is thus in bijection with the set $\Pi$ of all partitions via the assignment $\lambda = (\lambda_1 \geq \cdots \geq  \lambda_t) \mapsto |\lambda) := |\lambda_1) \oplus \cdots \oplus |\lambda_t)$. For $x \notin S$, the assignment $\mathscr{O}_{x,\bullet}^{(m)}:= (\mathscr{O}_{X,x}^{(m)})_\bullet \mapsto |m)$ gives rise to a Morita equivalence of $\PTor_x(X)$ and the category $\Rep^{nil}_{\mathbf{k}}(C_1)$ of nilpotent representations of $C_1$ over $\mathbf{k}$.

Now consider the quiver $C_n$ for $n>1$. Denote by $(\epsilon_i)_{i \in \mathbb{Z} / n\mathbb{Z}}$ the canonical basis of $\mathbb{Z}^{\mathbb{Z}/n\mathbb{Z}}$. For each $i \in \mathbb{Z} / n \mathbb{Z}$ and $k \in \mathbb{N}^*$, define the \textit{cyclic segment} $[i;k)$ to be the image of the projection to $\mathbb{Z}/n\mathbb{Z}$ of the segment $[i';i'-k+1]$ for any $i' \in \mathbb{Z}$ with $i' \equiv i \mod n $. A \textit{cyclic multisegment} is a finite linear combination $\mathbf{m}= \sum_{i,k} a_{i,k} [i;k)$ with $a_{i,k} \in \mathbb{N}$. The isomorphism classes of nilpotent representations (resp. indecomposable representations) of $\C_n$ over $\mathbf{k}$ is in bijection with the set of cyclic multisegments (resp. cyclic segments), c.f. \cite{S2}. For $p \in S$, the assigment $\mathscr{O}_{p,\bullet}^{(k)}(i\epsilon_p) \mapsto [i;k)$ give rise to a Morita equivalence between $\PTor_p(X)$ and $\Rep^{nil}_{\mathbf{k}}(C_{w_p})$. 

\subsection{The numerical Grothendieck group}\label{S3:Kgroup} The Grothendieck group $K(\PCoh(X))$ can be described as follows: recall from \autoref{S3:defpara} that every parabolic coherent sheaf is a direct sum of a parabolic vector bundle and a parabolic torsion sheaf. By the Proposition \ref{P:vecfilt}, the class of a parabolic vector bundle in $K(\PCoh(X))$ is determined by the classes of parabolic line bundles appeared in the filtration \eqref{E:vecfilt} and every parabolic torsion sheaf is an extension of simple parabolic torsion sheaves. Therefore the classes $[\mathscr{L}_\bullet]$ of parabolic line bundles $\mathscr{L}_\bullet$ form a system of generators for $K(\PCoh(X))$. We denote by $K^+(\PCoh(X)) \subset K(\PCoh(X))$ the subset of $K(\PCoh(X))$ consisting of the classes of parabolic coherent sheaves. We will denote by $\langle \; , \; \rangle_{par}$ the corresponding Euler form on $\PCoh(X)$ to distinguish the usual Euler form $\langle \; , \; \rangle$ on the category $\Coh(X)$ of coherent sheaves without parabolic structures, i.e., for any two parabolic coherent sheaves $\mathscr{F}_\bullet, \mathscr{G}_\bullet$, 
\begin{equation*}
\langle \mathscr{F}_\bullet, \mathscr{G}_\bullet \rangle_{par} := \dim \PHom(\mathscr{F}_\bullet,\mathscr{G}_\bullet) - \dim \PExt^1(\mathscr{F}_\bullet, \mathscr{G}_\bullet)
\end{equation*}
Let us first compute some easy cases of the Euler form:
\begin{lem}\label{L:euler}
For any $p,p' \in S, 0 \leq i < w_p$ and $ 0 \leq j < w_{p'}$, we have 
\begin{equation*}
\begin{split}
&\langle \mathscr{O}_\bullet, (\mathscr{O}(p) )_\bullet \rangle_{par} = 2-g, \; \langle \mathscr{O}_\bullet, \mathscr{O}(i\epsilon_p) \rangle_{par} = 1-g \\
&\langle \mathscr{O}_\bullet, \mathscr{O}^{(1)}_{p,\bullet}(i\epsilon_p) \rangle_{par} = 
\begin{cases}
1 & \textit{ if } i=0,\\
0 & \textit{ otherwise }
\end{cases} \\
&\langle \mathscr{O}^{(1)}_{p,\bullet}(i \epsilon_p) , \mathscr{O}_\bullet \rangle_{par} = 
\begin{cases}
-1 & \textit{ if } i = 1, \\
0 & \textit{ otherwise }
\end{cases} \\
&\langle \mathscr{O}^{(1)}_{p,\bullet}(i\epsilon_p) , \mathscr{O}^{(1)}_{p',\bullet}(j \epsilon_{p'}) \rangle_{par} =
\begin{cases}
1 & \textit{ if } p=p', i=j \\
-1 & \textit{ if } p = p', i \equiv j+1 \mod w_p \\
0 & \textit{ otherwise }
\end{cases}
\end{split}
\end{equation*}
\end{lem}
\noindent
\textit{Proof.} By the adjunction properties \eqref{E:adjhom}, \eqref{E:adjext} and the ordinary Riemann-Roch formula \eqref{C:RR}, we have 
\begin{equation*}
\langle \mathscr{O}_\bullet, (\mathscr{O}(p) )_\bullet \rangle_{par} = \langle \mathscr{O}_X, \mathscr{O}_X(p) \rangle = 2-g.
\end{equation*}
Similarly, since $0 \leq i < w_p$, $(\mathscr{O}_\bullet(i\epsilon_p) )_{\mathbf{0}} =( \mathscr{O}_\bullet(-i\epsilon_p))_{\mathbf{w}} = \mathscr{O}_X$, we have
\begin{equation*}
\langle \mathscr{O}_\bullet, \mathscr{O}_\bullet(i\epsilon_p) \rangle_{par}  = \langle \mathscr{O}_X, \mathscr{O}_X \rangle = 1-g.
\end{equation*}
Note that $(\mathscr{O}^{(1)}_{p,\bullet}(i\epsilon_p))_{j\epsilon_p} = \delta_{w_p-i,j} \mathscr{O}_{X,p}$. Again by the adjunction properties, we have
\begin{equation*}
\begin{split}
&\PHom(\mathscr{O}^{(1)}_{p,\bullet}(i\epsilon_p), \mathscr{O}_\bullet) = \PExt(\mathscr{O}_\bullet, \mathscr{O}^{(1)}_{p,\bullet}(i\epsilon_p) ) = 0, \\
& \PHom(\mathscr{O}_\bullet, \mathscr{O}^{(1)}_{p,\bullet}(i\epsilon_p) ) = 
\begin{cases}
\Hom(\mathscr{O}_\X, \mathscr{O}_{X,p}) & \textit{ if } i = 0 \\
\{ 0 \} & \textit{ if } i \neq 0 
\end{cases}, \\
& \PExt(\mathscr{O}^{(1)}_{p,\bullet} (i\epsilon_p) , \mathscr{O}_\bullet) = 
\begin{cases}
\Ext( \mathscr{O}_{X,p}, \mathscr{O}_X ) & \textit{ if } i=1 \\
\{ 0 \} & \textit{ if } i \neq 1
\end{cases}.
\end{split}
\end{equation*}
Finally, the computation of the last assertion is identical to the corresponding cyclic quiver (c.f. \cite[Lecture 3]{S1}) and we omit it. 
\qed

\vspace{.1in}

In general, we can deduce the Riemann-Roch formula for the parabolic coherent sheaves from the above computations:
\begin{cor}[Riemann-Roch formula]\label{C:parRR}
For any two parabolic coherent sheaves $\mathscr{F}_\bullet, \mathscr{G}_\bullet$, we have
\begin{equation}\label{E:parRR}
\langle \mathscr{F}_\bullet, \mathscr{G}_\bullet \rangle_{par} = \langle \mathscr{F}_{\mathbf{0}}, \mathscr{G}_{\mathbf{0}} \rangle 
+ \sum_{p \in S } \sum_{i = 1}^{w_p-1} \big( \deg \mathscr{F}_{i\epsilon_p} - \deg \mathscr{F}_{\mathbf{0}} \big) \cdot \big( \deg \mathscr{G}_{i\epsilon_p} - \deg \mathscr{G}_{(i+1)\epsilon_p} \big).
\end{equation}
\end{cor}
\noindent
\textit{Proof.} Note that the Euler form $\langle \mathscr{F}_\bullet, \mathscr{G}_\bullet \rangle_{par}$ of $\mathscr{F}_\bullet$ and $\mathscr{G}_\bullet$ depends only on their classes in $K(\PCoh(X))$. From the observation at the beginning of this subsection and the fact that any parabolic line bundle is a shift of a constant parabolic line bundle (see \autoref{S3:tensor}), we may write
\begin{equation*}
\begin{split}
&[\mathscr{F}_\bullet] = \sum_{j=1}^{\rank \mathscr{F}_\bullet} [\mathscr{L}^j_\bullet] + \deg \mathscr{F}_{\mathbf{0}} \delta + \sum_{p \in S} \sum_{i=1}^{w_p-1} \big( \deg \mathscr{F}_{i\epsilon_p} -\deg \mathscr{F}_{\mathbf{0}} \big) [ \mathscr{O}^{(1)}_{p,\bullet}(-i\epsilon_p) ] \\
&[\mathscr{G}_\bullet] = \sum_{j=1}^{\rank \mathscr{G}_\bullet} [\mathscr{L}'^j_\bullet] + \deg \mathscr{G}_{\mathbf{0}} \delta + \sum_{p \in S} \sum_{i=1}^{w_p-1} \big( \deg \mathscr{G}_{i\epsilon_p} -\deg \mathscr{G}_{\mathbf{0}} \big) [ \mathscr{O}^{(1)}_{p,\bullet}(-i\epsilon_p) ] 
\end{split}
\end{equation*}
where $\mathscr{L}^j_\bullet, \mathscr{L}'^j_\bullet$ are constant parabolic line bundles of constant multi-degree $(0,\dots, 0 )$ for all $j,k$ and $\delta = \sum_{i=0}^{w_p-1} [ \mathscr{O}^{(1)}_{p,\bullet}(i\epsilon_p) ]$ for all $p \in S$. By Lemma \autoref{L:euler} and the bilinearity of the Euler form, we have
\begin{equation*}
\begin{split}
\langle \mathscr{F}_\bullet, \mathscr{G}_\bullet \rangle_{par}  &= (1-g) \rank \mathscr{F}_\bullet \rank \mathscr{G}_\bullet \\
&+ \rank \mathscr{F}_\bullet \deg \mathscr{G}_{\mathbf{0}} - \rank \mathscr{G}_\bullet \deg \mathscr{F}_{\mathbf{0}}  - \rank \mathscr{G}_\bullet \sum_{p \in S} \big( \deg \mathscr{F}_{(w_p-1)\epsilon_p} - \deg \mathscr{F}_{\mathbf{0}} \big) \\
&+ \sum_{p \in S} \sum_{i =1 }^{w_p-1} \big( \deg \mathscr{F}_{i \epsilon_p} - \deg \mathscr{F}_{\mathbf{0}} \big) \big( \deg \mathscr{G}_{i \epsilon_p} - \deg \mathscr{G}_{\mathbf{0}}  - \deg \mathscr{G}_{(i+1)\epsilon_p} + \deg \mathscr{G}_{\mathbf{0}} \big) \\
&+ \sum_{p \in S} \big( \deg \mathscr{F}_{(w_p-1)\epsilon_p} - \deg \mathscr{F}_{\mathbf{0}} \big) \big( \deg \mathscr{G}_{w_p\epsilon_p} - \deg \mathscr{G}_{\mathbf{0}} \big) \\
&= (1-g) \rank \mathscr{F}_\bullet \rank \mathscr{G}_\bullet + \rank \mathscr{F}_\bullet \deg \mathscr{G}_{\mathbf{0}}  - \rank \mathscr{G}_\bullet \deg \mathscr{F}_{\mathbf{0}} \\
&+ \sum_{p \in S} \sum_{i =1 }^{w_p-1} \big( \deg \mathscr{F}_{i \epsilon_p} - \deg \mathscr{F}_{\mathbf{0}} \big) \big( \deg \mathscr{G}_{i \epsilon_p}  - \deg \mathscr{G}_{(i+1)\epsilon_p} \big) \\
&= \langle \mathscr{F}_{\mathbf{0}}, \mathscr{G}_{\mathbf{0}} \rangle 
+ \sum_{p \in S } \sum_{i = 1}^{w_p-1} \big( \deg \mathscr{F}_{i\epsilon_p} - \deg \mathscr{F}_{\mathbf{0}} \big) \cdot \big( \deg \mathscr{G}_{i\epsilon_p} - \deg \mathscr{G}_{(i+1)\epsilon_p} \big).
\end{split}
\end{equation*}
\qed

Consequencely, the Euler form is determined by the ranks and multi-degrees of parabolic coherent sheaves involved. As in the case of $\Coh(X)$, We may rather consider the \textit{numerical} Grothendieck group $\mathcal{K} = \mathbb{Z} \times \mathbb{Z} \times \prod_{p \in S} \mathbb{Z}^{w_p-1}$ by the assignment 
\begin{equation*}
\begin{split}
\varphi: K(\PCoh(X)) & \to \mathcal{K} \\
[\mathscr{F}_\bullet] &\mapsto ( \rank \mathscr{F}_\bullet, \deg \mathscr{F}_\bullet) .
\end{split}
\end{equation*}
It may be convenient to set $\alpha_{p,nw_p+i} = \varphi([\mathscr{O}_{p,\bullet}^{(1)}(i)])$ for all $p \in S , \; 0 \leq i \leq w_p-1, \; n \in \mathbb{Z}$ and $\delta = \sum_{i=0}^{w_p-1} \alpha_{p,i}$ for all $p$. We denote by $\mathcal{K}^+$ the image of $K^+(\PCoh(X))$ via the above map. 

Let
\begin{equation*}
\mathcal{K}^{+,tor} := \{ \varphi([\mathscr{T}_\bullet]) \mid \mathscr{T}_\bullet \in \PTor(X) \} = \{ \sum_{p \in S} \sum_{0 \leq i \leq w_p-1} \mathbb{N} \alpha_{p,i}  \} \setminus \{(0,\dots,0) \}.
\end{equation*}
and 
\begin{equation*}
\begin{matrix}
\mathcal{K}^{+,vec}_r := \{ \varphi([\mathscr{E}_\bullet]) \mid \mathscr{E}_\bullet \in \PBun_r(X) \}, & \mathcal{K}^{+,vec} := \{ \varphi([\mathscr{E}_\bullet]) \mid \mathscr{E}_\bullet \in \PBun(X) \} = \cup_{r \geq 1} \mathcal{K}^{+,vec}_r
\end{matrix}
\end{equation*}
where $\PBun_r(X)$ stands for the set of (isomorphic classes of) all parabolic vector bundles of rank $r$. Note that $(r,\mathbf{d}) \in \mathcal{K}^{+,vec}_r$ only if 
\begin{equation}\label{D:flagtype}
d_0 \leq d_{\epsilon_p} \leq \cdots \leq d_{(w_p-1)\epsilon_p} \leq d_0+r
\end{equation}
for all $p \in S$. For any pair $(r,\mathbf{d}),(r,\mathbf{d}') \in \mathcal{K}^{+,vec}_r$, we may define
\begin{equation}\label{E:partialorder}
(r,\mathbf{d}) \leq (r,\mathbf{d}') \textit{ if } (0,\mathbf{d}'-\mathbf{d}) \in \mathcal{K}^{+,tor} \cup \{ (0,\dots,0) \}.
\end{equation}
It is easy to check that $\Par ((r,\mathbf{d})) \leq \Par ((r,\mathbf{d}'))$ if $(r,\mathbf{d}) \leq (r,\mathbf{d}')$.

\vspace{.1in}

For any $(r,\mathbf{d}) \in \mathcal{K}^{+,vec}_r$, we may sometimes write $\mathbf{d} = (d_{\mathbf{0}},\mathbf{m}^{\mathbf{d}}) := d_{\mathbf{0}} \delta + \sum_{p \in S} \sum_{1 \leq i \leq w_p-1} m^{\mathbf{d}}_{p,i} \alpha_{p,i}$ with $m^{\mathbf{d}}_{p,i} = d_{i\epsilon_p} - d_{\mathbf{0}}$ for all $p,i$. We may call $\mathbf{m}^{\mathbf{d}}$ the \textit{dimension type} or \textit{flag type} of $(r,\mathbf{d})$. By \eqref{D:flagtype}, the set $\mathcal{F}(r)$ of flag types $\mathbf{m}$ of rank $r$ is finite. For any two flag types $\mathbf{m},\mathbf{m}' \in \mathcal{F}(r)$ of rank $r$, we set
\begin{equation}\label{E:flagbilinearform}
\langle \mathbf{m},\mathbf{m}' \rangle := \langle (r,0,\mathbf{m}) , (r,0,\mathbf{m}') \rangle_{par} - (1-g)r^2.
\end{equation}
We define an ordering on $\mathcal{F}(r)$ via \eqref{E:partialorder}:
\begin{equation}\label{E:flagorder}
\mathbf{m} \leq \mathbf{m}' \iff (r,0,\mathbf{m}) \leq (r,0,\mathbf{m}').
\end{equation}
In addition, for any flag type $\mathbf{m} \in \mathcal{F}(1)$ of rank $1$, i.e., $\mathbf{m} = \mathbf{m}^{\mathbf{d}}$ for some $(1,\mathbf{d}) \in \mathcal{K}^{+,vec}_1$, we set $m_p = \sum_{1 \leq i \leq w_p-1} m_{p,i}$ for all $p \in S$ and $\underline{m} = (m_p)_{p \in S}$. 

\subsection{Harder-Narasimhan filtration}\label{S3:HN} Finally, for any fixed slope $\nu \in \mathbb{Q} \cup \{ \infty \}$ we let $\mathcal{C}^\nu$ stand for the full subcategory of $\PCoh(X)$ consisting of semi-stable parabolic coherent sheaves of slope $\nu$. As an example, $\mathcal{C}^\infty = \PTor(X)$. The fundamental properties of the categories $\mathcal{C}^\nu$ are listed below:
\begin{prop}\label{P:HNfil}
The following hold:
\begin{enumerate}
\item[(i)] The categories $\mathcal{C}^\nu$ are abelian, artinian and noetherian.
\item[(ii)] $\PHom(\mathscr{F}_\bullet, \mathscr{G}_\bullet) = \emptyset$ if $\mathscr{F}_\bullet \in \mathcal{C}^\nu, \mathscr{G}_\bullet \in \mathcal{C}^\mu $ and $\nu > \mu$.
\item[(iii)] Any parabolic coherent sheaf $\mathscr{F}_\bullet$ admits a unique filtration
\begin{equation}\label{E:HNfil}
0 \subsetneq \mathscr{F}^s_\bullet \subsetneq \cdots \subsetneq \mathscr{F}^1_\bullet = \mathscr{F}_\bullet
\end{equation} 
satisfying the following conditions: $\mathscr{F}^i_\bullet / \mathscr{F}^{i+1}_\bullet$ is semi-stable for all $i$ and
\begin{equation*}
\mu(\mathscr{F}^1_\bullet / \mathscr{F}^2_\bullet ) < \cdots < \mu( \mathscr{F}^{s-1}_\bullet / \mathscr{F}^s_\bullet )< \mu( \mathscr{F}^s_\bullet).
\end{equation*} 
\item[(iv)] If $\mathscr{F}_\bullet$ and $\mathscr{G}_\bullet$ are stable and $\Par \mathscr{F}_\bullet = \Par \mathscr{G}_\bullet$, then $\mathscr{F}_\bullet \simeq \mathscr{G}_\bullet$.
\item[(v)] Let $\mathscr{F}_\bullet$ be semi-stable. There exists a filtration of $\mathscr{F}_\bullet$
\begin{equation*}
0 \subsetneq \mathscr{F}^t_\bullet \subsetneq \cdots \subsetneq \mathscr{F}^1_\bullet = \mathscr{F}_\bullet
\end{equation*}
such that $\mathscr{F}^i_\bullet / \mathscr{F}^{i+1}_\bullet$ is stable and $\mu( \mathscr{F}^i_\bullet / \mathscr{F}^{i+1}_\bullet) = \mu( \mathscr{F}_\bullet )$ for all $i$.
\end{enumerate}
\end{prop}

\noindent
\textit{Proof.} The argument is standard, c.f. \cite[Proposition 11, Th\'eor\`eme 8 and 12, Chap. 3]{Se}. 
\qed

\vspace{.1in}

The filtration \eqref{E:HNfil} is called the \textit{Harder-Narasimhan} (or HN for short) filtration of $\mathscr{F}_\bullet$ and the factors $\mathscr{F}^1_\bullet, \dots, \mathscr{F}^s_\bullet$ are called the \textit{semi-stable factors} of $\mathscr{F}_\bullet$. We define the \textit{HN-type} of $\mathscr{F}_\bullet$ to be $HN(\mathscr{F}_\bullet) = (\alpha_1,\cdots,\alpha_s)$ with $\alpha_i = [\mathscr{F}^i_\bullet] - [\mathscr{F}^{i+1}_\bullet] \in \mathcal{K}^+$. Note that we have $\alpha = \alpha_1 + \cdots \alpha_s = [\mathscr{F}_\bullet]$. We may associate a convex polygon $HNP(\mathscr{F}_\bullet)$ (called \textit{Harder-Narasimhan polygon} of $\mathscr{F}_\bullet$ in $\mathbb{R}^2$): the vertices of $HNP(\mathscr{F}_\bullet)$ are points $(0,0), (\rank \alpha_1, \Par \alpha_1), (\rank (\alpha_1+\alpha_2) , \Par (\alpha_1+\alpha_2) ), \dots, (\rank \alpha, \Par \alpha)$. The polygons corresponding to parabolic coherent sheaves have a natural partial ordering, namely, $HNP(\mathscr{F}_\bullet) \leq HNP(\mathscr{G}_\bullet)$ if all the vertices of $HNP(\mathscr{F}_\bullet)$ lie on or below the polygon $HNP(\mathscr{G}_\bullet)$. 

Let $\nu \leq \nu' \in \mathbb{Q} \cup \{ \infty \}$ and $\PCoh^{[\nu,\nu']}(X)$ be the full subcategory of $\PCoh(X)$ consisting of parabolic coherent sheaves all of whose semi-stable factors are of slope comprised between $\nu$ and $\nu'$. We define in a similar fashion the full subcategories $\PCoh^{\geq \nu}(X), \PCoh^{>\nu}(X), \PCoh^{\leq \nu}(X), \PCoh^{< \nu}(X) \dots$etc. 
These categories are preserved in a natural sense by extension of the base field. 

From the existence of a scheme of finite type parametrizing isomorphism classes of semi-stable parabolic sheaves of a given rank and parabolic degree (c.f. \cite[Th\'eor\`eme 23]{Se}), we deduce:
\begin{lem}\label{L:s-s}
For any $\alpha \in \mathcal{K}^+$, there exists only finitely many isomorphism classes of semi-stable parabolic sheaves of class $\alpha$. 
\end{lem}

\begin{lem}\label{L:semistab}
For any pair $\nu < \nu' \in \mathbb{Q} \cup \{ \infty \}$, we have the following:
\begin{enumerate}
\item[(i)] $\PCoh^{[\nu,\nu']}(X)$ is stable under extensions and taking quotients.
\item[(ii)] For any $\alpha \in \mathcal{K}^+$, the number of isomorphism classes of parabolic coherent sheaves in $\PCoh^{[\nu,\nu']}_\alpha(X)$ is finite.
\item[(iii)] Let $\alpha \in \mathcal{K}^+$ and $\mathscr{F}_\bullet \in \PCoh(X)$. The number of isomorphism classes of parabolic subsheaves $\mathscr{H}_\bullet \subset \mathscr{F}_\bullet$ of class $\alpha$ is finite. 
\end{enumerate}
\end{lem}

\noindent
\textit{Proof.} The fact that $\PCoh^{[\nu,\nu']}(X)$ is stable under extensions follows by definition. To prove $\PCoh^{[\nu,\nu']}(X)$ is stable under quotients, we fix $\mathscr{F}_\bullet \in \PCoh^{[\nu,\nu']}(X)$. Let $\mathscr{H}_\bullet$ be a quotient of $\mathscr{F}_\bullet$ and 
\begin{equation*}
0 \subsetneq \mathscr{H}^t_\bullet \subsetneq \cdots \subsetneq \mathscr{H}^1_\bullet = \mathscr{H}_\bullet
\end{equation*}
be the HN filtration of $\mathscr{H}_\bullet$. Then on one hand we have $\mathscr{H}_\bullet \in \PCoh^{\geq \mu(\mathscr{H}_\bullet / \mathscr{H}^2_\bullet) }(X)$. On the other hand, there exists a surjective map $\mathscr{F}_\bullet \to \mathscr{H}_\bullet / \mathscr{H}^2_\bullet$ and $\mathscr{H}_\bullet / \mathscr{H}^2_\bullet$ is semi-stable. This implies $\mu(\mathscr{H}_\bullet / \mathscr{H}^2_\bullet) \geq \nu$. Hence $\mathscr{H}_\bullet \in \PCoh^{\geq \nu}(X)$. Similarly, $\mathscr{H}_\bullet \in \PCoh^{\leq \mu(\mathscr{H}^t_\bullet)}(X)$. Since $\mathscr{H}^t_\bullet$ is a parabolic subsheaf of $\mathscr{F}_\bullet$ and $\mathscr{H}^t_\bullet$ is semi-stable, $\mathscr{H}^t_\bullet$ is a semi-stable factor of $\mathscr{F}_\bullet$. It implies $\mu(\mathscr{H}^t_\bullet \leq \nu'$ and hence $\mathscr{H}_\bullet \in \PCoh^{[\nu,\nu']}(X)$. 

Fix $\alpha \in \mathcal{K}^+$. There exists only finitely many HN-types $(\alpha_1,\dots,\alpha_s)$ satisfying $\sum_{i} \alpha_i = \alpha$, $\mu(\alpha_1) \geq \nu$ and $\mu(\alpha_s) \leq \nu'$. For any fixed HN-type $\underline{\alpha}$, there exist finitely many isomorphism classes of $\mathscr{F}_\bullet \in \PCoh(X)$ such that $HN(\mathscr{F}_\bullet) = \underline{\alpha}$. Indeed, such a parabolic sheaf $\mathscr{F}_\bullet$ is an iterated extension of semi-stable parabolic sheaves of respective classes $\alpha_1, \dots, \alpha_s$.  By Lemma \ref{L:s-s}, there are finitely many isomorphism classes of semi-stable parabolic sheaves of a given class (in $\mathcal{K}^+$) and all the $\PExt$ groups in the category $\PCoh(X)$ are finite. This proves the statement (ii).

For the statement (iii), let $\mathscr{F}_\bullet \in \PCoh(X)$. There exists $\nu \in \mathbb{Q} \cup \{\infty \}$ such that $\mathscr{F}_\bullet \in \PCoh^{\geq \nu}(X)$. Fix $\alpha \in \mathcal{K}^+$. Any quotient $\mathscr{H}_\bullet$ of $\mathscr{F}_\bullet$ belongs to $\PCoh^{\geq \nu}$ by statement (i). By statement (ii), there are only finitely many isomorphism classes of $\mathscr{H}_\bullet$ of class $[\mathscr{F}_\bullet] - \alpha$. Hence there are only finitely many isomorphism classes of parabolic subsheaves $\mathscr{G}_\bullet$ of $\mathscr{F}_\bullet$ of a given class $\alpha$. 

\qed

\section{The Hall algebras}\label{C:parHall}
\subsection{Reminder of Hall algebras}\label{S:Hall}
We briefly recall here the definition of the Hall algebra of a finitary abelian category for reader's convenience and refer to \cite[Lecture 1]{S1} for its standard properties. 
We will say an abelian category $\mathcal{A}$ \textit{finitary} if for any two objects $\mathscr{F}$ and $\mathscr{G}$ of $\mathcal{A}$ all the group $\Ext_{\mathcal{A}}^i(\mathscr{F}, \mathscr{G})$ have finite cardinality and are zero for almost all $i$. For instance any abelian category $\mathcal{A}$ which is linear over a finite field $\mathbf{k}$ and which satisfies $\dim \Ext_{\mathcal{A}}^i(\mathscr{F}, \mathscr{G}) < \infty$ for all $i$ is finitary. 
For a $\mathbf{k}$-linear finitary abelian category $\mathcal{A}$ we denote by $K(\mathcal{A})$ its Grothendieck group defined over $\mathbb{Z}$. If $\mathcal{A}$ is of finite global dimension then we may consider the \textit{Euler form}
\begin{equation*}
\langle \mathscr{F}, \mathscr{G} \rangle = \sum_{i=0}^{\infty} (-1)^i \dim \Ext_{\mathcal{A}}^i (\mathscr{F}, \mathscr{G}).
\end{equation*}
Note that the Euler form only depend on the classes of $\mathscr{F}$ and $\mathscr{G}$ in the Grothendieck group and thus it descends to a bilinear form $\langle \; , \; \rangle: K(\mathcal{A}) \otimes K(\mathcal{A}) \to \mathbb{C}$. We also need the \textit{symmetrized version} $( \mathscr{F}, \mathscr{G} ) = \langle \mathscr{F} , \mathscr{G} \rangle + \langle \mathscr{G} , \mathscr{F} \rangle$. 

Our only interest here is when $\mathcal{A}$ is finitary, abelian and hereditary (i.e. of global dimension at most one). Let $\mathcal{I}$ (resp. $\mathcal{I}_\alpha$) be the set of isomorphism classes of $\mathcal{A}$ (resp. of class $\alpha \in K(\mathcal{A})$). Let us choose a square root $v$ of $q^{-1}$ where $q$ is the number of elements of the finite field $\mathbf{k}$. The \textit{Hall algebra} $\mathbf{H}_{\mathcal{A}}$ of $\mathcal{A}$ is defined in the following way: as a vector space we have
\begin{equation*}
\begin{matrix}
\mathbf{H}_{\mathcal{A}}[\alpha] = \{ f: \mathcal{I}_\alpha \to \mathbb{C} \mid supp(f) < \infty \} = \bigoplus_{\mathscr{F} \in \mathcal{I}_\alpha} \mathbb{C} 1_{\mathscr{F}} , &
\mathbf{H}_{\mathcal{A}} = \bigoplus_{\alpha \in K(\mathcal{A})} \mathbf{H}_{\mathcal{A}}[\alpha],
\end{matrix}
\end{equation*}
where $1_{\mathscr{F}}$ stands for the characteristic function of $\mathscr{F} \in \mathcal{I}_\alpha$. The multiplication $m$ is defined as
\begin{equation}\label{E:multi}
m(f \otimes g) (\mathscr{R}) = (f \cdot g) (\mathscr{R}) = \sum_{\mathscr{N} \subseteq \mathscr{R}} v^{-\langle \mathscr{R} / \mathscr{N}, \mathscr{N} \rangle} f(\mathscr{R}/\mathscr{N}) g(\mathscr{N}).
\end{equation}
and the comultiplication $\Delta$ is defined as
\begin{equation}\label{E:comulti}
\Delta(f)(\mathscr{M},\mathscr{N}) = \frac{v^{\langle \mathscr{M},\mathscr{N} \rangle}}{\# \Ext_{\mathcal{A}}^1(\mathscr{M},\mathscr{N})} \sum_{\xi \in \Ext_{\mathcal{A}}^1(\mathscr{M}, \mathscr{N})} f(\mathscr{X}_\xi),
\end{equation}
where $\mathscr{X}_\xi$ is the extension of $\mathscr{M}$ by $\mathscr{N}$ corresponding to $\xi$. Note that the comultiplication may take values in a completion $\mathbf{H}_{\mathcal{A}} \widehat{\otimes} \mathbf{H}_{\mathcal{A}}$ of the tensor space $\mathbf{H}_{\mathcal{A}} \otimes \mathbf{H}_{\mathcal{A}}$. Then $(\mathbf{H}_{\mathcal{A}}, m)$ (resp. $(\mathbf{H}_{\mathcal{A}}, \Delta)$) is an associative algebra (resp. topological coassociative coalgebra). It is easy to deduce from the definition that, for any $\mathscr{M},\mathscr{N} \in \mathcal{I}$, we have 
\begin{equation}\label{E:multi2}
1_{\mathscr{M}} \cdot 1_{\mathscr{N}} = v^{-\langle \mathscr{M} ,\mathscr{N} \rangle} \sum_{\mathscr{L}} \frac{1}{a_{\mathscr{M}} a_{\mathscr{N}}} \mathbf{P}^{\mathscr{L}}_{\mathscr{M}, \mathscr{N}} 1_{\mathscr{L}}
\end{equation}
where $\mathbf{P}^{\mathscr{L}}_{\mathscr{M},\mathscr{N}} $ is the number of exact sequences 
\[
0 \to \mathscr{N} \to \mathscr{L} \to \mathscr{M} \to 0
\]
and $a_{\mathscr{M}} = \# \Aut(\mathscr{M})$. Note that
\[
\frac{1}{a_{\mathscr{M}} a_{\mathscr{N}}} = \# \{ \mathscr{R} \subset \mathscr{L} \mid \mathscr{R} \simeq \mathscr{N} \textit{ and } \mathscr{L}/\mathscr{R} \simeq \mathscr{M} \}.
\]
Similarly, for any $\mathscr{L} \in \mathcal{I}$ we have
\begin{equation}\label{E:comulti2}
\Delta(1_{\mathscr{L}}) = \sum_{\mathscr{M}, \mathscr{N}} v^{-\langle \mathscr{M} , \mathscr{N} \rangle } \frac{1}{a_{\mathscr{L}} } \mathbf{P}^{\mathscr{L}}_{\mathscr{M}, \mathscr{N}} 1_{\mathscr{M} } \otimes 1_{\mathscr{N}}.
\end{equation}
By definition the Hall algebra $\mathbf{H}_{\mathcal{A}}$ is gradded by the Grothendieck group $K(\mathcal{A})$. We will sometimes write $\Delta_{\alpha, \beta}$ in order to specify the graded components of the comultiplication. 

The Hall algebra $\mathbf{H}_{\mathcal{A}}$ will become a (topological) bialgebra if we suitably twist the coproduct: let $\mathbf{K}= \mathbb{C}[K(\mathcal{A})]$ be the group algebra of $K(\mathcal{A})$ and for any class $\alpha \in K(\mathcal{A})$ we denote by $\kappa_\alpha$ the corresponding element in $\mathbf{K}$. We define an \textit{extend Hall algebra} $\mathbb{H}_{\mathcal{A}} = \mathbf{H}_\mathcal{A} \otimes \mathbf{K}$ with relations
\begin{equation*}
\begin{matrix}
\kappa_\alpha \kappa_\beta = \kappa_{\alpha + \beta}, & \kappa_0 = 1, & \kappa_\alpha 1_{\mathscr{F}} \kappa_\alpha^{-1} = v^{- (\alpha, [\mathscr{F}])} 1_{\mathscr{F}},
\end{matrix}
\end{equation*}
where $[\mathscr{F}]$ stands for the class of $\mathscr{F}$ in $K(\mathcal{A})$. The new coproduct defined on $\mathbb{H}_{\mathcal{A}}$ is given by the formula
\begin{equation*}
\begin{matrix}
\widetilde{\Delta}(\kappa_\alpha) = \kappa_\alpha \otimes \kappa_\alpha,\\
\widetilde{\Delta}(f) = \sum_{\mathscr{M},\mathscr{N}} \Delta(f)(\mathscr{M}, \mathscr{N}) 1_{\mathscr{M}} \kappa_{[\mathscr{N}]} \otimes 1_{\mathscr{N}}.
\end{matrix}
\end{equation*}
Then we have the following Green's Theorem
\begin{theo}\label{T:Green}\cite{Gre}
Let $\mathcal{A}$ be a finitary abelian category of global dimension at most one. Then $(\mathbb{H}_{\mathcal{A}} , m, \widetilde{\Delta})$ is a (topological) biaglebra.
\end{theo}

Let
\begin{equation*}
( \quad , \quad )_G : \mathbb{H}_{\mathcal{A}} \otimes \mathbb{H}_{\mathcal{A}} \to \mathbb{C}
\end{equation*}
be the Green's Hermiltan scalar product defined by 
\begin{equation*}
(1_{\mathscr{M}} \kappa_\alpha , 1_{\mathscr{N}} \kappa_\beta)_G = \frac{\delta_{\mathscr{M},\mathscr{N}}}{a_{\mathscr{M}}} v^{-( \alpha, \beta)}. 
\end{equation*}
This scalar product is a \textit{Hopf pairing}, i.e., we have
\begin{equation}\label{E:hopf}
\begin{matrix}
(ab,c)_G = (a \otimes b, \widetilde{\Delta}(c) )_G ,& a,b,c \in \mathbb{H}_{\mathcal{A}}. 
\end{matrix}
\end{equation}
Moreover, the restriction of $( \quad , \quad  )_G$ to $\mathbf{H}_{\mathcal{A}}$ is nondegenerate. 

Finally, for any class $\gamma \in K(\mathcal{A})$, we set
\begin{equation*}
\mathbf{1}_\gamma = \sum_{\mathscr{M}, [\mathscr{M}] = \gamma} 1_\mathscr{M},
\end{equation*}
where the sum ranges over all the isoclasses of objects $\mathscr{M}$ of class $\gamma$. This sum may be infinite for some categories, so strickly speaking, $\mathbf{1}_\gamma$ belongs to the formal completion $\widehat{\mathbf{H}}_{\mathcal{A}}$ of $\mathbf{H}_{\mathcal{A}}$. The coproduct extends to a map 
\[
\Delta: \widehat{\mathbf{H}}_{\mathcal{A}} \to \mathbf{H}_{\mathcal{A}} \widehat{\otimes} \mathbf{H}_{\mathcal{A}}.
\]
We have the following useful Lemma
\begin{lem}\label{E:useful}\cite[Lemma 1.7]{S1}
\begin{equation*}
\Delta(\mathbf{1}_\gamma) = \sum_{\gamma = \alpha + \beta} v^{\langle \alpha , \beta \rangle } \mathbf{1}_\alpha \otimes \mathbf{1}_\beta.
\end{equation*}
\end{lem}

\subsection{Torsion subalgebra}\label{S4:halltor} 
We fix a smooth projective curve $X$ of genus $g$ over a finite field $\mathbf{k} = \mathbb{F}_q$ of $q$ elements, a finite set $S = \{p_1,\dots,p_N \} \subset X(\mathbf{k})$ of rational points of $X$, a collection $(w_p)_{p \in S}$ of positive integers and a collection $\underline{a} = (a_{p,i})_{p \in S, 1 \leq i \leq w_p-1}$ of real numbers satisfying (\ref{E:stab}). In this section we begin the study of the Hall algebra $\mathbf{H}:= \mathbf{H}_{\PCoh(X)}$ (resp. extended Hall algebra $\mathbb{H}:=\mathbb{H}_{\PCoh(X)}$) of the category $\PCoh(X)$. 

\vspace{.1in}

\noindent
\textit{Conventions:} We will always denote the same by $\mathcal{A}$ the set of isomorphism classes of objects in the category $\mathcal{A}$ for $\mathcal{A} = \PCoh(X),\PBun(X),\dots$etc to simplify the notations.  

The algebra $\mathbf{H}$ comes with a $K(\PCoh(X))$-grading, but we will mostly use the numerical $K$-group $\mathcal{K}$ only (c.f. \autoref{S3:Kgroup}):
\begin{equation*}
\mathbf{H} = \bigoplus_{(r,\mathbf{d}) \in \mathcal{K}} \mathbf{H}[(r,\mathbf{d})].
\end{equation*}
Of course, $\mathbf{H}[(r,\mathbf{d})] = \{ 0 \}$ unless $(r,\mathbf{d}) \in \mathcal{K}^+$. Let $\mathbf{K} = \mathbb{C}[\mathcal{K}]$ be the group algebra of $\mathcal{K}$.

Let $\mathbf{H}^{Tor} := \mathbf{H}_{\PTor(X)}$ be the Hall algebra of the subcategory $\PTor(X)$. Analogous to \eqref{E:torsiondec} there is a decomposition 
$
\mathbf{H}^{Tor} = \bigotimes'_{x \in |X|} \mathbf{H}_{\PTor_x(X)}
$
where $\mathbf{H}_{\PTor_x(X)}$ is the Hall algebra of $\PTor_x(X)$. Here by $\bigotimes'$ we mean the restricted tensor product, i.e., the linear space spanned by all vectors of the form $h_{x_1} \otimes \cdots \otimes h_{x_n}$, for $n \in \mathbb{N}, h_{x_i} \in \mathbf{H}_{\PTor_{x_i}(X)}$ and set of distinct points $\{x_1,\dots,x_n \}$. 

\vspace{.1in}

If $x \notin S$, since $\PTor_x(X) \simeq \mathscr{O}_{X,x}$-mod, $\mathbf{H}_{\PTor_x(X)}$ is isomorphic to the \emph{classical Hall algebra} defined over the residue field $\mathbf{k}_x$ at $x$, and is therefore identified with the algebra of symmetric functions:
\begin{equation*}
\begin{split}
\Psi_x: \mathbf{H}_{\PTor_x(X)} &\xrightarrow{\sim} \Lambda_{v_x} := \Lambda \otimes \mathbb{C}_{v_x} \\
1_{(\mathscr{O}^{(1^d)}_{X,x})_\bullet} &\mapsto v_x^{d(d-1)} e_d
\end{split}
\end{equation*}
where $v_x = v^{\deg(x)}$, $\mathbb{C}_{v_x} = \mathbb{C}[v_x, v_x^{-1}]$, $\Lambda = \mathbb{C}[y_1,y_2,\dots]^{\mathfrak{S}_\infty}$ is the Macdonald's ring of symmetric functions and $\{e_d\}_{d \in \mathbb{N}}$ denotes the elementary symmetric polynomials. Under this identification, for any partition $\lambda= (\lambda_1 \geq \cdots \geq \lambda_l) \in \Pi$, the element $1_{x,\lambda} :=1_{\mathscr{O}^{(\lambda)}_{x,\bullet}}$ corresponds to $v_x^{2n(\lambda)}P_\lambda(v_x^2)$, where $n(\lambda) = \sum_i (i-1)\lambda_i$ and $P_\lambda(t)$ is the Hall-Littlewood polynomial. We might set $I_{x,\lambda}=1_{\mathscr{O}_{x,\bullet}^{(\lambda)}}$. Let $d$ be any positive integer and $h_d, p_d \in \Lambda_{v_x}$ be the corresponding complete symmetric function and the power-sum symmetric function respectively, we define the elements $\mathbf{1}_{x,d}, T_{x,d} \in \mathbf{H}_{\PTor_x(X)}$ by
\begin{equation*}
\begin{matrix}
\mathbf{1}_{x,d} = \Psi^{-1}_x (h_d), & T_{x,d} = \Psi^{-1}_x (p_d).
\end{matrix}
\end{equation*}
Explicitly, they are given by the formula
\begin{equation*}
\begin{matrix}
\mathbf{1}_{x,d} = \sum_{ \begin{matrix} \mathscr{T}_\bullet \in \PTor_{x}(X) \\ [\mathscr{T}_\bullet ] = (0,d\delta) \end{matrix}} 1_{\mathscr{T}_\bullet} = \sum_{\lambda \vdash d} 1_{x,\lambda}, &T_{x,d} = \sum_{\lambda \vdash d} n_{v_x}(l(\lambda)-1) 1_{x,\lambda}
\end{matrix}
\end{equation*}
where $n_{v_x}(m) = \prod_{i=1}^m (1-v_x^{-2i})$. In fact, $\mathbf{1}_{x,d}$ and $T_{x,d}$ satisfy the relation
\begin{equation*}
\begin{matrix}
1 + \sum_{d >0 } \mathbf{1}_{x,d} s^{d \deg(x)} = \exp\big( \sum_{d >0} \frac{T_{x,d}}{d}s^{d \deg(x)} \big)
\end{matrix}
\end{equation*}
By construction, we have $\Delta(T_{x,d}) = T_{x,d} \otimes 1 + 1 \otimes T_{x,d}$ and the collection of elements $\{ T_{x,d} \mid d \in \mathbb{N} \}$ freely generates the commutative algebra $\mathbf{H}_{\PTor_x(X)}$. 

\vspace{.1in}

Let $p \in S$, $0 \leq i < w_p$ and $k \in \mathbb{N}^*$, we set 
\begin{equation*}
\begin{matrix}
I_{p,k}(i) : = 1_{\mathscr{O}^{(k)}_{p,\bullet}(i \epsilon_p)}, & I_{p,k}(-i) := I_{p,k}(w_p-i).
\end{matrix}
\end{equation*}
It is easy to see that
\begin{lem}\label{L:cyclicquiver} Let $p \in S$ and $ 0 \leq i \leq w_p-1$. Then
\begin{enumerate}
\item[(i)] For $n \in \mathbb{N}$ and $1 \leq j \leq w_p-1$, we have
\begin{equation}\label{indecomp}
I_{p,nw_p+j}(i) = \begin{cases}
I_{p,j}(i) & \textit{ if } n=0\\
I_{p,nw_p}(i) I_{p,j}(i) - v^2 I_{p,j}(i) I_{p,nw_p}(i) & \textit{ if } n > 0
\end{cases}.
\end{equation}
\item[(ii)] For $1 < j < w_p$, 
\begin{equation*}
I_{p,j}(i) = v^{-1}I_{p,j-1}(i) I_{p,1}(i-j+1) - I_{p,1}(i-j+1)I_{p,j-1}(i).
\end{equation*}
\end{enumerate}
\end{lem}

\noindent
\textit{Proof.} Recall from~\autoref{S3:tor} that $\PTor_p(X)$ is equivalent to $\Rep^{nil}_{\mathbf{k}}(C_{w_p})$. The proof is identical to the corresponding cyclic quiver, c.f. \cite[Lemma 5.16]{BS} for assertion i) and \cite[Sect. 3.3]{S1} for assertion ii).

\qed

Moreover, it follows from \cite{R} that the assignment $E_i \mapsto I_{p,1}(i)$ defines an embedding 
\begin{equation}\label{E:quantumgroup}
\mathbf{U}^+_v(\widehat{\mathfrak{sl}}_{w_p}) \hookrightarrow \mathbf{H}_{\PTor_p(X)}. 
\end{equation}
For any partition $\lambda = (\lambda_1 \geq \cdots \geq \lambda_r) \in \Pi$ and $0 \leq i < w_p$, we associate
\begin{equation*}
1_{p,\lambda}(i) :=1_{\bigoplus_j \mathscr{O}^{(\lambda_j w_p)}_{p,\bullet}(i)} \in \mathbf{H}_{\PTor_p(X)} 
\end{equation*}
and we set 
\begin{equation*}
\mathbf{H}^i_{\PTor_p(X)} :=  \bigoplus_{\lambda \in \Pi} \mathbb{C}_v 1_{p,\lambda}(i).
\end{equation*}
Then the assignment $1_{p,\lambda}(i) \mapsto v^{2n(\lambda)}P_\lambda(v^2)$ extends to an algebra isomorphism $\Psi_p(i): \mathbf{H}^i_{\PTor_p(X)} \xrightarrow{\sim} \Lambda_v$. Note that $1_{p,(d)}(i) = I_{p,d w_p}(i)$ and $\mathbf{H}^i_{\PTor_p(X)}$ is a commutative subalgebra (but not central) of $\mathbf{H}_{\PTor_p(X)}$ freely generated by the set $\{ I_{p,dw_p}(i) \}_{d \in \mathbb{N}^*}$. The multiplication maps $\mathbf{U}^+_v(\widehat{\mathfrak{sl}}_{w_p}) \otimes \mathbf{H}^0_{\PTor_p(X)} \to \mathbf{H}_{\PTor_p(X)}$ and $\mathbf{H}^0_{\PTor_p(X)} \otimes \mathbf{U}^+_v(\widehat{\mathfrak{sl}}_{w_p}) \to \mathbf{H}_{\PTor_p(X)}$ induce isomorphisms of $\mathbb{C}_v$-modules (c.f. \cite[Lemma 4.3]{S2}). 

\vspace{.1in}
\noindent
\textbf{Remarks.} However the embedding \eqref{E:quantumgroup} can be extended to an isomorphism 
\begin{equation*}
\mathbf{H}_{\PTor_p(X)} \xrightarrow{\sim} \mathbf{U}^+_v(\widehat{\mathfrak{sl}}_{w_p}) \otimes_{\mathbb{C}_v} \mathcal{Z}
\end{equation*}
where $\mathcal{Z} = \mathbb{C}_v[x_1,x_2,\dots]$ is a central subalgebra (see \cite[Theorem 4.2]{S2}). 

\vspace{.1in}
Similarly, for $d \geq 1$, we set $\mathbf{1}_{p,d}(i):= \Psi_p(i)^{-1}(h_d), T_{p,d}(i) := \Psi_p(i)^{-1}(p_d)$ . 

\subsection{Global Hecke algebras}\label{S4:hallbun} The full subcategory $\PBun(X)$ of parabolic vector bundles on $X$ is not abelian, but it is exact and stable under extensions. Hence it also gives rise to a subalgebra
\begin{equation*}
\mathbf{H}^{vec} := \bigoplus_{\mathscr{E}_\bullet \in \PBun(X)} \mathbb{C} 1_{\mathscr{E}_\bullet}.
\end{equation*}
Recall that every parabolic coherent sheaf $\mathscr{F}_\bullet$ is isomorphic to a (unique) direct sum $\mathscr{F}^{vec}_\bullet \oplus \mathscr{F}^{tor}_\bullet$ of a parabolic vector bundle $\mathscr{F}^{vec}_\bullet$ and a parabolic torsion sheaf $\mathscr{F}^{tor}_\bullet$. In addition, $\PHom(\mathscr{F}^{tor}_\bullet, \mathscr{F}^{vec}_\bullet) = \PExt(\mathscr{F}^{vec}_\bullet, \mathscr{F}^{tor}_\bullet) = 0$. It follows that
\begin{equation*}
1_{\mathscr{F}^{vec}_\bullet} \cdot 1_{\mathscr{F}^{tor}_\bullet} = v^{-\langle \mathscr{F}^{vec}_\bullet, \mathscr{F}^{tor}_\bullet \rangle_{par}} 1_{\mathscr{F}}.
\end{equation*}
We deduce that the multiplication map defines vector space isomorphisms
\begin{equation}\label{E:vsdecomp}
\begin{matrix}
m: \mathbf{H}^{vec} \otimes \mathbf{H}^{tor} \xrightarrow{\sim} \mathbf{H}, & \mathbf{H}^{vec} \otimes \mathbf{H}^{tor} \otimes \mathbf{K} \xrightarrow{\sim} \mathbb{H}.
\end{matrix}
\end{equation}
We denote by $\varpi: \mathbb{H} \to \mathbf{H}^{vec}$ the projection which satisfies
\begin{equation*}
\varpi(u_v u_t \kappa_{r,\mathbf{d}}) = \begin{cases}
u_v & \textit{ if } u_t =1 \\
0 & \textit{ otherwise}
\end{cases}, \quad u_v \in \mathbf{H}^{vec}, \quad u_t \in \mathbf{H}^{tor}, \quad (r,\mathbf{d}) \in \mathcal{K}.
\end{equation*}

In fact, one may view $\mathbf{H}$ as some kind of semi-direct product of $\mathbf{H}^{vec}$ and $\mathbf{H}^{tor}$. The collection of primitive elements $Prim= \{T_{x,d},T_{p,d}(0), I_{p,1}(i) \mid x \in |X| \setminus S, p \in S, d \in \mathbb{N}^*, 0 \leq i < w_p\}$ generates the algebra $\mathbf{H}^{tor}$. We have
\begin{prop}\label{P:hecke}
Let $x \in |X| \setminus S, p\in S, d \in \mathbb{N}^*$ and $0 \leq i < w_p$, we have
\begin{equation*}
\begin{matrix}
ad_v(T_{x,d}) \cdot \mathbf{H}^{vec} \subseteq \mathbf{H}^{vec}, & ad_v(T_{p,d}(0)) \cdot \mathbf{H}^{vec} \subseteq \mathbf{H}^{vec}, & ad_v(I_{p,1}(i)) \cdot \mathbf{H}^{vec} \subseteq \mathbf{H}^{vec}
\end{matrix},
\end{equation*}
where 
\[
ad_v(T)(1_{\mathscr{F}_\bullet}) := T \cdot 1_{\mathscr{F}_\bullet} - v^{([T],[\mathscr{E}_\bullet])_{par}} 1_{\mathscr{F}_\bullet} \cdot T,
\] 
$[T]$ denotes the corresponding class of $T$ in $\mathcal{K}^{+,tor}$ and $([T],[\mathscr{E}_\bullet])_{par}$ is the symmetrized Euler form.
\end{prop}

\noindent
\textit{Proof.} Let $\mathscr{E}_\bullet$ be a parabolic vector bundle over $X$ and $T \in Prim$. We will prove that $ad(T) \cdot 1_{\mathscr{E}_\bullet} \in \mathbf{H}^{vec}$ by showing that $( ad(T) \cdot 1_{\mathscr{E}_\bullet} , 1_{\mathscr{F}_\bullet} )_G = 0$ for all parabolic coherent sheaves $\mathscr{F}_\bullet \notin \PBun(X)$. Here $( \; , \; )_G$ is the Green's scalar product which is non-degenerate on $\mathbf{H}$ and a Hopf pairing.

Let $\mathscr{F}_\bullet$ be one such parabolic coherent sheaf of rank $r = \rank \mathscr{E}_\bullet$ and as usual we write $\mathscr{F}_\bullet = \mathscr{F}^{vec}_\bullet \oplus \mathscr{F}^{tor}_\bullet$. Then 
\begin{equation}\label{E:left}
\begin{split}
(T \cdot 1_{\mathscr{E}_\bullet} , 1_{\mathscr{F}_\bullet})_G &= \big(T \cdot 1_{\mathscr{E}_\bullet}, v^{\langle \mathscr{F}^{vec}_\bullet, \mathscr{F}^{tor}_\bullet \rangle_{par} } 1_{\mathscr{F}^{vec}_\bullet} \cdot 1_{\mathscr{F}^{tor}_\bullet} \big)_G \\
&= v^{\langle \mathscr{F}^{vec}_\bullet, \mathscr{F}^{tor}_\bullet \rangle_{par} } \big( \Delta(T) \star \Delta(1_{\mathscr{E}_\bullet}), 1_{\mathscr{F}^{vec}_\bullet} \otimes 1_{\mathscr{F}^{tor}_\bullet} \big)_G \\
&= v^{\langle \mathscr{F}^{vec}_\bullet, \mathscr{F}^{tor}_\bullet \rangle_{par} } \big( T \cdot 1_{\mathscr{E}_\bullet} \otimes 1 + 1_{\mathscr{E}_\bullet} \otimes T, 1_{\mathscr{F}^{vec}_\bullet} \otimes 1_{\mathscr{F}^{tor}_\bullet} \big)_G
\end{split}
\end{equation}
where $\star$ denotes the twisted multiplication on $\mathbf{H} \otimes \mathbf{H}$. The last equality holds since $T$ is primitive, i.e., $\Delta(T) = T \otimes 1 + 1 \otimes T$, for all $T \in Prim$ and 
\begin{equation*}
 \Delta_{(r_{\mathscr{E}_\bullet}, \mathbf{d}), (0,\mathbf{d}')} (1_{\mathscr{E}_\bullet}) = \begin{cases}
1_{\mathscr{E}_\bullet} \otimes 1 & \textit{ if } \mathbf{d} = \deg \mathscr{E}_\bullet \textit{ and } \mathbf{d}' = 0 \\
0 & \textit{ otherwise }
\end{cases}.
\end{equation*}
Similarly,
\begin{equation}\label{E:right}
v^{([T],[\mathscr{E}_\bullet])_{par}} (1_{\mathscr{E}_\bullet} \cdot T , 1_{\mathscr{F}_\bullet})_G = v^{\langle \mathscr{F}^{vec}_\bullet, \mathscr{F}^{tor}_\bullet \rangle_{par} } \big(  1_{\mathscr{E}_\bullet} \cdot T \otimes 1 + 1_{\mathscr{E}_\bullet} \otimes T, 1_{\mathscr{F}^{vec}_\bullet} \otimes 1_{\mathscr{F}^{tor}_\bullet} \big)_G
\end{equation}
By assumption that $\mathscr{F}_\bullet$ is not a parabolic vector bundle, $\deg \mathscr{F}^{tor}_\bullet \neq 0$. \eqref{E:left} or \eqref{E:right} does not vanish if and only if 
\[
\deg \mathscr{F}^{tor}_\bullet = \begin{cases}
d\deg(x)\delta & \textit{ if } T= T_{x,d} \\
d\delta & \textit{ if } T = T_{p,d}(0) \\
\alpha_{p,i} & \textit{ if } T = I_{p,1}(i)
\end{cases},
\]
and in these cases, \eqref{E:left} and \eqref{E:right} coincide. 

\qed

The Proposition \autoref{P:hecke} allows us to define an action
\begin{equation*}
\sigma: \mathbf{H}^{tor} \otimes \mathbf{H}^{vec} \to \mathbf{H}^{vec}
\end{equation*}
by simply extending the actions of generators linearly and multiplicatively.

\subsection{The completion}\label{S4:HN} We will now use the Harder-Narasimhan filtration (c.f. \autoref{S3:HN}) to define a completion of the Hall algebras $\mathbf{H}$ and $\mathbf{H} \otimes \mathbf{H}$. The set of all isomorphism classes of parabolic coherent sheaves on $X$ can be stratified by all possible HN-types. For an HN-type $\underline{\alpha} = (\alpha_1,\dots,\alpha_s)$, we denote by $\mathcal{S}_{\underline{\alpha}}$ the stratum of HN-type $\underline{\alpha}$ and set
\begin{equation*}
1_{\mathcal{S}_{\underline{\alpha}}} = \sum_{\mathscr{F}_\bullet \in \mathcal{S}_{\underline{\alpha}}} 1_{\mathscr{F}_\bullet}.
\end{equation*}
If $\underline{\alpha} = (\alpha)$ for $\alpha \in \mathcal{K}^+$, we simply set
\begin{equation*}
\mathbf{1}^{ss}_{\alpha} = 1_{\mathcal{S}_{(\alpha)}} = \sum_{\mathscr{F}_\bullet \in \PCoh^{\mu(\alpha)}_{\alpha}(X)} 1_{\mathscr{F}_\bullet}.
\end{equation*}
Recall from \autoref{S:Hall} that for any $\alpha \in \mathcal{K}^+$ we define $\mathbf{1}_{\alpha} = \sum_{\mathscr{F}_\bullet \in \PCoh_{\alpha}(X)} 1_{\mathscr{F}_\bullet}$. By the uniqueness of HN filtration (c.f. Proposition \autoref{P:HNfil}) of a given parabolic sheaf, we can easy deduce 
\begin{lem}\label{L:HNstra}
For any HN-type $\underline{\alpha} = (\alpha_1, \dots, \alpha_s)$, we have
\begin{equation*}
1_{\mathcal{S}_{\underline{\alpha}}} = v^{\sum_{i <j} \langle \alpha_i, \alpha_j \rangle_{par}} \mathbf{1}^{ss}_{\alpha_1} \cdots \mathbf{1}^{ss}_{\alpha_s}.
\end{equation*}
\end{lem}

For any pair $n \in \mathbb{Z} $, the Hall algebra 
\begin{equation*}
\mathbf{H}^{\geq n} := \bigoplus_{\mathscr{F}_\bullet \in \PCoh^{\geq n}(X)} \mathbb{C} \cdot 1_{\mathscr{F}_\bullet}
\end{equation*}
is naturally a subalgebra of $\mathbf{H}$ (but not a coalgebra). We set also
\begin{equation*}
\mathbf{H}^{\ngeq n} = \bigoplus_{\mathscr{F}_\bullet \notin \PCoh^{\geq n}(X)} \mathbb{C} \cdot 1_{\mathscr{F}_\bullet}.
\end{equation*}
Then we have $\mathbf{H} = \mathbf{H}^{\geq n} \oplus \mathbf{H}^{\ngeq n}$. For any $\alpha \in \mathcal{K}^+$, there is a surjective linear map of vector spaces $jet_n: \mathbf{H}[\alpha] \to \mathbf{H}^{\geq n}[\alpha] = \mathbf{H}^{\geq n} \cap \mathbf{H}[\alpha]$ which induces an isomorphism $\pi_n: \mathbf{H}[\alpha] / \mathbf{H}^{\ngeq n}[\alpha] \xrightarrow{\sim} \mathbf{H}^{\geq n}[\alpha]$. The canonical embedding $\mathbf{H}^{\ngeq m}[\alpha] \to \mathbf{H}^{\ngeq n}[\alpha]$ for any $m \leq n$ induces a commutative diagram
\begin{equation*}
\xymatrix{
\mathbf{H}[\alpha] / \mathbf{H}^{\ngeq m}[\alpha] \ar[r]^-{\pi_m} \ar[d] & \mathbf{H}^{\geq m} [\alpha] \ar[d]^-{\phi_{m,n}} \\
\mathbf{H}[\alpha] / \mathbf{H}^{\ngeq n}[\alpha] \ar[r]^-{\pi_n} & \mathbf{H}^{\geq n}[\alpha]
}.
\end{equation*}
Obviously $(\mathbf{H}^{\geq n}[\alpha], \phi_{m,n})$ forms a projective system and hence we can define
\begin{equation*}
\widehat{\mathbf{H}}[\alpha] := \varprojlim_{n} \mathbf{H}^{\geq n} [\alpha].
\end{equation*}
By Lemma \autoref{L:semistab} (ii) each $\mathbf{H}^{\geq n}[\alpha] $ is finite dimensional and since $\mathbf{H}[\alpha] = \cup_n \mathbf{H}^{\geq n}[\alpha]$, we may view $\widehat{\mathbf{H}}[\alpha]$ as the set of infinite sums $\sum_{\mathscr{F}_\bullet} u_{\mathscr{F}_\bullet} 1_{\mathscr{F}_\bullet}$ with $u_{\mathscr{F}_\bullet} \in \mathbb{C}, [\mathscr{F}_\bullet] = \alpha$, i.e., $\widehat{\mathbf{H}}[\alpha]  = \{ f: \PCoh_{\alpha}(X) \to \mathbb{C} \} = \prod_{\mathscr{F}_\bullet \in \PCoh_{\alpha}(X)} \mathbb{C} 1_{\mathscr{F}_\bullet}$ as a vector space. Finally, we define
\begin{equation*}
\widehat{\mathbf{H}} := \bigoplus_{\alpha \in \mathcal{K}^+} \widehat{\mathbf{H}}[\alpha].
\end{equation*}
Similarly, for $\alpha, \beta \in \mathcal{K}^+$, the sequence of vector spaces $(\mathbf{H}^{\geq n}[\alpha] \otimes \mathbf{H}^{\geq m}[\beta])$ forms a projective system and we put
\begin{equation*}
\widehat{\mathbf{H}}[\alpha] \otimes \widehat{\mathbf{H}}[\beta] := \varprojlim_{n,m} \mathbf{H}^{\geq n}[\alpha] \otimes \mathbf{H}^{\geq m}[\beta]
\end{equation*}
and
\begin{equation*}
(\widehat{\mathbf{H}} \widehat{\otimes} \widehat{\mathbf{H}})[\gamma] := \prod_{\begin{matrix}
\alpha,\beta \in \mathcal{K}^+\\
\alpha+\beta=\gamma
\end{matrix}} \widehat{\mathbf{H}}[\alpha] \otimes \widehat{\mathbf{H}}[\beta].
\end{equation*}
Finally, we define
\begin{equation*}
\widehat{\mathbf{H}} \widehat{\otimes} \widehat{ \mathbf{H}} := \bigoplus_{\gamma \in \mathcal{K}^+} ( \widehat{\mathbf{H}} \widehat{\otimes} \widehat{\mathbf{H}})[\gamma].
\end{equation*}

\begin{prop}\cite[Sect. 2]{BS2}
In the notation as above the following hold:
\begin{enumerate}
\item[(i)] $\widehat{\mathbf{H}}$ and $\widehat{\mathbf{H}} \widehat{\otimes} \widehat{\mathbf{H}}$ are associative algebras.
\item[(ii)] For any $\alpha , \beta \in \mathcal{K}^+$, we have $\Delta_{\alpha,\beta}( \widehat{\mathbf{H}}[\alpha+\beta] ) \subset \widehat{\mathbf{H}}[\alpha] \otimes \widehat{\mathbf{H}}[\beta]$. 
\end{enumerate}
\end{prop}

\section{The spherical Hall algebra}\label{C:SHA}
In this section, we introduce a natural subalgebra $\mathbf{U}$ of $\mathbf{H}$ which is much easier to study than the whole $\mathbf{H}$. It is analogous to the ``composition subalgebra" for Hall algebras of quivers. We will state some of its properties and give a shuffle presentation, valid for any arbitrary curve $X$ and any given parabolic datum $(S,\mathbf{w}, \underline{a})$. 

\subsection{Harder-Narasimhan strata}\label{S5:intro} For any $\alpha = (r,\mathbf{d}) \in \mathcal{K}^+$, we set
\begin{equation*}
\begin{matrix}
\mathbf{1}_{r,\mathbf{d}} = \sum_{\mathscr{F}_\bullet \in \PCoh_{r,\mathbf{d}}(X)} 1_{\mathscr{F}_\bullet} \in \widehat{\mathbf{H}}[(r,\mathbf{d})], & \mathbf{1}^{vec}_{r,\mathbf{d}} = \sum_{\mathscr{E}_\bullet \in \PBun_{r,\mathbf{d}}(X)} 1_{\mathscr{E}_\bullet} \in \widehat{\mathbf{H}}^{vec}[(r,\mathbf{d})]
\end{matrix}.
\end{equation*}
We define $\mathbf{U}^0$ (resp. $\mathbf{U}^>$) to be the subalgebra of $\mathbf{H}$ generated by $\mathbf{1}_{0,\mathbf{d}}$ for $(0,\mathbf{d}) \in \mathcal{K}^+$ (resp. by $\mathbf{1}^{vec}_{1,\mathbf{d}}$ for $(1,\mathbf{d}) \in \mathcal{K}^+$) and define $\mathbf{U}$ to be the subalgebra of $\mathbf{H}$ generated by $\mathbf{U}^0$ and $\mathbf{U}^>$. We call $\mathbf{U}$ the \emph{spherical Hall algebra} of $X$ with the fixed parabolic datum $(S,\mathbf{w},\underline{a})$. Note that $\mathbf{1}^{vec}_{r,\mathbf{d}} \in \mathbf{H}$ and $\mathbf{1}_{0,\mathbf{d}}$ are given by \textit{finite} sums, so that there are no completions involved in the definition of $\mathbf{U}$. However it makes sense to also consider the completion $\widehat{\mathbf{U}}$ of $\widehat{\mathbf{H}}$ defined as the closure of $\mathbf{U}$ in $\widehat{\mathbf{H}}$ with respect to the adic norm induced by the filtration $\{\mathbf{H}^{\geq n} \}_{n}$ of $\mathbf{H}$. We will first show that

\begin{theo}\label{T:HN}
For any $\alpha \in \mathcal{K}^+$, we have $\mathbf{1}^{ss}_\alpha \in \mathbf{U}$.
\end{theo}
Let us first prove the following:
\begin{prop}\label{P:T:HN}
For any $\alpha \in \mathcal{K}^+$, we have $\mathbf{1}^{ss}_\alpha \in \widehat{\mathbf{U}}$.
\end{prop}

\noindent
\textit{Proof.} By Reineke's inversion formula (see \cite{Rei}), 
\begin{equation*}
\mathbf{1}^{ss}_\alpha = \sum_{\underline{\beta}} (-1)^{s-1} v^{\sum_{i<j} \langle \beta_i,\beta_j \rangle_{par}} \mathbf{1}_{\beta_1} \cdots \mathbf{1}_{\beta_s}
\end{equation*}
where the sum ranges over all tuples $\underline{\beta}= (\beta_1,\dots,\beta_s)$ of elements in $\mathcal{K}^+$ such that $\mu (\sum_{t=k}^s \beta_t) > \mu (\alpha)$ for $k = 2, \dots, s$. The above sum converges in $\widehat{\mathbf{H}}$. Since $\widehat{\mathbf{U}}$ is a subalgebra of $\widehat{\mathbf{H}}$, the proposition will be proved if we can show $\mathbf{1}_\alpha \in \widehat{\mathbf{U}}$ for all $\alpha \in \mathcal{K}^+$. Furthermore, any $\mathscr{F}_\bullet \in \PCoh_\alpha(X)$ can be decomposed into a direct sum of a parabolic vector bundle and a parabolic torsion. Thus
\begin{equation}\label{E:lem}
\mathbf{1}_{r,\mathbf{d}} = \sum_{(0,\mathbf{d}') \in \mathcal{K}^{+,tor}} v^{ \langle (r,\mathbf{d} - \mathbf{d}'),(0,\mathbf{d}') \rangle_{par}} \mathbf{1}^{vec}_{r,\mathbf{d}-\mathbf{d}'} \mathbf{1}_{0,\mathbf{d}'},
\end{equation}
and $\mathbf{1}_{0,\mathbf{d}'} \in \mathbf{U}$ for all $(0,\mathbf{d}') \in \mathcal{K}^{+,tor}$ (\cite[Proposition 5.13]{BS}) and $\mathbf{1}^{vec}_{r,\mathbf{d}} = 0$ if $(r,\mathbf{d}) \notin \mathcal{K}^{+,vec}$, it suffices in fact to prove that $\mathbf{1}^{vec}_{r,\mathbf{d}} \in \widehat{\mathbf{U}}$ for all $(r,\mathbf{d}) \in \mathcal{K}^{+,vec}$. 

We will argue by induction on the rank $r$. The case $r=1$ is obvious by definition. Let $r >1$ and let us assume the proposition holds for all $r' < r$. We have to show that for any $(r,\mathbf{d}) \in \mathcal{K}^{+,vec}_r$ and any $n \in \mathbb{Z}$ we have
\begin{equation*}
\mathbf{1}^{vec}_{r,\mathbf{d}} \in \mathbf{U} + \mathbf{H}^{\ngeq n}.
\end{equation*}
Let us fix $n$ and argue by induction on $(r,\mathbf{d}) \in \mathcal{K}^{+,vec}_r$ with respect to the partial order \eqref{E:partialorder}. If $(r,\mathbf{d}) \in \mathcal{K}^{+,vec}_r$ satisfies $d_{\mathbf{0}} < nr - \sum_{p \in S} (w_p-1)r$, then any $\mathscr{E}_\bullet \in \PCoh_{r,\mathbf{d}}(X)$ has parabolic degree $< nr$. Hence $\mathscr{E}_\bullet \notin \PCoh^{\ngeq n}$, i.e., $\mathbf{1}_{r,\mathbf{d}} \in \mathbf{H}^{\ngeq n}$. Let us fix some $(r,\mathbf{d}) \in \mathcal{K}^{+,vec}_r$ and assume that \eqref{E:lem} holds for all $(r,\mathbf{d}' ) < (r,\mathbf{d})$. We choose $M < n -K$ where $K = \Par \omega_\bullet + \sum_{p \in S} (w_p -1)$ and let $(1,\mathbf{d}') \in \mathcal{K}^{+,vec}_1 $ be such that $d'_{\mathbf{0}} = M$ and $(r-1,\mathbf{d} - \mathbf{d}') \in \mathcal{K}^{+,vec}_{r-1}$. Consider the product
\begin{equation*}
\mathbf{1}_{r-1,\mathbf{d} - \mathbf{d}'} \cdot \mathbf{1}^{vec}_{1,\mathbf{d}'}= \sum_{\mathscr{F}_\bullet} c_{\mathscr{F}_\bullet} 1_{\mathscr{F}_\bullet}
\end{equation*} 
where 
\begin{equation*}
\begin{split}
c_{\mathscr{F}_\bullet} & = v^{-\langle (r-1,\mathbf{d} - \mathbf{d}') , (1,\mathbf{d}') \rangle_{par}} \sum_{\mathscr{L}_\bullet \in \PBun_{1,\mathbf{d}'}} \frac{\# \{ \mathscr{L}_\bullet \hookrightarrow \mathscr{F}_\bullet \} }{\# \Aut \mathscr{L}_\bullet} \\
&= v^{-\langle (r-1,\mathbf{d} - \mathbf{d}') , (1,\mathbf{d}') \rangle_{par}} \sum_{\mathscr{L}_\bullet \in \PBun_{1,\mathbf{d}'}} \frac{\# \{ \mathscr{L}_\bullet \hookrightarrow \mathscr{F}_\bullet \} }{v^{-2}-1}
\end{split}
\end{equation*}
Let us assume $\mathscr{F}_\bullet \in \PCoh^{\geq n}(X)$. Then $\mathscr{F}_\bullet \in \PCoh^{> M+K}(X)$. Since $\Par (\mathscr{L}_\bullet \otimes \omega_\bullet ) \leq M+K$, by the Serre duality \eqref{P:Serre} and Proposition \autoref{P:HNfil} (ii), we have $\PExt(\mathscr{L}_\bullet, \mathscr{F}_\bullet) = \emptyset$. Then
\begin{equation*}
\dim \PHom(\mathscr{L}_\bullet , \mathscr{F}_\bullet ) = \langle (1,\mathbf{d}'), (r,\mathbf{d}) \rangle_{par}.
\end{equation*}
Since any nonzero map from a parabolic line bundle into a parabolic vector bundle is an embedding, we deduce
\begin{equation*}
\# \{ \mathscr{L}_\bullet \hookrightarrow \mathscr{F}_\bullet \} = v^{-2 \dim \PHom(\mathscr{L}_\bullet, \mathscr{F}_\bullet)} - v^{-2 \dim \PHom(\mathscr{L}_\bullet, \mathscr{F}^{tor}_\bullet)}
\end{equation*}
and once the class of $\mathscr{L}_\bullet$ is fixed as $(1,\mathbf{d}')$, the above number depends only on the multi-degree of $\mathscr{F}^{tor}_\bullet$ by the Riemann-Roch \eqref{C:RR}. Hence there exist nonzero constants $c_{\mathbf{d}''}$ for $(0,\mathbf{d}'') \in \mathcal{K}^{+,tor}$ such that
\begin{equation*}
\mathbf{1}_{r-1,\mathbf{d}-\mathbf{d}'} \cdot \mathbf{1}^{vec}_{1,\mathbf{d}''} \in c_0 \mathbf{1}^{vec}_{r,\mathbf{d}} + \sum_{\mathbf{d}'' \in \mathcal{K}^{+,tor}} \mathbf{1}^{vec}_{r-1,\mathbf{d}-\mathbf{d}''} \cdot \mathbf{1}_{0,\mathbf{d}''} + \mathbf{H}^{\ngeq n}.
\end{equation*}
We can rewrite above as
\begin{equation*}
c_0\mathbf{1}^{vec}_{r,\mathbf{d}} \in \mathbf{1}_{r-1,\mathbf{d}-\mathbf{d}'} \cdot \mathbf{1}^{vec}_{1,\mathbf{d}'} - \sum c_{\mathbf{d}''} \mathbf{1}^{vec}_{r,\mathbf{d} - \mathbf{d}''} \cdot \mathbf{1}_{0,\mathbf{d}''} + \mathbf{H}^{\ngeq n}
\end{equation*}
and now by the induction hypothesis $\mathbf{1}_{r-1,\mathbf{d}-\mathbf{d}'} \in \widehat{\mathbf{U}}$ and $\mathbf{1}^{vec}_{r,\mathbf{d} - \mathbf{d}''} \in \widehat{\mathbf{U}}$. Hence \eqref{E:lem} holds. 

\qed

By Proposition \autoref{P:T:HN}, there exists for all $n$ an element $v_n \in \mathbf{H}^{\ngeq n}$ such that $u_n := \mathbf{1}^{ss}_\alpha + v_n \in \mathbf{U}$. We may further decompose $v_n = \sum_{\underline{\alpha}} v_{n,\underline{\alpha}}$ according to the HN type $\underline{\alpha}$. The set of $\underline{\alpha}$ for which $v_{n,\underline{\alpha}}$ is nonzero is a finite set since $v_n \in \mathbf{H}$. For those $\underline{\alpha}$ we have $\mu(\alpha_1) < n$. In order to show that $\mathbf{1}^{ss}_\alpha $ actually belongs to $\mathbf{U}$ and not only to $\widehat{\mathbf{U}}$, we need the following lemmas:

\begin{lem}\label{L:HN:1}
There exists $n \ll 0$ such that for any HN-type $\underline{\alpha} = (\alpha_1,\dots,\alpha_s)$ satisfying $\mu(\alpha_1) <n$, we have $\mu(\alpha_{i+1} )- \mu(\alpha_i) > \Par \omega_\bullet$ for some $1 \leq i \leq s$. 
\end{lem}

\noindent
\textit{Proof.} If $\Par \omega_\bullet \leq 0$, the statement is trivial. Suppose now $\Par \omega_\bullet >0$. Let $\underline{\alpha}$ be as above. We have $\Par (\alpha) = \rank(\alpha_1) \mu(\alpha_1) + \cdots \rank(\alpha_s) \mu(\alpha_s)$, where $\alpha = \alpha_1 + \cdots + \alpha_s$. If $\mu(\alpha_1) < n$ and $\mu(\alpha_{i+1}) - \mu (\alpha_i) < \Par \omega_\bullet$ for all $i$, then
\begin{equation*}
\begin{split}
\Par (\alpha) &< \rank(\alpha_1)n + \rank(\alpha_2)(n+ \Par \omega_\bullet) + \\
&+ \cdots + \rank(\alpha_s) \big( n + (s-1) \Par \omega_\bullet \big) \\
&= \rank(\alpha)n + \sum_{i=2}^{s} (i-1)\rank(\alpha_i)\Par \omega_\bullet \\
&< \rank(\alpha)\big( n + (\sum_{k=1}^{\rank(\alpha)} k ) \Par \omega_\bullet \big).
\end{split}
\end{equation*}
This is impossible for $n$ sufficiently negative.

\qed

\begin{lem}\label{L:HN:2}
Let $\mathscr{F}_\bullet \in \PCoh_{\alpha}(X)$ be of HN type $(\alpha_1,\dots,\alpha_s)$. Assume that $\mu(\alpha_{i+1}) - \mu(\alpha_i) > \Par \omega_\bullet$ for some $i$. Then $1_{\mathscr{F}_\bullet} = m \circ \Delta_{\beta,\gamma}(1_{\mathscr{F}_\bullet})$ for $\beta = \alpha_1 + \cdots + \alpha_i$ and $\gamma = \alpha_{i+1} + \cdots + \alpha_s$.
\end{lem}

\noindent
\textit{Proof.} Let $\mathscr{F}^s_\bullet \subset \cdots \subset \mathscr{F}^1_\bullet = \mathscr{F}_\bullet$ be the HN filtration of $\mathscr{F}_\bullet$. Since $\mathscr{F}^{i+1}_\bullet \in \PCoh^{\geq \mu(\alpha_{i+1})}(X)$ and $\mathscr{F}_\bullet / \mathscr{F}_\bullet^{i+1} \in \PCoh^{\leq \mu(\alpha_i)}(X)$ while $\mu(\alpha_{i+1}) - \mu(\alpha_i) > \Par \omega_\bullet$, we have
\begin{equation*}
\PExt(\mathscr{F}^{i+1}_\bullet, \mathscr{F}_\bullet / \mathscr{F}^{i+1}_\bullet) = \emptyset.
\end{equation*}
It follows that $\mathscr{F}_\bullet \simeq \mathscr{F}^{i+1}_\bullet \oplus \mathscr{F}_\bullet / \mathscr{F}^{i+1}_\bullet$, Moreover, $1_{\mathscr{F}_\bullet/ \mathscr{F}^{i+1}_\bullet} \cdot 1_{\mathscr{F}^{i+1}_\bullet} = v^{- \langle \mathscr{F}_\bullet / \mathscr{F}^{i+1}_\bullet, \mathscr{F}^{i+1}_\bullet \rangle_{par}} 1_{\mathscr{F}_\bullet}$ since there is a unique subsheaf of $\mathscr{F}_\bullet$ isomorphic to $\mathscr{F}^{i+1}_\bullet$. Hence the Lemma will be proved once we can show that $\Delta_{\beta,\gamma}(1_{\mathscr{F}_\bullet}) = v^{\langle \mathscr{F}_\bullet / \mathscr{F}^{i+1}_\bullet, \mathscr{F}^{i+1}_\bullet \rangle_{par}} 1_{\mathscr{F}_\bullet / \mathscr{F}^{i+1}_\bullet} \otimes 1_{\mathscr{F}^{i+1}_\bullet}$. However the last equation is a consequence of the same fact that there exists a unique subsheaf of $\mathscr{F}_\bullet$ of class $\gamma$, i.e., $\mathscr{F}^{i+1}_\bullet$.

\qed

\noindent
\textit{Proof of Theorem \autoref{T:HN}.} Let us choose $n \ll 0$ as in Lemma \autoref{L:HN:1}. Let $A$ be the finite set of all $\underline{\alpha}$ for which $v_{n,\underline{\alpha}}$ is nonzero and let $\underline{\alpha}^0$ be such that $HNP(\underline{\alpha}^0)$ is the lower boundary of the convex hull of elements of $A$. Thus $HNP(\underline{\alpha}^0) = ((0,0),(\rank(\alpha^0_1) \Par(\alpha^0_1)), \dots, ((\rank(\alpha^0_s), \Par(\alpha_s^0)) )$ is also a convex path in $\mathbb{R}^2$ of weight $\alpha$. Moreover $\mu(\alpha_1^0) < \infty$ so that the conclusion of Lemma \ref{L:HN:1} applies. Let $i$ be such that $\mu(\alpha^0_{i+1}) - \mu(\alpha^0_{i} > \Par \omega_\bullet$ and set $\beta = \alpha^0_1 + \cdots + \alpha^0_i$ and $\gamma = \alpha^0_{i+1} + \cdots + \alpha^0_s$. By Lemma \ref{L:HN:2}, $\Delta_{\beta,\gamma}(1_{\mathscr{F}_\bullet}) =0$ for all parabolic sheaves $\mathscr{F}_\bullet$ whose HN polygon does not lie below the segmant $((0,0),(\rank (\beta), \Par (\beta)))$. This implies $\Delta_{\beta,\gamma}(v_{n,\underline{\alpha}} = 0$ for all HN type $\underline{\alpha}$ whose associated HN polygon does not pass through the point $(\rank(\beta), \Par (\beta) )$. Furthermore, by Lemma \ref{L:HN:2} again, $m \circ \Delta_{\beta,\gamma}(v_{n,\underline{\alpha}}) = v_{n,\underline{\alpha}}$ for any HN type $\underline{\alpha}$ whose HN polygon does pass through $(\rank (\beta), \Par(\beta))$. Hence
\begin{equation*}
m \circ \Delta_{\beta,\gamma}(u_n) = m \circ \Delta_{\beta,\gamma} \big( \mathbf{1}^{ss}_\alpha + \sum_{\underline{\alpha}} v_{n,\underline{\alpha}} \big) = \sum_{\underline{\alpha} \in Z_\beta} v_{n,\underline{\alpha}}
\end{equation*} 
where $Z_\beta$ is the set of all HN types passing through $(\rank (\beta), \Par(\beta))$. Since $u_n \in \mathbf{U}$, which is stable under the coproduct, we deduce that $\sum_{\underline{\alpha} \in Z_\beta} v_{n,\underline{\alpha}} \in \mathbf{U}$. The same holds for $u'_n = \mathbf{1}^{ss}_\alpha + \sum_{\underline{\alpha} \notin Z_\beta} v_{n,\underline{\alpha}}$. Note that $u'_n$ contains strictly fewer terms that $u_n$. Arguing as above repeatly and the Theorem \autoref{T:HN} is proved. 

\qed

By Lemma \autoref{L:HNstra} and Theorem \autoref{T:HN}, we deduce

\begin{cor}\label{C:HN}
For any HN type $\underline{\alpha}$, we have $1_{\mathcal{S}_{\underline{\alpha}}} \in \mathbf{U}$.
\end{cor}

\subsection{Hecke actions}\label{S5:tor} The torsion part $\mathbf{U}^0$ of the spherical Hall algebra $\mathbf{U}$ has been studied in \cite{BS, S2}. The aim of this paragraph is to determine the Hecke action of $\mathbf{U}^0$ on the subspace of rank one functions $\mathbf{U}^>[1] := \mathbf{U}^> \cap \mathbf{H}[1]$, where $\mathbf{H}[1] := \bigoplus_{(1,\mathbf{d}) \in \mathcal{K}^{+,vec}_1} \mathbf{H}[(1,\mathbf{d})]$. For any sequence $\underline{i} = (i_p)_{p \in S}$ with $0 \leq i_{p} \leq w_{p}-1$ for all $p$, we consider the elements $T_{0,d}(\underline{i})$ defined by the generating series
\begin{equation*}
\exp \big( \sum_{d \geq 1} \frac{T_{0,d}(\underline{i})}{[d]_v} z^d \big) = \prod_{x \in |X| \setminus S} \exp \big( \sum_{n \geq 1} \frac{T_{x,n}}{n} z^{n\deg(x)} \big) \cdot \prod_{p \in S} \exp \big( \sum_{m \geq 1} \frac{T_{p,m}(i_p)}{m} z^m \big),
\end{equation*}
where $T_{x,d},T_{p,m}(i_p)$ are defined as in \autoref{S4:halltor} and 
\begin{equation*}
[s]_v := \frac{v^{-s}-v^s}{v^{-1}-v}
\end{equation*}
is the $q$-integer. We define also the elements $\mathbf{1}_{0,d}(\underline{i}),\theta_{0,d}(\underline{i})$ by
\begin{equation*}
\begin{split}
1+ \sum_{d \geq 1} \mathbf{1}_{0,d}(\underline{i}) z^d & = \exp \big( \sum_{m \geq 1} \frac{T_{0,m}(\underline{i})}{[m]_v} z^m \big) \\
1+ \sum_{d \geq 1} \theta_{0,d}(\underline{i}) z^d  &= \exp\big( (v^{-1} - v) \sum_{m \geq 1} T_{0,m}(\underline{i}) z^m \big)
\end{split}.
\end{equation*}

\begin{prop}\label{P:quantum}\cite[Proposition 5.13]{BS}
$\mathbf{U}^0$ is generated by the elements $\{ T_{0,d}(\mathbf{0}), I_{p,1}(i_p) \mid d \in \mathbb{N}^*, \; p \in S, \; 0 \leq i_p \leq w_p-1 \}$ .
\end{prop}

\begin{cor}\label{C:quamtum}
$\mathbf{U}^0 \cong \Lambda_v \otimes \mathbf{U}^+_v(\widehat{\mathfrak{sl}}_{w_{p_1}}) \otimes \cdots \otimes \mathbf{U}^+_v(\widehat{\mathfrak{sl}}_{w_{p_N}})$
\end{cor}

Alternatively, it will be useful to consider the following elements: for any subset $S' \subseteq S$ and $(i_p)_{p \in S \setminus S'} \in \prod_{p \in S \setminus S'} \mathbb{Z} /w_p\mathbb{Z}$, we define the elements $\mathbf{1}^{S'}_{0,d}, \theta^{S'}_{0,d}, T^{S'}_{0,d}$ by the generating series
\begin{equation*}
\begin{split}
\exp \big( \sum_{d \geq 1} \frac{T^{S'}_{0,d}(\underline{i})}{[d]_v} z^d \big) &= \prod_{x \in |X| \setminus S} \exp \big( \sum_{n \geq 1} \frac{T_{x,n}}{n} z^{n\deg(x)} \big) \cdot \prod_{p \in S \setminus S'} \exp \big( \sum_{m \geq 1} \frac{T_{p,m}(i_p)}{m} z^m \big), \\
1+ \sum_{d \geq 1} \mathbf{1}^{S'}_{0,d}(\underline{i}) z^d & = \exp \big( \sum_{m \geq 1} \frac{T^{S'}_{0,m}(\underline{i})}{[m]_v} z^m \big) \\
1+ \sum_{d \geq 1} \theta^{S'}_{0,d}(\underline{i}) z^d  &= \exp\big( (v^{-1} - v) \sum_{m \geq 1} T^{S'}_{0,m}(\underline{i}) z^m \big).
\end{split}
\end{equation*}

\begin{lem}\label{L4:TTT}
For any $S' \subseteq S$, $k \in \mathbb{N}$ and $d \in \mathbb{Z}$, one has
\begin{equation*}
[ \frac{T^{S'}_{0,k}(\underline{0})}{[k]_v}, \mathbf{1}^{vec}_{1,d\delta}]_v = v^{k}\frac{\# (X \setminus S')(\mathbb{F}_{q^k})}{k} \mathbf{1}^{vec}_{1,(d+k)\delta}
\end{equation*}
\end{lem}
\noindent
\textit{Proof.} Let us write
\begin{equation}\label{E4:TTT}
[ \frac{T^{S'}_{0,k}(\underline{0})}{[k]_v}, \mathbf{1}^{vec}_{1,d\delta}]_v = \sum_{\mathscr{L}_\bullet \in \PBun_{1,(d+k)\delta}(X)} c_{\mathscr{L}_\bullet} 1_{\mathscr{L}_\bullet}.
\end{equation}
By the definition of $T^{S'}_{0,k}(\underline{0})$ we see that
\begin{equation*}
c_{\mathscr{L}_\bullet} = v^k \sum_{d | k} \sum_{x \in X \setminus S'; \; \deg(x)=d} \frac{d}{k}\sum_{\lambda \vdash \frac{k}{d}} n_{v_x}(l(\lambda)-1) \frac{\# \PHom^{surj}(\mathscr{L}_\bullet, \mathscr{O}^{(\lambda)}_{x,\bullet}(\mathbf{0}))}{a_{\mathscr{O}^{(\lambda)}_{x,\bullet}(\mathbf{0})}}
\end{equation*}
It is easy to see that there is no epicmorphism $\mathscr{L}_\bullet \twoheadrightarrow \mathscr{O}^{(\lambda)}_{x,\bullet}(\mathbf{0})$ when $l(\lambda) >1$. Thus \eqref{E4:TTT} reduces to
\begin{equation*}
\begin{split}
c_{\mathscr{L}_\bullet} &= v^k \sum_{d | k} \sum_{x \in X \setminus S'; \; \deg(x)=d} \frac{d}{k}\frac{\# \PHom^{surj}(\mathscr{L}_\bullet, \mathscr{O}^{(\frac{k}{d})}_{x,\bullet}(\mathbf{0}))}{a_{\mathscr{O}^{(\frac{k}{d})}_{x,\bullet}(\mathbf{0})}} \\
&= v^k \sum_{d | k} \sum_{x \in X \setminus S'; \; \deg(x)=d} \frac{d}{k} \frac{v^{-2\frac{k}{d}} - v^{-2(\frac{k}{d}-1)}}{v^{-2\frac{k}{d}} - v^{-2(\frac{k}{d}-1)}}\\
&= \frac{v^k}{k} \sum_{d | k} \sum_{x \in X \setminus S'; \; \deg(x)=d} d \\
&= v^{k}\frac{\# (X \setminus S')(\mathbb{F}_{q^k})}{k}.
\end{split}
\end{equation*}
\qed

\begin{lem}\label{L5:simple}
For any $p \in S$, $0 \leq i < w_p$ and $ 0 < j < w_p$, we have the following relation:
\begin{equation*}
\varpi\big( I_{p,j}(i)  \mathbf{1}_{1,\mathbf{d}}^{vec} \big) = \begin{cases}
v\mathbf{1}^{vec}_{1,\mathbf{d} + \sum_{s=0}^{j-1} \alpha_{p,i+s} } & \textit{ if }  m_p^{\mathbf{d}} +j - i \in w_p\mathbb{Z} \\
0 & \textit{ otherwise } 
\end{cases}.
\end{equation*}
\end{lem}

\noindent
\textit{Proof.} It follows from Lemma \autoref{App:L:B}, see \autoref{App1}. \qed

\begin{lem}\label{L5:delta}
For any $(1,\mathbf{d}) \in \mathcal{K}^{+,vec}_1, p \in S, 0 \leq i < w_p$ such that $i \neq m^{\mathbf{d}}_p$ and any $n \geq 0$, we have
\begin{equation*}
\varpi \big( I_{p,n}(i) \mathbf{1}^{vec}_{1,\mathbf{d}} \big) = 0.
\end{equation*}
\end{lem}

\noindent
\textit{Proof.} By tensoring with a parabolic line bundle, we may assume that $\mathbf{d}=\mathbf{0}$. Then the Lemma is equivalent to show that there is no extension of $\mathscr{O}_\bullet$ by $\mathscr{O}_{p,\bullet}^{(nw_p)}(i\epsilon_p)$ belongs to $\PBun(X)$. Let $\mathscr{F}_\bullet$ be an extension of $\mathscr{O}_\bullet$ by $\mathscr{O}_{p,\bullet}^{(nw_p)}(i\epsilon_p)$ as:
\begin{equation*}
\xymatrix{
\cdots \ar[r] \ar[d] & \mathscr{O} \ar[r] \ar[d] & \mathscr{O} \ar[r] \ar[d] & \cdots \ar[r] \ar[d] & \mathscr{O} \ar[r] \ar[d] & \mathscr{O}(p) \ar[d] \\
\cdots \ar[r] \ar[d] & \mathscr{F}_{(j-1)\epsilon_p} \ar[r]^-{g_j} \ar[d] & \mathscr{F}_{j\epsilon_p} \ar[r] \ar[d] & \cdots \ar[r] \ar[d] & \mathscr{F}_{(w_p-1)\epsilon_p} \ar[r] \ar[d] & \mathscr{F}_{\mathbf{0}}(p) \ar[d] \\
\cdots \ar[r] & \mathscr{O}_p^{(n)} \ar[r]^-{f_j} & \mathscr{O}_p^{(n)} \ar[r] & \cdots \ar[r] & \mathscr{O}_p^{(n)} \ar[r] & \mathscr{O}_p^{(n)}
}
\end{equation*}
where $j = w_p -i$ and all the maps but $f_{j+sw_p}$ in the lower sequence are identity maps and $f_{j+sw_p}$ is nilpotent for all $s \in \mathbb{Z}$. Note that by assumption $i \neq 0$. Since $f_j$ is nilpotent, we have $\mathscr{F}_{(j-1)\epsilon_p} \ncong \mathscr{F}_{j\epsilon_p}$. But they have the same degree and rank. Hence at least one of them is not a vector bundle. It implies $\mathscr{F}_\bullet$ is not a parabolic vector bundle. \qed

\begin{cor}\label{C5:indecomp}
For any $p \in S$, $n \in \mathbb{N}$, $0 \leq i < w_p$ and $0 \leq r_p \leq w_p-1$, 
\begin{equation*}
\varpi \big( I_{p,nw_p+r_p}(i) \mathbf{1}^{vec}_{1,\mathbf{d}} \big) = \begin{cases} 
v^{n+(1-\delta_{r_p,0})} \mathbf{1}^{vec}_{1,\mathbf{d}+n\delta + \sum_{s=0}^{r_p-1} \alpha_{p,i+s}} & \textit{ if } m_p^{\mathbf{d}} +r_p - i \in w_p\mathbb{Z} \\
0 & \textit{ otherwise }
\end{cases}.
\end{equation*}
\end{cor}
\noindent
\textit{Proof.} 
Let us first suppose that $r_p = 0$. By Lemma \ref{L5:delta}, we may assume that $m_p^{\mathbf{d}} = i = 0$. Let us write
\begin{equation*}
\begin{split}
\varpi\big(I_{p,nw_p}(0)\mathbf{1}^{vec}_{1,\mathbf{d}}\big) = \sum_{\mathscr{L}_\bullet \in \PBun_{1,\mathbf{d}+n\delta}} c_{\mathscr{L}_\bullet} 1_{\mathscr{L}_\bullet}
\end{split}
\end{equation*}
and the coefficient $c_{\mathscr{L}_\bullet}$ is given by
\begin{equation*}
\begin{split}
c_{\mathscr{L}_\bullet} &= v^{-\langle (0,n\delta), (1,\mathbf{d}) \rangle_{par}} \frac{\# \PHom^{surj}(\mathscr{L}_\bullet, \mathscr{O}_{p,\bullet}^{(nw_p)}(0))}{a_{\mathscr{O}_{p,\bullet}^{(nw_p)}(0)}}
= v^n \frac{v^{-2n} -v^{-2(n-1)}}{v^{-2n}-v^{-2(n-1)}} =v^n
\end{split}.
\end{equation*}
Hence, for $i=m_p^{\mathbf{d}}$, we have
\begin{equation*}
\varpi\big(I_{p,nw_p}(i)\mathbf{1}^{vec}_{1,\mathbf{d}}\big) = v^n \mathbf{1}^{vec}_{1,\mathbf{d}+n\delta}.
\end{equation*}

Now let us assume that $r_p \neq 0$. By Lemma \autoref{L:cyclicquiver} (i), we have
\begin{equation}\label{C411:eq1}
\begin{split}
\varpi \big( I_{p,nw_p+r_p}(i) \mathbf{1}^{vec}_{1,\mathbf{d}} \big) &= \varpi \big( I_{p,nw_p}(i)I_{p,r_p}(i) \mathbf{1}^{vec}_{1,\mathbf{d}} - v^2I_{p,r_p}(i) I_{p,nw_p}(i) \mathbf{1}^{vec}_{1,\mathbf{d}}\big) \\
&= \varpi \big( I_{p,nw_p}(i)I_{p,r_p}(i) \mathbf{1}^{vec}_{1,\mathbf{d}} \big) - v^2 \varpi\big( I_{p,r_p}(i) I_{p,nw_p}(i) \mathbf{1}^{vec}_{1,\mathbf{d}}\big). 
\end{split}
\end{equation}
The second term of the last equality is always zero by Lemma \ref{L5:simple} and \ref{L5:delta} and the first one does not vanish only if the condition $m_p^{\mathbf{d}} +r_p - i \in w_p\mathbb{Z}$ is satisfied by Lemma \ref{L5:simple}. Hence \eqref{C411:eq1} becomes
\begin{equation*}
\varpi \big( I_{p,nw_p+r_p}(i) \mathbf{1}^{vec}_{1,\mathbf{d}} \big) = v \varpi\big( I_{p,nw_p}(i)\mathbf{1}^{vec}_{1,\mathbf{d} + \sum_{s=0}^{r_p-1} \alpha_{p,i+s}} \big) =v^{n+1} \mathbf{1}^{vec}_{1,\mathbf{d} + n\delta+\sum_{s=0}^{r_p-1} \alpha_{p,i+s}}. 
\end{equation*}
\qed

\noindent
\textbf{Remark.} Similar to the first part of the proof above, we have
\begin{equation*}
\varpi\big(1_{x,(n)} \mathbf{1}^{vec}_{1,\mathbf{d}}\big) = v^{n\deg(x)} \mathbf{1}^{vec}_{1,\mathbf{d}+n\deg(x)\delta}
\end{equation*}
for any $x \in |X| \setminus S$ and $n \in \mathbb{N}^*$.

Consequencely, it easy to see that the Hecke action of the subalgebra $\mathbf{U}^{+}_v(\widehat{\mathfrak{sl}}_{w_{p_1}}) \otimes \cdots \otimes \mathbf{U}^{+}_v(\widehat{\mathfrak{sl}}_{w_{p_N}})$ of $\mathbf{U}^0$ on $\mathbf{U}^>[1]$ is the restriction of the quantum affine algebra $\mathbf{U}_v(\widehat{\mathfrak{sl}}_{w_{p_1}}) \otimes \cdots \otimes \mathbf{U}_v(\widehat{\mathfrak{sl}}_{w_{p_N}})$'s action on the affine extremal weight module $V(\varpi_{p_1,1}) \otimes \cdots \otimes V(\varpi_{p_N,1})$ of extremal weight $(\varpi_{p_1,1}, \dots, \varpi_{p_N,1})$, where $\varpi_{p_i,1}$'s are the fundamental weights of the corresponding quantum affine algebras $\mathbf{U}_v(\widehat{\mathfrak{sl}}_{w_i})$ (c.f. \cite[Section 5]{Kash} for instance). Moreover, as an $\mathbf{U}_v(\widehat{\mathfrak{sl}}_{w_{p_1}}) \otimes \cdots \otimes \mathbf{U}_v(\widehat{\mathfrak{sl}}_{w_{p_N}})$-module, it is irreducible and admits a global base (\cite[Proposition 5.16]{Kash}). To be more precisely, for any $p \in S$, let us set $V_p = \bigoplus_{i=0}^{w_p-1} \mathbb{C}v_{p,i}$ and let $v_p^{(i)} = v_{p,0} + v_{p,1} + \cdots v_{p,i}$. Let
\begin{equation*}
V:= \mathbb{C}[z^\pm] \otimes \bigotimes_{p \in S} V_p. 
\end{equation*}
Identify $\mathbf{U}^>[1]$ with $V$ as vector spaces via the assignment
\begin{equation*}
\mathbf{1}^{vec}_{1,\mathbf{d}} \mapsto z^{\mathbf{d}_{\mathbf{0}}} \otimes v_{p_1}^{(m_{p_1}^{\mathbf{d}})} \otimes \cdots \otimes v_{p_N}^{(m_{p_N}^{\mathbf{d}})},
\end{equation*}
we may write
\begin{equation*}
z^d \otimes v_{p_1}^{(m_{p_1}+n_{p_1}w_{p_1})}\otimes \cdots  \otimes v_{p_N}^{(m_{p_N} + n_{p_N} w_{p_N})} = z^{d+\sum_{p \in S} n_p} \otimes v_{p_1}^{(m_{p_1})} \otimes \cdots \otimes v_{p_N}^{(m_{p_N})}.
\end{equation*}
Then the Hecke action of $\mathbf{U}^0$ on $\mathbf{U}^>[1]$ can be described as
\begin{equation}
\begin{split}
I_{p_i,1}(j) \bullet z^{d} v_{p_1}^{(m_{p_1})} \cdots v_{p_N}^{(m_{p_N})}  &= \begin{cases}
v z^{d} v_{p_1}^{(m_{p_1})} \cdots v_{p_i}^{(m_{p_i}+1)} \cdots v_{p_N}^{(m_{p_N})}  & \textit{ if } j \equiv m_{p_i} +1 \mod w_{p_i} \\
0 & \textit{ otherwise}
\end{cases}, \\
T_{0,d'}(\underline{i}) \bullet z^{d} v_{p_1}^{(m_{p_1})} \cdots v_{p_N}^{(m_{p_N})} &= \begin{cases}
\frac{v^{d'}\# X( \mathbb{F}_{q^{d'}}) [d']_{v}}{d'} z^{d+d'} v_{p_1}^{(m_{p_1})} \cdots v_{p_N}^{(m_{p_N})} & \textit{ if } i_p = m_p \quad \forall p \in S \\
0 & \textit{ otherwise }
\end{cases},
\end{split}
\end{equation}

\subsection{Spherical subalgebra of parabolic bundles}\label{S5:bun}
The rest of this section is devoted to investigate the structure of $\mathbf{U}^>$. Let us begin with a couple of lemmas. 

\begin{prop}\label{P5:theta}
\begin{equation*}
\Delta(\mathbf{1}^{vec}_{1,\mathbf{d}} )= \mathbf{1}^{vec}_{1,\mathbf{d}} \otimes 1 + \sum_{(0,\mathbf{d}') \in \mathcal{K}^{+,tor} } \theta^{\underline{m}^{\mathbf{d}}}_{0,\mathbf{d}'} \otimes \mathbf{1}^{vec}_{1,\mathbf{d}-\mathbf{d}'}.
\end{equation*}
where the non-vanished summands are $(0,\mathbf{d}')$ such that $\mathbf{d}' = d' \delta + \sum_{p \in S} \sum_{0 \leq i < r_p} \alpha_{p, m^{\mathbf{d}}_p + i}$ for $d' \geq 0$ and $0 \leq r_p < w_p$ for all $p \in S$ and
\begin{equation*}
\theta_{0,\mathbf{d}'}^{\underline{m}^{\mathbf{d}}} = \sum_{S' = S'_1 \sqcup S'_2} \bigg(   \prod_{p \in S'_2} \big(  -vI_{p,r_p}(m^{\mathbf{d}}_p) \big) \times \theta_{0,d'}(\underline{m}^{\mathbf{d}}) \times  \prod_{p' \in S'_1} \big( v^{-1}I_{p',r_{p'}}(m^{\mathbf{d}}_{p'}) \big) \bigg).
\end{equation*}
\end{prop}

\noindent
\textit{Proof.} We may assume first that $\mathbf{d} = \mathbf{0}$. Let $\mathscr{T}_\bullet \in \PTor(X)$ and $\mathscr{L}_\bullet \in \PPic(X)$. It is easy to see that
\begin{equation*}
\Delta(\mathbf{1}^{vec}_{1,\mathbf{0}} )(\mathscr{L}_\bullet, \mathscr{T}_\bullet) = \begin{cases}
1 & \textit{ if } \mathscr{T}_\bullet = 0 \textit{ and } \mathscr{L}_\bullet \in \PBun_{1,\mathbf{0}}(X) \\
0 & \textit{ otherwise}
\end{cases}.
\end{equation*}
While 
\begin{equation*}
\Delta(\mathbf{1}^{vec}_{1,\mathbf{0}})(\mathscr{T}_\bullet,\mathscr{L}_\bullet) = v^{-\langle [\mathscr{T}_\bullet],[\mathscr{L}_\bullet] \rangle_{par}} \sum_{\mathscr{L}'_\bullet \in \PBun_{1,\mathbf{0}}(X)} \mathbf{P}^{\mathscr{L}'_\bullet}_{\mathscr{T}_\bullet,\mathscr{L}_\bullet} \frac{a_{\mathscr{L}'_\bullet} a_{\mathscr{T}_\bullet}}{a_{\mathscr{L}_\bullet}}.
\end{equation*}
Note that for any $\mathscr{N} \in \Pic^0(X)$ we have $\mathscr{N} \otimes \mathscr{T}_\bullet \simeq \mathscr{T}_\bullet$. Hence
\begin{equation*}
\sum_{\mathscr{L}'_\bullet \in \PBun_{1,\mathbf{0}}(X)} \mathbf{P}^{\mathscr{L}'_\bullet}_{\mathscr{T}_\bullet, \mathscr{L}_\bullet} 
= \sum_{\mathscr{L}'_\bullet \in \PBun_{1,\mathbf{0}}(X)} \mathbf{P}^{\mathscr{L}'_\bullet \otimes \mathscr{N}}_{\mathscr{T}_\bullet, \mathscr{L}_\bullet \otimes \mathscr{N}} 
= \sum_{\mathscr{L}'_\bullet \in \PBun_{1,\mathbf{0}}(X)} \mathbf{P}^{\mathscr{L}'_\bullet}_{\mathscr{T}_\bullet, \mathscr{L}_\bullet \otimes \mathscr{N}}  
\end{equation*}
Thus $\Delta(\mathbf{1}^{vec}_{1,\mathbf{0}})(\mathscr{T}_\bullet, \mathscr{L}_\bullet)$ depends only on the class of $\mathscr{L}_\bullet$ in $\mathcal{K}^+$ and we may write
\begin{equation*}
\Delta(\mathbf{1}^{vec}_{1,\mathbf{0}}) = \mathbf{1}^{vec}_{1,\mathbf{0}} \otimes 1 + \sum_{\begin{matrix} (0,\mathbf{d}) \in \mathcal{K}^{+,tor} \\ (1,-\mathbf{d}) \in \mathcal{K}^{+,vec}_1 \end{matrix}} \theta_{0,\mathbf{d}} \otimes \mathbf{1}^{vec}_{1,-\mathbf{d}}.
\end{equation*}
Notice that in this case, we always have $\mathbf{d} = d \delta + \sum_{p \in S} \sum_{0 \leq i < r_p} \alpha_{p,i}$ for some $0 \leq r_p < w_p$. It remains to determine $\theta_{0,\mathbf{d}}$ explicitly. 

Let us fix $\mathbf{d} = d\delta + \sum_{p \in S} \sum_{0 \leq i <r_p} \alpha_{p,i} $ and let $\mathscr{L}'_\bullet \in \PBun_{1,\mathbf{0}}(X)$. Any corresponding torsion quotient $\mathscr{T}_\bullet$ of $\mathscr{L}_\bullet$ is of the form
\begin{equation}\label{E:L:tor}
\mathscr{T}_\bullet \simeq \bigoplus_{x \in |X| \setminus S} \mathscr{O}_{x,\bullet}^{(n_x)} \oplus \bigoplus_{p \in S} \mathscr{O}_{p,\bullet}^{(n_p w_p + r_p)}(\mathbf{0})
\end{equation}
for some nonnegative integers $(n_x)_{x \in |X| \setminus S}$ and $(n_p,r_p)_{p \in S}$ with $\sum_{x \in |X| \setminus S} n_x \deg(x) + \sum_{p \in S} n_p = d$ and $0 \leq r_p < w_p$ for all $p \in S$. We agree that $\mathscr{O}_{x,\bullet}^{(0)} = 0$. We set $S' = \{ p \in S \mid r_p \neq 0 \} \subset S$. For such $\mathscr{T}_\bullet$ as \eqref{E:L:tor}, $\mathbf{P}^{\mathscr{L}'_\bullet}_{\mathscr{T}_\bullet, \mathscr{L}_\bullet}$ is nonzero only if $\mathscr{L}_\bullet \in \PBun_{1,-\mathbf{d}}(X)$. Arguing as above, we have that $\sum_{\mathscr{L}_\bullet \in \PBun_{1,-\mathbf{d]}}(X)} \mathbf{P}^{\mathscr{L}'_\bullet}_{\mathscr{T}_\bullet, \mathscr{L}_\bullet}$ is independent of the choice of $\mathscr{L}'_\bullet \in \PBun_{1,\mathbf{0}}(X)$. But then
\begin{equation*}
\sum_{\mathscr{L}_\bullet \in \PBun_{1,-\mathbf{d}}(X)} \mathbf{P}^{\mathscr{L}'_\bullet}_{\mathscr{T}_\bullet,\mathscr{L}_\bullet} = \frac{1}{\# \Pic^0(X)} \sum_{\mathscr{L}_\bullet,\mathscr{L}'_\bullet} \mathbf{P}^{\mathscr{L}'_\bullet}_{\mathscr{T}_\bullet,\mathscr{L}_\bullet} = \sum_{\mathscr{L}'_\bullet} \mathbf{P}^{\mathscr{L}'_\bullet}_{\mathscr{T}_\bullet,\mathscr{L}_\bullet}.
\end{equation*}
Therefore
\begin{equation*}
\begin{split}
\Delta(\mathbf{1}^{vec}_{1,\mathbf{0}} )(\mathscr{T}_\bullet,\mathscr{L}_\bullet) 
&= v^{-\langle [\mathscr{T}_\bullet],[\mathscr{L}_\bullet] \rangle_{par}} \sum_{\mathscr{L}'_\bullet \in \PBun_{1,\mathbf{0}}(X)} \mathbf{P}^{\mathscr{L}'_\bullet}_{\mathscr{T}_\bullet,\mathscr{L}_\bullet} \frac{a_{\mathscr{L}_\bullet} a_{\mathscr{T}_\bullet}}{a_{\mathscr{L}'_\bullet}} \\
&= v^{-\langle [\mathscr{T}_\bullet],[\mathscr{L}_\bullet] \rangle_{par}}  \mathbf{P}^{\mathscr{O}_\bullet}_{\mathscr{T}_\bullet,\mathscr{L}_\bullet} a_{\mathscr{T}_\bullet} \\
&= v^{-\langle [\mathscr{T}_\bullet],[\mathscr{L}_\bullet] \rangle_{par}}  \# \PHom^{surj}(\mathscr{O}_\bullet, \mathscr{T}_\bullet).
\end{split}
\end{equation*}
Note that a morphism $\mathscr{O}_\bullet \to \mathscr{T}_\bullet$ is surjective if and only if all its components are surjectives. By the adjointness of $( - )_\bullet$ and $( - )_{\mathbf{0}}$, we have
\begin{equation*}
\# \PHom^{surj}(\mathscr{O}_\bullet, \mathscr{O}_{x,\bullet}^{(n_x)}) = \# \Hom^{surj}(\mathscr{O}_X, \mathscr{O}_{X,x}^{(n_x)}) = q^{\deg(x)n_x} - q^{\deg(x)(n_x-1)}
\end{equation*}
for $n_x \neq 0, x \in |X| \setminus (S \setminus S')$ and similarly for $p \in S'$
\begin{equation*}
\# \PHom ( \mathscr{O}_\bullet, \mathscr{O}_{p,\bullet}^{(n_p w_p + r_p)}(\mathbf{0})) = \# \Hom^{surj}(\mathscr{O}_X, \mathscr{O}_{X,p}^{(n_p+1)}) = q^{n_p+1} - q^{n_p}.
\end{equation*}
Moreover, $\langle [\mathscr{T}_\bullet],[\mathscr{L}_\bullet] \rangle_{par} = - \langle [\mathscr{L}'_\bullet],[\mathscr{T}_\bullet] \rangle_{par} = - \sum_{x \in |X|} n_x \deg(x)- \# S'$. We have
\begin{equation}\label{goodtheta}
\begin{split}
\theta_{0,\mathbf{d}}^{\underline{0}} &= v^{- \sum_{x \in |X|} n_x \deg(x) -\#S'} \sum_{\underline{n}} c_{\underline{n}} \prod_{x \in |X| \setminus S} 1_{x,(n_x)} \times \prod_{p \in S \setminus S'} I_{p,n_pw_p}(0) \times \prod_{p \in S'} I_{p,n_p w_p + r_p}(0)
\end{split}
\end{equation}
where $c_{\underline{n}} = \prod_{x \in |X|, n_x \neq 0} (1 - v^{2 \deg(x)})$ and the sum is taken over all decomposition types \eqref{E:L:tor} such that $\sum_{x \in |X|} n_x \deg(x) = d$. Using Lemma \ref{L:cyclicquiver} (i) and \eqref{goodtheta}, we may rewrite above as
\begin{equation*}
\begin{split}
\theta_{0,\mathbf{d}}^{\underline{0}} 
&= v^{-d - \#S'} \sum_{c_{\underline{n}}} c_{\underline{n}} \prod_{x \in |X| \setminus S} 1_{x,n_x} \times \prod_{y \in S \setminus S'} I_{y,n_yw_y}(0) \times \prod_{p \in S'} \big( I_{p,n_pw_p}(0) \cdot I_{p,r_p}(0) - v^2 I_{p,r_p}(0) I_{p,n_pw_p}(0) \big) \\
&= \sum_{S' = S'_1 \sqcup S'_2}  \prod_{p \in S'_2} \big( -v I_{p,r_p}(0) \big) \times \big( v^{-d} \sum_{\underline{n}} c_{\underline{n}} \prod_{x \in |X| \setminus S} 1_{x,n_x} \times \prod_{y \in S} I_{y,n_yw_y}(0) \big)  \times \prod_{p' \in S'_1} \big( v^{-1} I_{p',r_{p'}}(0) \big) \\
&= \sum_{S' = S'_1 \sqcup S'_2}  \bigg(   \prod_{p \in S'_2} \big(  -vI_{p,r_p}(0) \big) \times \theta_{0,d}(0) \times  \prod_{p' \in S'_1} \big( v^{-1}I_{p',r_{p'}}(0) \big) \bigg)
\end{split}
\end{equation*}

In general, for any $(1,\mathbf{d}) \in \mathcal{K}^{+,vec}_1$,  the possible $\mathbf{d}'$ such that $(0,\mathbf{d}') \in \mathcal{K}^{+,tor}, (1,\mathbf{d} - \mathbf{d}') \in \mathcal{K}^{+,vec}_1$ is of the form $\mathbf{d} = d' \delta + \sum_{p \in S} \sum_{0 \leq i < r_p} \alpha_{p, m^{\mathbf{d}}_p + i}$ for some $0 \leq r_p < w_p$ and in this case
\begin{equation}\label{E:theta}
\theta_{0,\mathbf{d}'}^{\underline{m}^{\mathbf{d}}} = \sum_{S' = S'_1 \sqcup S'_2} \bigg(   \prod_{p \in S'_2} \big(  -vI_{p,r_p}(m^{\mathbf{d}}_p) \big) \times \theta_{0,d'}(\underline{m}^{\mathbf{d}}) \times  \prod_{p' \in S'_1} \big( v^{-1}I_{p',r_{p'}}(m^{\mathbf{d}}_{p'}) \big) \bigg)
\end{equation}
\qed

 To compute the $\varpi\big( \theta_{0,\mathbf{d}'}^{\underline{i}} \mathbf{1}^{vec}_{1,\mathbf{d}} \big)$, using the expression \eqref{goodtheta} and by Corollary \autoref{C5:indecomp}, we have
\begin{equation}\label{E:theta1}
\begin{split}
\varpi\big(\theta_{0,\mathbf{d}'}^{\underline{i}}\mathbf{1}^{vec}_{1,\mathbf{d}}\big) &= v^{- \sum_{x \in |X|} n_x \deg(x) -\#S'} \sum_{\underline{n}} c_{\underline{n}} \varpi\bigg(\prod_{x \in |X| \setminus S} 1_{x,(n_x)} \times \prod_{p \in S \setminus S'} I_{p,n_pw_p}(i_p) \times \prod_{p \in S'} I_{p,n_p w_p + r'_p}(i_p) \mathbf{1}^{vec}_{1,\mathbf{d}} \bigg) \\
&= v^{- \sum_{p \in S} n_p -\#S'} \sum_{\underline{n}} c_{\underline{n}} \bigg( \prod_{p \in S \setminus S'} (\delta_{i_p,m_p^{\mathbf{d}}} v)^{n_p} \prod_{p \in S'} (\delta_{\overline{i_p-r'_p},m_p^{\mathbf{d}}} v)^{n_p+1} \bigg)\mathbf{1}^{vec}_{1,\mathbf{d}+\mathbf{d}'} \\
&= \sum_{\underline{n}} c_{\underline{n}} \prod_{p \in S \setminus S'} (\delta_{i_p,m_p^{\mathbf{d}}} )^{n_p} \prod_{p \in  S'} (\delta_{\overline{i_p-r'_p},m_p^{\mathbf{d}}} )^{n_p+1} \mathbf{1}^{vec}_{1,\mathbf{d}+\mathbf{d}'}  \\
\end{split}
\end{equation}
where the sums are taken over all $\underline{n}$ such that $\sum_{x \in |X|} n_x \deg(x) = d'$ and $c_{\underline{n}} = \prod_{x \in |X|, n_x \neq 0} (1 - v^{2\deg(x)})$ as in the proof of Proposition \autoref{P5:theta} and $\overline{i_p - r'_p}$ is known as the residue of $i_p-r'_p$ modulo $w_p$. Let us do a little bit more on \eqref{E:theta1}
\begin{equation}\label{E:theta2}
\begin{split}
\sum_{S'' \subseteq S \setminus S'} &\bigg( \sum_{\underline{n}; \; n_p = 0 \; \forall p \in S''} c_{\underline{n}} \prod_{p \in S \setminus (S' \sqcup S'')} (\delta_{i_p,m_p^{\mathbf{d}}} )^{n_p} \prod_{p \in  S' } \delta_{\overline{i_p-r'_p},m_p^{\mathbf{d}}}  \\ 
&+ \sum_{\underline{n}  ; \; n_p \neq 0 \; \forall p \in S'} c_{\underline{n}} \prod_{p \in S \setminus (S' \sqcup S'')} (\delta_{i_p,m_p^{\mathbf{d}}} )^{n_p}  \prod_{p \in S' \sqcup S''} \delta_{\overline{i_p-r'_p},m_p^{\mathbf{d}}}    \bigg)\mathbf{1}^{vec}_{1,\mathbf{d}+\mathbf{d}'} 
\end{split}
\end{equation}
and write
\begin{equation*}
\xi^{S''}_{d'} = \sum_{\underline{n}; \; n_p = 0, \; \forall p \in S'' } c_{\underline{n}}. 
\end{equation*}
Then we have
\begin{lem}\label{mainlemmatheta} 
For any $(1,\mathbf{d}) \in \mathcal{K}^{+,vec}_1$, $(0,\mathbf{d}') \in \mathcal{K}^{+,tor}$ with $\mathbf{d}' = d'\delta + \sum_{p \in S} \sum_{0 \leq i < r'_p} \alpha_{p,i_p+i}$ for some $\underline{i} = (i_p)_{p\in S} \in \prod_{p \in S} \mathbb{Z} / w_p\mathbb{Z}$, 
\begin{equation*}
\varpi\big( \theta_{0,\mathbf{d}'}^{\underline{i}} \mathbf{1}^{vec}_{1,\mathbf{d}} \big) = \begin{cases}
\xi^{S''}_{d''} \mathbf{1}^{vec}_{1,\mathbf{d}+\mathbf{d}'} & \textit{ if } (1,\mathbf{d} + \mathbf{d}') \in \mathcal{K}^{+,vec}_1 \\
0 & \textit{ otherwise}
\end{cases}
\end{equation*}
where $S'' = \{ p \in S \mid r'_p =0, \;  m_p^{\mathbf{d}} - i_p \notin w_p\mathbb{Z} \}$.
\end{lem}
\noindent
\textit{Proof.} The condition $(1,\mathbf{d}+\mathbf{d}' )\in \mathcal{K}^{+,vec}_1$ is equivalent to $m_p^{\mathbf{d}}+r'_p - i_p \in w_p \mathbb{Z}$ for all $p \in S$ such that $r'_p \neq 0$. If $S'' = \emptyset$, then the Lemma follows by \eqref{E:theta1}. Otherwise, in \eqref{E:theta2} the only summand does not vanish is the first term of $S'' = \{ p \in S \mid r'_p =0, \;  m_p^{\mathbf{d}} - i_p \notin w_p\mathbb{Z} \}$. Hence the Lemma follows.
\qed

\vspace{.1in}
\noindent
Moreover, we have 
\begin{lem}\label{L:zeta} As a series in $\mathbb{C}[[z]]$, we have
\begin{equation*} 
\sum_{d' \geq 0} \xi^{S''}_{d'} z^{d'}= \frac{\zeta_{X \setminus S''}(z)}{\zeta_{X\setminus S''}(v^2z)}
\end{equation*}
\end{lem}
\noindent
\textit{Proof.} It is easy to see that $\varpi\big(\theta^{S''}_{0,d'}(\underline{0}) \mathbf{1}^{vec}_{1,d\delta} \big)= \xi^{S''}_{d'} \mathbf{1}^{vec}_{1,(d+d')\delta}$. By Lemma \ref{L4:TTT} and the definition of $\theta^{S''}_{0,d'}$, we have
\begin{equation*}
\begin{split}
\sum_{d' \geq 0 } \xi^{S''}_{0,d'} z^{d'} & = \exp \bigg( (v^{-1}-v) \sum_{d \geq 1} \# (X \setminus S'') (\mathbb{F}_{q^d}) \frac{[d']_v v^{d'}}{d'} z^{d'} \bigg)\\
&= \exp\bigg( \sum_{d \geq 1} \# (X \setminus S'') (\mathbb{F}_{q^d}) (1-v^{2d'}) \frac{z^{d'}}{d'} \bigg) = \frac{\zeta_{X \setminus S''}(z)}{\zeta_{X \setminus S''}(v^2z)}.
\end{split}
\end{equation*}

\qed

\vspace{.1in}

As in \cite{SV2}, we introduce the so-called \emph{constant term map} as follows. For $r \geq 1$, let
\begin{equation*}
\begin{matrix}
\mathbf{H}[r] := \bigoplus_{(r,\mathbf{d}) \in \mathcal{K}^+} \mathbf{H}[(r,\mathbf{d})], & \mathbf{U}^>[r] := \mathbf{U} \cap \mathbf{H}[r]
\end{matrix}
\end{equation*}
and we set
\begin{equation*}
\begin{split}
J_r: \mathbf{U}^>[r] &\to \mathbf{U}^>[1] \widehat{\otimes} \cdots \widehat{\otimes} \mathbf{U}^>[1] \\
u &\mapsto \big( \varpi \otimes \cdots \otimes \varpi \big) \widetilde{\Delta}_{1,\dots,1}(u)
\end{split}
\end{equation*}
where $\widetilde{\Delta}_{1,\dots,1}$ is the iterated coproduct $ (\widetilde{\Delta}_{1,1} \otimes \Id \otimes \cdots \otimes \Id)  \circ \cdots \circ (\widetilde{\Delta}_{1,1} \otimes \Id ) \circ \widetilde{\Delta}_{1,1}$ of the component $\widetilde{\Delta}_{1,1} = \sum_{\mathbf{d},\mathbf{d}'} \widetilde{\Delta}_{(1,\mathbf{d}),(1,\mathbf{d}')}$ of the (extended) coproduct $\widetilde{\Delta}$. We denote by $J: \mathbf{U}^> \to \bigoplus_{r \geq 0} \big( \mathbf{U}^>[1] \big)^{\widehat{\otimes}r}$ the sum of the maps $J_r$. Writing 
\begin{equation*}
J(u) = \sum_{\mathscr{L}_\bullet^1,\dots,\mathscr{L}_\bullet^r} u\big(\mathscr{L}_\bullet^1,\dots,\mathscr{L}_\bullet^r \big) 1_{\mathscr{L}^1_\bullet} \otimes \cdots \otimes 1_{\mathscr{L}^r_\bullet}
\end{equation*}
we have
\begin{equation}\label{E:5:constant}
u\big(\mathscr{L}_\bullet^1,\dots,\mathscr{L}_\bullet^r \big)  = \frac{1}{(q-1)^r}\big( J_r(u), 1_{\mathscr{L}^1_\bullet} \otimes \cdots \otimes 1_{\mathscr{L}^r_\bullet} \big)_G = \frac{1}{(q-1)^r}\big( u, 1_{\mathscr{L}^1_\bullet}  \cdots  1_{\mathscr{L}^r_\bullet} \big)_G
\end{equation}
which coincides with the standard notion of constant term in the theory of automorphic forms. Since the image of $J_r$ lies in $(\mathbf{U}^>[1])^{\widehat{\otimes}r}$ and $\mathbf{U}^> [1] = \bigoplus_{(1,\mathbf{d}) \in \mathcal{K}^{+,vec}_1} \mathbf{1}^{vec}_{1,\mathbf{d}}$, the function $u\big(\mathscr{L}_\bullet^1,\dots,\mathscr{L}_\bullet^r \big)$ only depends on the respective multi-degree $\mathbf{d}_1,\dots,\mathbf{d}_r$ of the parabolic line bundles $\mathscr{L}^1_\bullet,\dots,\mathscr{L}^r_\bullet$. 

\begin{lem}
The constant term map $J: \mathbf{U}^> \to \bigoplus_{r} \big( \mathbf{U}^>[1] \big)^{\widehat{\otimes}r}$ is injective.
\end{lem}
\noindent
\textit{Proof.} Let $u \in \mathbf{U}^>[r]$ be in the kernel of $J_r$. By \eqref{E:5:constant} we have $(u,1_{\mathscr{L}^1_\bullet}  \cdots  1_{\mathscr{L}^r_\bullet} )_G =0$ for all $r$-tuple of parabolic line bundles $(\mathscr{L}^1_\bullet,\dots,\mathscr{L}^r_\bullet)$. In particular, $(u,\mathbf{1}^{vec}_{1,\mathbf{d}^1} \cdots \mathbf{1}^{vec}_{1,\mathbf{d}^r})_G = 0$ for all sequence $(\mathbf{d}^1,\dots,\mathbf{d}^r) \in ( \mathcal{K}^{+,vec}_1 )^r$. Since $\mathbf{U}^>$ is generated by the elements $\mathbf{1}^{vec}_{1,\mathbf{d}}$, we have $(u,\mathbf{U}^>)_G=0$. But the restriction of $( \;, \; )_G$ to $\mathbf{U}^>$ is non-degenerate. Therefore $u=0$. 
\qed
\vspace{.1in}

\subsection{Shuffle-like algebra}\label{S5:mainthm}

Let $\mathfrak{S}_r$ be the symmetric group on $r$ letters. For $w \in \mathfrak{S}_r$, $v = v_1 \otimes \cdots \otimes v_r \in V^{\otimes r}$ for some finite dimensional vector space $V$ and $P(z_1,\dots,z_r)$ a function in $r$ variables, we set $w(v) = (v_{w(1)} \otimes \cdots \otimes v_{w(r)})$ and $wP(z_1,\dots,z_r) = P(z_{w(1)} ,\dots, z_{w(r)})$. Let
\begin{equation*}
Sh_{r,s} = \{ w \in \mathfrak{S}_{r+s} \mid w(i) < w(j) \textit{ if } 1 \leq i < j \leq r \textit{ or } r < i < j \leq r+s \}
\end{equation*}
be the set of $(r,s)$-shuffles, i.e., the set of minimal length representatives of the left cosets in $\mathfrak{S}_{r+s} / \mathfrak{S}_r \times \mathfrak{S}_s$. For any $w \in Sh_{r,s}$, we write 
\begin{equation*}
I_w = \{ (i,j) \mid 1\leq i < j \leq r+s, \; w^{-1}(i) > r \geq w^{-1}(j) \}.
\end{equation*}

\vspace{.1in}

Recall from \autoref{S5:tor} that we may identify $\mathbf{U}^>[1]$ with $V= \mathbb{C}[z^{\pm 1}] \otimes \bigotimes_{p \in S} V_p$ via the assignment
\begin{equation*}
\mathbf{1}^{vec}_{1,\mathbf{d}} \mapsto z^{\mathbf{d}_{\mathbf{0}}} \otimes v_{p_1}^{(m^{\mathbf{d}}_{p_1})} \otimes \cdots \otimes v_{p_N}^{(m^{\mathbf{d}}_{p_N})}, \quad (1,\mathbf{d}) \in \mathcal{K}^{+,vec}_1.
\end{equation*}
Thus we have
\begin{equation*}
\big( \mathbf{U}^>[1] \big)^{\widehat{\otimes} r} \simeq \mathbb{C}[z_1^{\pm 1}, \dots, z_r^{\pm 1} ] [[ z_1/z_2, \dots, z_{r-1}/z_r ]] \otimes \bigotimes_{p \in S} V_p^{\otimes r} =: V^{\widehat{\otimes}r} .
\end{equation*}
Let $h(z) \in \mathbb{C}(z)$ be a fixed rational function. We define a linear oeprator $\Gamma_{h}(z/z') : V \otimes V \to V \widehat{\otimes} V$ by
\begin{equation}\label{S5:Gamma}
\begin{split}
&\Gamma_h(z/z')\big( (p(z)\otimes v_{p_1}^{(i_1)} \otimes \cdots \otimes v_{p_N}^{(i_N)} )  \otimes (q(z') \otimes v_{p_1}^{(j_1)} \otimes \cdots \otimes v_{p_N}^{(j_N)}\big) \\ 
= &\sum_{S' \subseteq S} h(z/z') \big(\frac{1-z/z'}{1-v^2z/z'} v\big)^{\# S'} ((p(z)\otimes v_{p_1}^{(i'_1)} \otimes \cdots \otimes v_{p_N}^{(i'_N)} )  \otimes (q(z') \otimes v_{p_1}^{(j'_1)} \otimes \cdots \otimes v_{p_N}^{(j'_N)}) 
\end{split}
\end{equation}
where, in each summand,
\begin{equation*}
\begin{matrix}
i'_p = \begin{cases}
i_p + \overline{j_p - i_p} & \textit{ if } p \in S \setminus S' \\
i_p & \textit{ if } p \notin S \setminus S'
\end{cases}, &
j'_p = \begin{cases}
j_p - \overline{j_p - i_p} & \textit{ if } p \in S \setminus S' \\
j_p & \textit{ if } p \notin S \setminus S'
\end{cases}
\end{matrix}
\end{equation*}
and the $\overline{j_p - i_p}$ is known as the residue of $j_p - i_p$ modulo $w_p$. 
 
We define an associative algebra $\mathbf{Sh}_{h(z)}(V)$ as follows: $\mathbf{Sh}_{h(z)}(V)$ is isomorphic, as a vector space, to $\bigoplus_{r \geq 0} V^{\widehat{\otimes} r}$.
The product is defined on $\mathbf{Sh}_{h(z)}(V)$ by:
\begin{equation*}
\begin{split}
P(z_1,\dots,z_r) \otimes v \star Q(z_1,\dots ,z_s) \otimes w =
 \sum_{\sigma \in Sh_{r,s}} \prod^{\longrightarrow}_{(i,j) \in I_\sigma} \Gamma^{ij}_h (z_i/z_j)\sigma\big( P(z_1,\dots,z_r) Q(z_{r+1},\dots,z_{r+s})  \otimes v \otimes w \big) 
\end{split}
\end{equation*}
where, $\Gamma_h^{i,j}(z_i/z_j)$ is known as acting on the pair of $i$th and $j$th compononts, the oriented product is with respect to the lexicographical order on the pairs $(i,j)$ and the rational functions $h(z_i/z_j)\big(\frac{1-z_i/z_j}{1-v^2z_i/z_j} v\big)^{\# S'} $ are developed as Laurent series in $z_1/z_2, \dots , z_{r-1}/z_r$. We equip also $\mathbf{Sh}_{h(z)}(V)$ with a coproduct $\Delta : \mathbf{Sh}_{h(z)}(V) \to \mathbf{Sh}_{h(z)}(V) \widehat{\otimes} \mathbf{Sh}_{h(z)}(V)$ defined by
\begin{equation*}
\begin{matrix}
\Delta_{s,t}([v_1 \cdots v_r]) =  [ v_1\cdots v_s] \otimes [v_{s+1} \cdots v_{s+t}], & \Delta = \bigoplus_{r=s+t} \Delta_{s,t}.
\end{matrix}
\end{equation*}
where we denote by $[v_1 \cdots v_s]$ the element $v_1 \otimes \cdots \otimes v_s$ in $V^{\widehat{\otimes}s}$. Let $\mathbf{S}_{h(z)} $ be the subalgebra of $\mathbf{Sh}_{h(z)}(V)$ generated by $V$.

\begin{theo}\label{T:maintheo} Set $h_X(z) = v^{2g-2} \frac{\zeta_X(z)}{\zeta_x(v^2z)}$. Then the constant term map $J: \mathbf{U}^> \to \mathbf{Sh}_{h_X(z)}(V)$ is an algebra morphism such that $J(\mathbf{U}^>) \simeq \mathbf{S}_{h_X(z)}$. Moreover we have
\begin{equation}\label{T:comulti}
(J_s \otimes J_t) \circ ( \varpi \otimes \varpi ) \circ \widetilde{\Delta}_{s,t} = \Delta_{s,t} \circ J_r, \quad r = s + t.
\end{equation}
\end{theo}
\noindent
\textit{Proof.} Since $\mathbf{U}^>$ is generated by $\mathbf{U}^>[1]$, it is enough to show that, for any tuple $((1,\mathbf{d}^1),\dots,(1,\mathbf{d}^r)) \in (\mathcal{K}^{+,vec}_1)^r$, we have
\begin{equation*}
J_r( \mathbf{1}^{vec}_{1,\mathbf{d}^1} \cdots \mathbf{1}^{vec}_{1,\mathbf{d}^r}) = J_1(\mathbf{1}^{vec}_{1,\mathbf{d}^1}) \star \cdots \star J_1(\mathbf{1}^{vec}_{1,\mathbf{d}^r}).
\end{equation*}
Let us first compute explicitly the $J_2$. Let $(1,\mathbf{d}^1)= (1,d_1,\mathbf{m}^1), (1,\mathbf{d}^2)=(1,d_2,\mathbf{m}^2) \in \mathcal{K}^{+,vec}_1$, we have
\begin{equation}\label{E5:E1}
\begin{split}
J_2( \mathbf{1}^{vec}_{1,\mathbf{d}^1} \mathbf{1}^{vec}_{1,\mathbf{d}^2}) = \varpi^{\otimes 2} \big( \widetilde{\Delta}_{1,1}(\mathbf{1}^{vec}_{1,\mathbf{d}^1} \mathbf{1}^{vec}_{1,\mathbf{d}^2} ) \big) = \varpi^{\otimes 2} \big( \sum_{\sigma \in \mathfrak{S}_2} \widetilde{\Delta}_{\delta_{\sigma(1)}} (\mathbf{1}^{vec}_{1,\mathbf{d}^1}) \widetilde{\Delta}_{\delta_{\sigma(2)}} (\mathbf{1}^{vec}_{1,\mathbf{d}^2}) \big)
\end{split}
\end{equation}
where we denote by $(\delta_1,\delta_2)$ the standard basis of $\mathbb{Z}^2$ and the second equality we have made use of that $\widetilde{\Delta}$ is a morphism of algebras. By the Proposition \ref{P5:theta}, we have
\begin{equation*}
\widetilde{\Delta}_{\delta_{k+1}}(\mathbf{1}^{vec}_{1,\mathbf{d}}) = \sum_{(0,\mathbf{d}_1),\dots,(0,\mathbf{d}_k) \in \mathcal{K}^{+,tor}} \bigg( \theta_{0,\mathbf{d}_1}^{\underline{m}^{\mathbf{d}}} \kappa_{1,\mathbf{d}-\mathbf{d}_1} \otimes \cdots \otimes \theta_{0,\mathbf{d}_k}^{\underline{m}^{\mathbf{d}-\sum_{j=1}^{k-1}\mathbf{d}_j}} \kappa_{1,\mathbf{d}-\sum_{j=1}^{k}\mathbf{d}_j} \otimes \mathbf{1}^{vec}_{1,\mathbf{d}-\sum_{j=1}^{k}\mathbf{d}_j} \otimes 1 \otimes \cdots \otimes 1 \bigg).
\end{equation*}
Thus \eqref{E5:E1} becomes
\begin{equation}\label{E5:E2}
\mathbf{1}_{1,\mathbf{d}^1}^{vec} \otimes \mathbf{1}_{1,\mathbf{d}^2}^{vec} + \varpi^{\otimes 2} \big( \sum_{(0,\mathbf{d}) \in \mathcal{K}^{+,tor}} \theta_{0,\mathbf{d}}^{\underline{m}^{\mathbf{d}_1}}\kappa_{1,\mathbf{d}^1-\mathbf{d}} \mathbf{1}^{vec}_{1,\mathbf{d}^2} \otimes \mathbf{1}^{vec}_{1,\mathbf{d}^1-\mathbf{d}} \big)
\end{equation}
The summand will not vanish under $\varpi^{\otimes 2}$ if and only if both $(1,\mathbf{d}^1-\mathbf{d})$ and $ (1,\mathbf{d}^2+\mathbf{d}) $ belong to $ \mathcal{K}^{+,vec}_1$, or equivalently, $\mathbf{d}$ is of the forms $\mathbf{d}_{\mathbf{0}}\delta + \sum_{p \in S} \sum_{0 \leq i < r_p} \alpha_{p,m^{\mathbf{d}^1}_p+i}$ with $\mathbf{d}_{\mathbf{0}} \geq 0$ and 
\begin{equation*}
r_p = \begin{cases}
m^{\mathbf{d}^1}_p - m^{\mathbf{d}^2}_p & \textit{ if } m^{\mathbf{d}^1}_p \geq m^{\mathbf{d}^2}_p \\
w_p +m^{\mathbf{d}^1}_p - m^{\mathbf{d}^2}_p & \textit{ if } m^{\mathbf{d}^1}_p < m^{\mathbf{d}^2}_p
\end{cases} \textit{ or } 0.
\end{equation*}
Therefore, for such a $\mathbf{d}$, let us write
\begin{equation*}
S^{\mathbf{d}^1,\mathbf{d}^2} = \{ p \in S \mid m_p^{\mathbf{d}^1} \neq m_p^{\mathbf{d}^2} \} \supseteq
S^{\mathbf{d}^1,\mathbf{d}^2}_{\mathbf{d}} = \{ p \in S^{\mathbf{d}^1,\mathbf{d}^2} \mid  r_p = 0 \}.
\end{equation*}
We have
\begin{equation*}
\begin{matrix}
(1,\mathbf{d}^1-\mathbf{d}) = (1,\mathbf{d}^1 - \mathbf{d}_{\mathbf{0}} \delta - \sum_{p \in S \setminus S^{\mathbf{d}^1,\mathbf{d}^2}_{\mathbf{d}}} \sum_{0 \leq i < r_p} \alpha_{p,m^{\mathbf{d}^1}_p+i} ), \\
(1,\mathbf{d}^2+\mathbf{d}) = (1,\mathbf{d}^2 + \mathbf{d}_{\mathbf{0}}\delta  + \sum_{p \in S \setminus S^{\mathbf{d}^1,\mathbf{d}^2}_{\mathbf{d}}} \sum_{0 \leq i < r_p} \alpha_{p,m^{\mathbf{d}^1}_p+i} )
\end{matrix}
\end{equation*}
and by Riemann-Roch
\begin{equation*}
( (1,\mathbf{d}^1-\mathbf{d}), (1,\mathbf{d}^2) )_{par} = 2(1-g) - \# S^{\mathbf{d}^1,\mathbf{d}^2}_{\mathbf{d}} .
\end{equation*}
Therefore the second term of \eqref{E5:E2} becomes
\begin{equation}\label{E5:E3}
v^{2g-2} \sum_{\mathbf{d}_{\mathbf{0}} \geq 0} \sum_{S' \subseteq S^{\mathbf{d}^1,\mathbf{d}^2}} v^{\# S'} \xi^{S'}_{\mathbf{d}_{\mathbf{0}}} \mathbf{1}^{vec}_{1,\mathbf{d}^2 + \mathbf{d}_{\mathbf{0}}\delta  + \sum_{p \in S \setminus S'} \sum_{0 \leq i < r_p} \alpha_{p,m^{\mathbf{d}^1}_p+i} } \otimes \mathbf{1}^{vec}_{1,\mathbf{d}^1 - \mathbf{d}_{\mathbf{0}} \delta - \sum_{p \in S \setminus S'} \sum_{0 \leq i < r_p} \alpha_{p,m^{\mathbf{d}^1}_p+i} }. 
\end{equation}
We introduce the automorphism $\gamma^d$ of $\mathbf{U}^>[1]$ for all $d \in \mathbb{Z}$ by $\gamma^d( \mathbf{1}^{vec}_{1,\mathbf{d}} ) = \mathbf{1}^{vec}_{1,\mathbf{d}+d\delta}$. Then we can rewrite \eqref{E5:E3} as
\begin{equation}\label{E5:E4}
\begin{split}
& v^{2g-2} \sum_{S' \subseteq S^{\mathbf{d}^1,\mathbf{d}^2}} \sum_{\mathbf{d}_0 \geq 0} \xi^{S'}_{\mathbf{d}_0}(\frac{\gamma_1}{\gamma_2})^{\mathbf{d}_0} v^{\# S'} \mathbf{1}^{vec}_{1,\mathbf{d}^2  + \sum_{p \in S \setminus S'} \sum_{0 \leq i < r_p} \alpha_{p,m^{\mathbf{d}^1}_p+i} } \otimes \mathbf{1}^{vec}_{1,\mathbf{d}^1 - \sum_{p \in S \setminus S'} \sum_{0 \leq i < r_p} \alpha_{p,m^{\mathbf{d}^1}_p+i} } \\
&=v^{2g-2} \sum_{S' \subseteq S^{\mathbf{d}^1,\mathbf{d}^2}} \frac{\zeta_{X \setminus S'}(\frac{\gamma_1}{\gamma_2})}{\zeta_{X \setminus S'}(v^2\frac{\gamma_1}{\gamma_2})} v^{\# S'} \mathbf{1}^{vec}_{1,\mathbf{d}^2  + \sum_{p \in S \setminus S'} \sum_{0 \leq i < r_p} \alpha_{p,m^{\mathbf{d}^1}_p+i} } \otimes \mathbf{1}^{vec}_{1,\mathbf{d}^1 - \sum_{p \in S \setminus S'} \sum_{0 \leq i < r_p} \alpha_{p,m^{\mathbf{d}^1}_p+i} }.
\end{split}
\end{equation}
Put
\begin{equation*}
h_{X \setminus S'}(z) = v^{2g-2} \frac{\zeta_{X \setminus S'}(z)}{\zeta_{X \setminus S'}(v^2z)},
\end{equation*}
and using $\zeta_{X \setminus S'}(z) = \zeta_X(z) (1-z)^{\#S'}$, we have
\begin{equation}
h_{X \setminus S'}(z) = h_X(z)\cdot \big( \frac{1-z}{1-v^2z} \big)^{\#S'}.
\end{equation}
Then \eqref{E5:E2} becomes
\begin{equation*}
\begin{split}
&\mathbf{1}^{vec}_{1,\mathbf{d}^1} \otimes \mathbf{1}^{vec}_{1,\mathbf{d}^2} +\sum_{S' \subseteq S^{\mathbf{d}^1,\mathbf{d}^2}} h_{X \setminus S'}(\frac{\gamma_1}{\gamma_2})v^{\# S'} \mathbf{1}^{vec}_{1,\mathbf{d}^2  + \sum_{p \in S \setminus S'} \sum_{0 \leq i < r_p} \alpha_{p,m^{\mathbf{d}^1}_p+i} } \otimes \mathbf{1}^{vec}_{1,\mathbf{d}^1 - \sum_{p \in S \setminus S'} \sum_{0 \leq i < r_p} \alpha_{p,m^{\mathbf{d}^1}_p+i} }.\\
& = \mathbf{1}^{vec}_{1,\mathbf{d}^1} \otimes \mathbf{1}^{vec}_{1,\mathbf{d}^2}  + \Gamma_{h_X}^{1,2}(\frac{\gamma_1}{\gamma_2}) \cdot \mathbf{1}^{vec}_{1,\mathbf{d}^2} \otimes \mathbf{1}^{vec}_{1,\mathbf{d}^1} 
\end{split}
\end{equation*}
We can easily deduce from the above computations that
\begin{equation*}
\begin{split}
J_r( \mathbf{1}^{vec}_{1,\mathbf{d}^1} \cdots \mathbf{1}^{vec}_{1,\mathbf{d}^r}) &= \sum_{\sigma \in \mathfrak{S}_r} \prod_{(i,j) \in I_\sigma}^{\longrightarrow}\Gamma_{h_X}^{i,j}(\frac{\gamma_i}{\gamma_j}) \mathbf{1}^{vec}_{1,\mathbf{d}^{\sigma(1)}} \otimes \cdots \otimes \mathbf{1}^{vec}_{1,\mathbf{d}^{\sigma(r)}} \\
&= J_1(\mathbf{1}^{vec}_{1,\mathbf{d}^1}) \star \cdots \star J_1(\mathbf{1}^{vec}_{1,\mathbf{d}^r}).
\end{split}
\end{equation*}
where the oriented product is with respect to the lexicographical order. 

To check the equation \eqref{T:comulti}, by coassociativity of the Hall comultiplication $\widetilde{\Delta}$, we have
\begin{equation*}
\begin{split}
\Delta_{s,t} \circ J_r &= (\varpi^{\otimes s} \otimes \varpi^{\otimes t}) \circ \widetilde{\Delta}_{(1^s),(1^t)} \circ \widetilde{\Delta}_{s,t} \\
&=  (\varpi^{\otimes s} \otimes \varpi^{\otimes t}) \circ \widetilde{\Delta}_{(1^s),(1^t)} \circ (\varpi \otimes \varpi ) \circ \widetilde{\Delta}_{s,t}  \\
&= (J_s \otimes J_t ) \circ (\varpi \otimes \varpi ) \circ \widetilde{\Delta}_{s,t} .
\end{split}
\end{equation*}

\qed

\subsection{The generic form}\label{S5:generic} Finally, for any fixed smooth projective curve $X$, observe that the spherical Hall algebra $\mathbf{U}_X := \mathbf{U}$ depends only on the genus $g$, the Weil numbers $\alpha_1, \dots, \alpha_{2g}$ of $X$ and the fixed parabolic datum $(S,\mathbf{w},\underline{a})$. Therefore it possesses a \textit{generic form} in the following sense. Let us fix $g \geq 0$ and consider the torus
\begin{equation*}
T_g = \{ ( \eta_1, \dots, \eta_{2g} ) \in ( \mathbb{C}^{\times} )^{2g} \mid \eta_{2i-1} \eta_{2i} = \eta_{2j-1} \eta_{2j} , \quad \forall i,j \}.
\end{equation*}
The Weyl group 
\begin{equation*}
W_g = \mathfrak{S}_g \rtimes (\mathfrak{S}_2)^{g}
\end{equation*}
naturally acts on $T_g$ and the collection $(\alpha_1,\dots,\alpha_{2g})$ defines a canonical element $\alpha_X$ in the quotient $T_g(\mathbb{C}) / W_g$. Let $R_g = \mathbb{Q}[T_g]^{W_g}$ and $K_g$ be its localization at the multiplicative set generated by $\{ q^s -1 \mid s \geq 1 \}$ where by definition $q (\alpha_1, \dots, \alpha_{2g}) = \alpha_{2i-1} \alpha_{2i}$ for any $1 \leq i \leq g$. For any choice of smooth projective curve $X$ of genus $g$ there is a natural map $K_g \to \mathbb{C}, \quad f \mapsto f(\alpha_X)$. We define, using the construction of $\mathbf{Sh}_{h(t)}(V)$, the $K_g$-algebra $\mathbf{U}^>_{K_g}$. Note that $\mathbf{U}^0$ has an obvious generic form $\mathbf{U}^0_{K_g}$ (c.f. \cite{Lu}). We can define also $\mathbf{U}_{K_g} = \mathbf{U}^>_{K_g} \otimes \mathbf{U}^0_{K_g}$. The (twisted) bialgebra structure and Green's scalar product both depend polynomially on the $\{ \alpha_1, \dots, \alpha_{2g} \}$ and hence may be defined over $K_g$. Let $\mathbf{U}^>_{R_g}$ be the $R_g$-subalgebra of $U^>_{K_g}$ generated by $R_g[z^{\pm 1}] \otimes \bigotimes_{p \in S} V_p \subset U^>_{K_g}$. By construction, $\mathbf{U}^>_{R_g}$ is a torsion-free integral form of $\mathbf{U}^>_{K_g}$ in the sense that $\mathbf{U}^>_{R_g} \otimes_{R_g} K_g = \mathbf{U}^>_{K_g}$. Moreover, there exists a specialization map
\begin{equation*}
\mathbf{U}^>_{R_g} \to \mathbf{U}^>_X, \quad 
\end{equation*}    
to a fixed curve $X$ of genus $g$. To summarize, Theorems \ref{T:HN} and \ref{T:maintheo} imply
\begin{cor}\label{S5:maincor}
There exists a (twisted) Hopf $R_g$-algebra $\mathbf{U}_{R_g}$ equipped with a Hopf pairing 
\begin{equation*}
( \; , \; )_G : \mathbf{U}_{R_g} \otimes \mathbf{U}_{R_g} \to K_g
\end{equation*}
generated by the elements $\mathbf{1}_{0,\mathbf{d},R_g},\mathbf{1}^{vec}_{1,\mathbf{e},R_g}$, $(0,\mathbf{d} ) \in \mathcal{K}^{+,tor}, (1,\mathbf{e}) \in \mathcal{K}^{+,vec}_1$, containing elements $1_{\mathcal{S}_{\underline{\alpha}},R_g}$ for any HN strata $\mathcal{S}_{\underline{\alpha}}$ and having the following property: for any smooth connected projective curve $X$ of genus $g$ defined over a finite field $\mathbf{k}$ there is a specialization morphism of (twisted) Hopf algebras 
\begin{equation*}
\Upsilon_X: \mathbf{U}_{R_g} \otimes_{R_g} \mathbb{C} \twoheadrightarrow \mathbf{U}_X
\end{equation*}
such that 
\begin{equation*}
\Upsilon_X ( 1_{\mathcal{S}_{\underline{\alpha}},R_g}) = 1_{\mathcal{S}_{\underline{\alpha}}}
\end{equation*}
for any HN stratum $\mathcal{S}_{\underline{\alpha}}$. 
\end{cor}

An immediate application of Corollary \autoref{S5:maincor} is the following theorem by Shiffmann:

\begin{cor}\cite[Theorem 7.1]{S5}
For any fixed genus $g$, any positive integer $N >0$, any collection of positive integers $\mathbf{w} = (w_1, \dots, w_N)$ and any tuple $\alpha= (r,\mathbf{d}) \in \mathcal{K}^{+,vec}_r$, there exists a unique polynomial $A_{g,\mathbf{w}, \alpha} \in K_g$ such that for any smooth projective curve $X$ of genus $g$ defined over a finite field, for any divisor $D = \sum_{i} w_i p_i$ with $p_i \in X(\mathbb{F}_q)$ we have 
\begin{equation}
A_{g,\mathbf{w},\alpha}(\sigma_X) = \mathcal{A}_{\alpha}(X),
\end{equation}
where $\mathcal{A}_{\alpha}(X)$ denotes the number of isomorphism classes of absolutely indecomposable vector bundles over $X$ of rank $r$ with $\mathbf{w}$-step quasi-parabolic structure over $p_1,\dots ,p_N$ and of multi-degree $\mathbf{d}$.
\end{cor}

\subsection{Remarks on counting stable parabolic Higgs bundles}

The proof of Schiffmann's Theorem goes through the same lines in \cite[Sect. 4.5]{S5}. One can go further to extend the Theroem in \cite{S5} to the parabolic cases by using Mozgovoy-Schiffmann's method in \cite{MSch}. Namely,
\begin{theo}\label{T:Higgs}
Let $\PHiggs^{ss,vec}_{r,\mathbf{d}}(X)$ be the moduli space of stable parabolic Higgs bundles of class $(r,\mathbf{d})$ over $X$ with $(r,\mathbf{d})$ generic. Then 
\begin{equation}
\mathcal{A}_{g,\mathbf{w},r,\mathbf{d}}(\sigma_X)(q) = v^{\dim \PHiggs^{ss,vec}_{r,\mathbf{d}}(X)} \# \PHiggs^{ss,vec}_{r,\mathbf{d}}(X)(\mathbb{F}_q).
\end{equation}
\end{theo}

However, a more general result has recently been provided in \cite[Theorem 1.4.5]{DGT}. We only sketch the approach following \cite{MSch} for our cases:

\begin{enumerate}
\item Since we feel more comfortable in an abelian category rather than exact one. Using the reformulation of \cite{GK}, we define a parabolic Higgs sheaf to be a pair $\overline{\mathscr{F}}=(\mathscr{F}, \theta)$ with $\mathscr{F} \in \PCoh(X)$ and $\theta \in \PHom(\mathscr{F}, \mathscr{F} \otimes \omega_\bullet)$. Let $\PHiggs(X)$ be the category of parabolic Higgs sheaves over $X$. Then the category $\PHiggs(X)$ is abelian and $2$-Calabi-Yau (the Serre duality can be found in \cite{Yo}). We can use the slope function $\mu(\overline{\mathscr{F}}) = \frac{\Par \mathscr{F}}{\rank \mathscr{F}}$ to define (semi-)stable objects in $\PHiggs(X)$.

\item Simple computation by using the Riemann-Roch (\ref{C:RR}) and long exact sequence in \cite[Theorem 4.1]{GK} implies that the Euler form $\chi$ on $\PHiggs(X)$ is symmetric. Therefore the quantum torus $\mathbb{T}=\mathbb{Q}(v^{-1})[K_0(\PHiggs(X))]$ equipped with the product
\begin{equation*}
e^{\alpha} \circ e^{\beta} = (-v^{-1})^{\chi(\alpha, \beta) - \chi(\beta,\alpha)} e^{\alpha+\beta} 
\end{equation*}
is then commutative. Moreover, the integration map $\int: \mathbf{H}_{\PHiggs(X)} \to \mathbb{T}$ by sending $1_{\overline{\mathscr{F}}} \mapsto (-v^{-1})^{\chi(\overline{\mathscr{F}},\overline{\mathscr{F}})}\frac{e^{[\overline{\mathscr{F}}]}}{\# \Aut \overline{\mathscr{F}}}$ will preserve the product of characteristic functions on semi-stable objects of decreasing slopes in the Hall algebra $\mathbf{H}_{\PHiggs(X)}$ inside the quantum torus $\mathbb{T}$. Therefore the uniqueness of Harder-Narasimhan filtration of a given parabolic Higgs bundle does imply the wall-crossing formula in this case. Note that the argument above will only make sense if we choose a suitable truncation using HN filtration as we did in \autoref{S5:intro}. 

\item The dimension reduction in our case follows with the same lines as in proof of \cite[Theorem 4.6]{MSch} and it shows that the number of points of the truncated substack is given by 
\[v^{-\dim \PHiggs^{ss,vec,\geq 0}_{r,\mathbf{d}}(X)}\mathcal{A}^{\geq 0}_{g,\mathbf{w},r,\mathbf{d}}(\sigma_X)(\mathbb{F}_q),\]
the number of absolutely indecomposable parabolic bundles of rank $r$ and multi-degree $\mathbf{d}$ on $X$ over $\mathbb{F}_q$. Twisting with a parabolic line bundle $\mathscr{L}$ of sufficient large multi-degree $\mathbf{d}'$, Lemma \ref{L:HN:2} implies that 
\[
\mathcal{A}^{\geq 0}_{g,\mathbf{w},r,\mathbf{d}+\mathbf{d}'}(\sigma_X)(\mathbb{F}_q) = \mathcal{A}_{g,\mathbf{w},r,\mathbf{d}+\mathbf{d}'}(\sigma_X)(\mathbb{F}_q).
\]
Finally, using $\PHiggs^{ss,vec}_{r,\mathbf{d}}(X) = \PHiggs^{ss,vec}_{r,\mathbf{d}+rd'\delta}(X)$, we have Theorem \autoref{T:Higgs}. 

\end{enumerate}

\noindent
\textbf{Remarks.}
\begin{enumerate}
\item Theorem \autoref{T:Higgs} may be generalized to more general twisted parabolic Higgs sheaves (in which the Euler form does not vanish) by using the Harder-Narasimhan recursion, c.f. \cite[Section 7]{MSch}.
\item It involves more complicated computations to obtain an explicit formula for the Kac polynomials $\mathcal{A}_{g,\mathbf{w},r,\mathbf{d}}$ of orbiford curves and it will be the subject in the companion paper.
\end{enumerate}

\appendix
\section{Proof of Lemma \ref{L5:simple}}\label{App1}
We prove the following Lemma insteaded: 

\begin{lem}\label{App:L:B} For any $p \in S$, $0\leq i < w_p$ and $0 < j < w_p$, we have
\begin{equation}\label{AppB:EQ1}
I_{p,j}(i) \mathbf{1}^{vec}_{1,\mathbf{0}} =
\begin{cases}
\mathbf{1}^{vec}_{1,\mathbf{0}} I_{p,j}(i) & \textit{ if } i > j \\
v^{-1} \mathbf{1}^{vec}_{1,\mathbf{0}} I_{p,j}(0) & \textit{ if } i=0 \\
v\mathbf{1}^{vec}_{1,\mathbf{0}}I_{p,j}(i) + v \mathbf{1}^{vec}_{1,\alpha_{p,1}+\cdots + \alpha_{p,i}} & \textit{ if } i=j\\
\mathbf{1}^{vec}_{1,\mathbf{0}} I_{p,j}(i) + (v^{-1} - v) \mathbf{1}^{vec}_{1,\alpha_{p,1} + \cdots + \alpha_{p,i}} I_{p,j-i}(0)& \textit{ if } 0 < i<j
\end{cases}.
\end{equation}
\end{lem}
\noindent
\textit{Proof.} Let us first suppose that $j=1$. If $i \neq 0,1$, we have 
\[
\begin{split}
\langle \mathscr{O}_\bullet, \mathscr{O}_{p,\bullet}^{(1)}(i\epsilon_p) \rangle_{par} &= \langle \mathscr{O}_{p,\bullet}^{(1)}(i\epsilon_p) , \mathscr{O}_\bullet \rangle_{par} =0 \\
\PHom(\mathscr{O}_\bullet, \mathscr{O}_{p,\bullet}^{(1)}(i\epsilon_p)) &= \PHom( \mathscr{O}_{p,\bullet}^{(1)}(i\epsilon_p),\mathscr{O}_\bullet) =0 \\
\PExt(\mathscr{O}_\bullet, \mathscr{O}_{p,\bullet}^{(1)}(i\epsilon_p)) &= \PExt( \mathscr{O}_{p,\bullet}^{(1)}(i\epsilon_p),\mathscr{O}_\bullet) =0 
\end{split}
\]
It is also true for any parabolic line bundle $\mathscr{L}_\bullet$ of multi-degree $\deg \mathscr{L}_\bullet = \mathbf{0}$ instead of $\mathscr{O}_\bullet$. Hence $I_{p,1}(i)$ and $\mathbf{1}^{vec}_{1,\mathbf{0}}$ commute. If $i=0$, then
\[
\begin{split}
\langle \mathscr{O}_\bullet, \mathscr{O}_{p,\bullet}^{(1)}(\mathbf{0}) \rangle_{par} = \deg(p) = 1, \quad  \langle \mathscr{O}_{p,\bullet}^{(1)}(\mathbf{0}) , \mathscr{O}_\bullet \rangle_{par} = 0, \quad \PHom(\mathscr{O}_\bullet, \mathscr{O}_{p,\bullet}^{(1)}(\mathbf{0}) = \mathbf{k}_p  \\
\PHom(\mathscr{O}_{p,\bullet}^{(1)}(\mathbf{0}),\mathscr{O}_\bullet) = \PExt(\mathscr{O}_\bullet, \mathscr{O}_{p,\bullet}^{(1)}(\mathbf{0} ) = \PExt(\mathscr{O}_{p,\bullet}^{(1)}(\mathbf{0}),\mathscr{O}_\bullet) = 0 
\end{split}.
\]
Thus
\[
I_{p,1}(0) 1_{\mathscr{O}_\bullet} = v^{-2}1_{\mathscr{L}_\bullet \oplus \mathscr{O}_{p,\bullet}^{(1)}(\mathbf{0})}, \quad 1_{\mathscr{O}_\bullet} I_{p,1}(0) = v^{-1} 1_{\mathscr{O}_\bullet \oplus \mathscr{O}_{p,\bullet}^{(1)}(\mathbf{0})}. 
\]
The same is true for any parabolic line bundle $\mathscr{L}_\bullet$ of multi-degree $\mathbf{0}$ and hence $I_{p,1}(0) \mathbf{1}^{vec}_{1,\mathbf{0}} = v^{-1} \mathbf{1}^{vec}_{1,\mathbf{0}} I_{p,1}(0)$. Similarly for $i=1$, we have
\[
\begin{split}
I_{p,1}(1) 1_{\mathscr{O}_\bullet} = v \big( 1_{\mathscr{O}_\bullet \oplus \mathscr{O}_{p,\bullet}^{(1)}(\epsilon_p)} + 1_{\mathscr{O}_\bullet(\epsilon_p)} \big), \quad 1_{\mathscr{O}_\bullet} I_{p,1}(1) = 1_{\mathscr{O}_\bullet \oplus \mathscr{O}_{p,\bullet}^{(1)}(\epsilon_p)}
\end{split}
\]
Therefore we have $I_{p,1}(1)\mathbf{1}^{vec}_{1,\mathbf{0}} = v\mathbf{1}^{vec}_{1,\mathbf{0}} I_{p,1}(1) + v \mathbf{1}^{vec}_{1,\alpha_{p,1}}$. Then \eqref{AppB:EQ1} holds for $j=1$. We will proceed the Lemma by induction on $j$. Let us assume \eqref{AppB:EQ1} holds for all $j'$ with $0 < j' < j$. For $i>j$, we have $i-j+1 > 1$ and hence $i > j-1$. By Lemma \ref{L:cyclicquiver} and induction hypothesis, we deduce
\[
\begin{split}
I_{p,j}(i) \mathbf{1}^{vec}_{1,\mathbf{0}} &= \big( v^{-1}I_{p,j-1}(i) I_{p,1}(i-j+1) - I_{p,1}(i-j+1) I_{p,j-1}(i) \big) \mathbf{1}^{vec}_{1,\mathbf{0}} \\
&= \mathbf{1}^{vec}_{1,\mathbf{0}} \big( v^{-1}I_{p,j-1}(i) I_{p,1}(i-j+1) - I_{p,1}(i-j+1) I_{p,j-1}(i) \big) = \mathbf{1}^{vec}_{1,\mathbf{0}} I_{p,j}(i).
\end{split}
\]
Similarly, if $i=0$, 
\[
\begin{split}
I_{p,j}(0) \mathbf{1}^{vec}_{1,\mathbf{0}} &= \big( v^{-1}I_{p,j-1}(0) I_{p,1}(1-j) - I_{p,1}(1-j) I_{p,j-1}(0) \big) \mathbf{1}^{vec}_{1,\mathbf{0}} \\
&=v^{-1} \mathbf{1}^{vec}_{1,\mathbf{0}} \big( v^{-1}I_{p,j-1}(0) I_{p,1}(1-j) - I_{p,1}(1-j) I_{p,j-1}(0) \big) = v^{-1} \mathbf{1}^{vec}_{1,\mathbf{0}} I_{p,j}(0).
\end{split}
\]
For $i=j$, on one hand we have
\[
\begin{split}
v^{-1} I_{p,j-1}(i) I_{p,1}(i-j+1) \mathbf{1}^{vec}_{1,\mathbf{0}} &= v^{-1}I_{p,j-1}(i)\big( v\mathbf{1}^{vec}_{1,\mathbf{0}} I_{p,1}(1) + v\mathbf{1}^{vec}_{1,\alpha_{p,1}} \big) \\
&=\mathbf{1}^{vec}_{1,\mathbf{0}} I_{p,j-1}(i)I_{p,1}(0) + \big( I_{p,j-1}(i-1) \mathbf{1}^{vec}_{1,\mathbf{0}} \big) \otimes \mathscr{O}_\bullet(\epsilon_p) \\
&= \mathbf{1}^{vec}_{1,\mathbf{0}} I_{p,j-1}(i) T_{p,1}(0) + v\mathbf{1}^{vec}_{1,\alpha_{p,1}} I_{p,j-1}(i) + v\mathbf{1}^{vec}_{1,\alpha_{p,1}+\cdots + \alpha_{p,i}}.
\end{split}
\]
On the other hand, 
\[
\begin{split}
I_{p,1}(i-j+1)I_{p,j-1}(i)\mathbf{1}^{vec}_{1,\mathbf{0}} &= I_{p,1}(i-j+1)\mathbf{1}^{vec}_{1,\mathbf{0}} I_{p,j-1}(i) \\
&= v\mathbf{1}^{vec}_{1,\mathbf{0}} I_{p,1}(1) I_{p,j-1}(0) + v\mathbf{1}^{vec}_{1,\alpha_{p,1}} I_{p,j-1}(i)
\end{split}
\]
Thus
\[
\begin{split}
I_{p,j}(i)\mathbf{1}^{vec}_{1,\mathbf{0}} &= v^{-1}I_{p,j-1}(i) I_{p,1}(i-j+1) \mathbf{1}^{vec}_{1,\mathbf{0}} - I_{p,1}(i-j+1) I_{p,j-1}(i) \mathbf{1}^{vec}_{1,\mathbf{0}} \\
&= \mathbf{1}^{vec}_{1,\mathbf{0}} \big( I_{p,j-1}(i) I_{p,1}(1) - vI_{p,1}(1) I_{p,j-1}(i) \big) + v\mathbf{1}^{vec}_{1,\alpha_{p,1}+\cdots + \alpha_{p,i}} \\
&= v\mathbf{1}^{vec}_{1,\mathbf{0}} T_{p,j}(i) + v \mathbf{1}^{vec}_{1,\alpha_{p,1} + \cdots + \alpha_{p,i}}.
\end{split}
\]
Now, if $0 < i = j-1$, then 
\[
\begin{split}
&v^{-1}I_{p,j-1}(i)I_{p,1}(0) \mathbf{1}^{vec}_{1,\mathbf{0}} = v^{-2}I_{p,j-1}(j-1)\mathbf{1}^{vec}_{1,\mathbf{0}}I_{p,1}(0) \\
&= v^{-1} \mathbf{1}^{vec}_{1,\mathbf{0}} I_{p,j-1}(j-1) + v^{-1}\mathbf{1}^{vec}_{1,\alpha_{p,1}+ \dots + \alpha_{p,j-1}} I_{p,1}(0).
\end{split}
\]
On the other hand,
\[
\begin{split}
I_{p,1}(0) I_{p,j-1}(j-1) \mathbf{1}^{vec}_{1,\mathbf{0}} &= vI_{p,1}(0)\mathbf{1}^{vec}_{1,\mathbf{0}}I_{p,j-1}(j-1) + vI_{p,1}(0)\mathbf{1}^{vec}_{1,\alpha_{p,1} + \cdots + \alpha_{p,j-1}} \\
&= \mathbf{1}^{vec}_{1,\mathbf{0}} I_{p,1}(0) I_{p,j-1}(j-1) + v\mathbf{1}^{vec}_{1,\alpha_{p,1} + \cdots + \alpha_{p,j-1}}I_{p,1}(0).
\end{split}
\]
Thus
\[
I_{p,j}(j-1) \mathbf{1}^{vec}_{1,\mathbf{0}} = \mathbf{1}^{vec}_{1,\mathbf{0}} I_{p,j}(j-1) + (v^{-1}-v)\mathbf{1}^{vec}_{1,\alpha_{p,1} + \cdots + \alpha_{p,j-1}} I_{p,1}(0).
\]
Finally, for $0 < i< j-1$ (hence $j-i>1$), we have
\[
\begin{split}
&v^{-1}I_{p,j-1}(i)I_{p,1}(i-j+1) \mathbf{1}^{vec}_{1,\mathbf{0}} = v^{-1}I_{p,j-1}(i) \mathbf{1}^{vec}_{1,\mathbf{0}}I_{p,1}(i-j+1) \\
&= v^{-1}\mathbf{1}^{vec}_{1,\mathbf{0}} I_{p,j-1}(i) I_{p,1}(i-j+1) + (v^{-1}-v)v^{-1}\mathbf{1}^{vec}_{1,\alpha_{p,1} + \cdots + \alpha_{p,i}} I_{p,j-i-1}(0) I_{p,1}(i-j+1).
\end{split}
\]
On the other hand, 
\[
\begin{split}
I_{p,1}(i-j+1) I_{p,j-1}(i) \mathbf{1}^{vec}_{1,\mathbf{0}} &= I_{p,1}(i-j+1) \mathbf{1}^{vec}_{1,\mathbf{0}}I_{p,j-1}(i) + (v^{-1}-v)I_{p,1}(i-j+1)\mathbf{1}^{vec}_{1,\alpha_{p,1} + \cdots + \alpha_{p,i}} T_{p,j-1}(0) \\
&= \mathbf{1}^{vec}_{1,\mathbf{0}} I_{p,1}(i-j+1) I_{p,j-1}(i) + (v^{-1}-v)\mathbf{1}^{vec}_{1,\alpha_{p,1} + \cdots + \alpha_{p,i}} I_{p,1}(i-j+1) I_{p,j-1}(0)
\end{split}.
\]
Thus, 
\[
I_{p,j}(i)\mathbf{1}^{vec}_{1,\mathbf{0}} = \mathbf{1}^{vec}_{1,\mathbf{0}} I_{p,j}(i) + (v^{-1}-v)\mathbf{1}^{vec}_{1,\alpha_{p,1}+ \cdots + \alpha_{p,i}} I_{p,j-i}(0).
\]
\qed

\small{

\vspace{.3in}

\noindent
Institute of Mathematics, Academia Sinica\\
6F, Astronomy-Mathematics Building, \\
No. 1, Sec. 4, Roosevelt Road, \\
Taipei 10617, TAIWAN \\
\texttt{jalin@math.sinica.edu.tw}
\end{document}